\documentclass[a4paper,11pt]{article}

\usepackage[utf8]{inputenc}
\usepackage[UKenglish]{babel}
\usepackage[T1]{fontenc}
\usepackage{amsmath,amsfonts,amssymb,amsthm}
\usepackage{graphicx}
\usepackage[top=2cm,bottom=3cm,left=2.5cm,right=2.5cm]{geometry}
\usepackage{setspace}
\usepackage{xcolor}
\usepackage{stmaryrd}
\usepackage{enumitem}
\usepackage{esint}
\usepackage{centernot}
\usepackage[format=hang,singlelinecheck=false]{caption}
\usepackage{framed}
\usepackage[normalem]{ulem}

\usepackage[colorlinks=true]{hyperref}

\title{Localisation/delocalisation of random height functions on (random) trees}
\author{Juhan Aru\thanks{École Polytechnique Fédérale de Lausanne, \href{mailto:juhan.aru@epfl.ch}{juhan.aru@epfl.ch}.} \and Guillaume Blanc\thanks{École Polytechnique Fédérale de Lausanne, \href{mailto:guillaume.blanc@epfl.ch}{guillaume.blanc@epfl.ch}, \url{https://sites.google.com/view/guillaume-blanc-math}.}}

\newif\iffigures
\figurestrue

\begin{document}

\theoremstyle{plain}
\newtheorem{thm}{Theorem}
\newtheorem*{thm*}{Theorem}
\newtheorem{lem}[thm]{Lemma}
\newtheorem{claim}[thm]{Claim}
\newtheorem{prop}[thm]{Proposition}
\newtheorem{cor}[thm]{Corollary}
\newtheorem*{lemtech}{Technical lemma}
\newtheorem{defn}[thm]{Definition}

\theoremstyle{remark}
\newtheorem{rem}[thm]{Remark}

\maketitle

\begin{abstract}
We consider integer-valued random height functions (gradient models) on trees pinned to zero on a random subset of leaves, and study the distribution of the height at the root of the tree.
More precisely, we consider the first $(n+1)$ generations of a supercritical Bienaymé tree, and select each leaf to have pinning independently with probability $p_n\in{]0,1]}$. We show that
a localisation/delocalisation phase transition occurs around $p_n=\left.n\middle/m^n\right.$, where $m\in{]1,\infty[}$ is the mean of the Bienaymé offspring distribution, and we establish a Gaussian scaling limit for the height at the root of the tree in the delocalised phase.
Our proof relies on a quenched concentration result for supercritical branching random walks.
\end{abstract}

\section*{Introduction and main results}

We study integer-valued random height functions on trees with pinning on a random subset of leaves.
One of our aims is to obtain a better understanding of the localisation/delocalisation phase transitions for height functions.
This has been studied recently for instance in \cite{peledsamotijyehudayoff,lammerstoninelli,logvariance,coquilledarioleny,hellerspinka}.
Another motivation is to probe the universality of the Gaussian scaling limit often appearing in the delocalised phase \cite{wu,chalhoubdariofaipeurleny}.

\quad

Let $T_n$ be a finite rooted tree with height $n\in\mathbb{N}^*=\{1,2,\ldots\}$ and with no internal leaves (i.e, such that each node $u\in T_n$ with $|u|\leq n-1$ has at least one child), and let $B_n$ be a non-empty subset of ${\partial T_n=\{u\in T_n:|u|=n\}}$.
For fixed $T_n$ and $B_n$, one can consider the following (gradient) model of integer-valued random height functions on $T_n$ with zero boundary conditions on $B_n$.

\paragraph{Integer-valued random height functions.}
Let $\mathbb{H}^0(T_n,B_n)$ denote the set of functions
\[h=(h_u)_{u\in T_n}\in\mathbb{Z}^{T_n}\]
such that $h_u=0$ for all $u\in B_n$, and fix a discrete probability measure $\nu$ on $\mathbb{Z}$.
We assume that $\nu$ is symmetric and has finite support.
We define a $\nu$-integer-valued random height function on $T_n$ with zero boundary conditions on $B_n$ as a random variable $H_n=(H_n(u))_{u\in T_n}$ with values in $\mathbb{H}^0(T_n,B_n)$, whose density with respect to the counting measure on $\mathbb{H}^0(T_n,B_n)$ is proportional to
\begin{equation}\label{eq:ivrhf}
(h_u)_{u\in T_n}\in\mathbb{H}^0(T_n,B_n)\longmapsto\prod_{(u,v)\in\vec{E}(T_n)}\nu(\{h_v-h_u\}),
\end{equation}
where $\vec{E}(T_n)=\left\{(u,v)\in T_n^2:\text{$v$ is a child of $u$}\right\}$ denotes the set of directed edges of $T_n$.
Note that the strict positivity of the partition function
\[Z=\sum_{h\in\mathbb{H}^0(T_n,B_n)}\prod_{(u,v)\in\vec{E}(T_n)}\nu(\{h_v-h_u\})\]
is ensured by the fact that $|u|=n$ for all $u\in B_n$, and its finiteness follows from the fact that $B_n$ is non-empty.\footnote{Observe that for any $w\in B_n$, we have
\[Z\leq\sum_{h\in\mathbb{H}^0(T_n,\{w\})}\prod_{(u,v)\in\vec{E}(T_n)}\nu(\{h_v-h_u\})=1.\]}
The symmetry assumption on $\nu$ ensures that the random variables $(H_n(u))_{u\in T_n}$ are centred.

\quad

\begin{figure}[h]
\centering
\iffigures\includegraphics[width=0.74\linewidth]{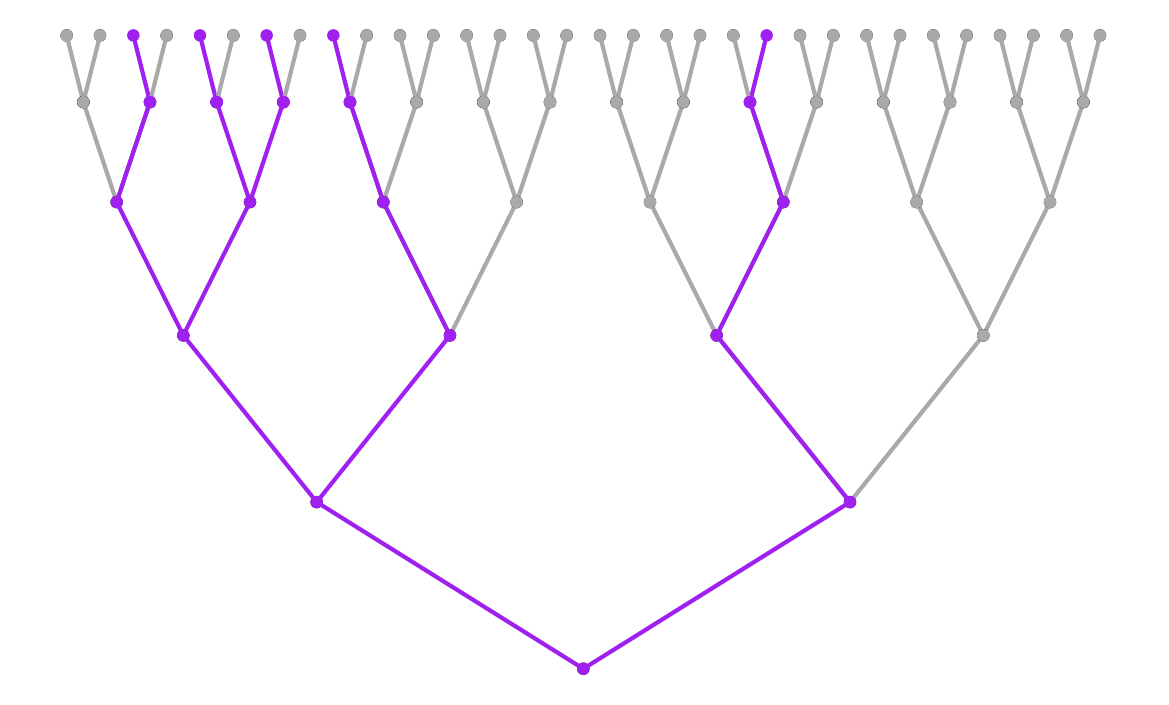}\fi
\caption{In gray, the tree $T_n$ consisting of the first $(n+1)$ generations of the binary tree, with ${n=5}$.
The boundary $B_n$ consists of the purple leaves.
The ``pruned'' tree $T_n'$ induced by $B_n$, in purple, consists of the nodes of $T_n$ that contain a purple leaf in their descendance.
For the models of random height functions considered in this paper, sampling a random height function on $T_n$ with zero boundary conditions on $B_n$ amounts to sampling a random height function on the pruned tree $T_n'$ with zero boundary conditions on $\partial T_n'=B_n$ (in purple), and sampling independent and identically distributed random variables on the remaining edges of $T_n$ (in gray).}\label{fig:tree}
\end{figure}

\paragraph{Supercritical Bienaymé trees with sparse boundary conditions.}
The set of (random) trees considered in this paper originates from supercritical Bienaymé trees, on which we work conditionally.
The (random) subset of leaves for pinning is then chosen using independent coin tosses. Let us make this more precise.
Fix a discrete probability measure $\mu$ on $\mathbb{N}=\{0,1,\ldots\}$, and let $\mathcal{T}$ be a Bienaymé tree with offspring distribution $\mu$.
We assume that $\mu(\{0\})=0$ and $\mu(\{1\})<1$ (in particular, the tree $\mathcal{T}$ is infinite and has no internal leaves).
Moreover, we assume that $\sum_{k\geq0}k^2\cdot\mu(\{k\})<\infty$, and we denote by ${m=\sum_{k\geq0}k\cdot\mu(k)\in{]1,\infty[}}$ the mean of $\mu$.
Next, for each $n\in\mathbb{N}^*$, fix $p_n\in{]0,1]}$, and conditionally on $\mathcal{T}$, let $(X_u)_{u\in\partial[\mathcal{T}]_n}$ be a collection of independent Bernoulli random variables with parameter $p_n$, where $\partial[\mathcal{T}]_n=\{u\in\mathcal{T}:|u|=n\}$ denotes the $n$-th generation of $\mathcal{T}$.\footnote{Formally, let $(X_u)_{u\in\partial[\mathbb{U}]_n}$ be a collection of independent Bernoulli random variables with parameter $p_n$, independent of $\mathcal{T}$, where $\partial[\mathbb{U}]_n=\{u\in\mathbb{U}:|u|=n\}$ denotes the $n$-th generation of the Ulam--Harris tree $\mathbb{U}$.}
Then, set ${\mathcal{B}_n=\{u\in\partial[\mathcal{T}]_n:X_u=1\}}$.
Conditionally on the event
\begin{equation}\label{eq:boundary}
(\mathcal{B}_n\neq\emptyset)=(\text{there exists $u\in\partial[\mathcal{T}]_n$ such that $X_u=1$}),
\end{equation}
we want to consider the distribution of the height at the root of the tree for our model of random height functions on $[\mathcal{T}]_n$ with zero boundary conditions on $\mathcal{B}_n$.
A first observation is that this distribution is the same as for random height functions on $\mathcal{T}_n'$ with zero boundary conditions on $\partial\mathcal{T}_n'=\mathcal{B}_n$, where $\mathcal{T}_n'$ is the ``pruned'' tree:
\[\mathcal{T}_n'=\{u\in[\mathcal{T}]_n:\text{there exists $v\in\mathcal{B}_n$ such that $v\succeq u$}\}.\]
See Proposition \ref{prop:prune} below for a more complete statement, and see Figure \ref{fig:tree} for an illustration.

We denote by $\mathcal{T}_n^{p_n}$ the random tree $\mathcal{T}_n'$ conditionally on the event \eqref{eq:boundary}.
In view of the previous observation, we will work from now on conditionally on the random tree $\mathcal{T}_n^{p_n}$.

\paragraph{Main result for integer-valued random height functions.}
To state our main result, conditionally on the random tree $\mathcal{T}_n^{p_n}$, let ${H_n=(H_n(u))_{u\in\mathcal{T}_n^{p_n}}}$ be a $\nu$-integer-valued random height function on $\mathcal{T}_n^{p_n}$ with zero boundary conditions on $\partial\mathcal{T}_n^{p_n}$.

\begin{thm}\label{thm:main}
The following localisation/delocalisation phase transition occurs:
\begin{itemize}
\item \emph{Delocalisation:} If $(m^n\cdot p_n)/n\rightarrow0$ as $n\to\infty$, then $\left.\#\partial\mathcal{T}_n^{p_n}\middle/n\right.\rightarrow0$ in probability as $n\to\infty$, and we have
\begin{equation}\label{eq:clt}
\mathbb{P}\left(\frac{\sqrt{\#\partial\mathcal{T}_n^{p_n}}}{\sqrt{n}}\cdot H_n(\varnothing)\in\cdot\,\middle|\,\mathcal{T}_n^{p_n}\right)\longrightarrow\mathcal{N}\left(0,\sigma^2\right)\quad\text{in probability as $n\to\infty$}
\end{equation}
in the Polish space of Borel probability measures on $\mathbb{R}$ endowed with the topology of weak convergence, where $\sigma^2=\sum_{x\in\mathbb{Z}}x^2\cdot\nu(\{x\})$ is the variance of $\nu$.

\item \emph{Localisation:} Suppose that $\nu$ is the uniform distribution on $\{-1,1\}$.
If
\[\inf_{n\geq1}\frac{m^n\cdot p_n}{n}>0,\]
then the sequence of random variables ${\left(\mathbb{V}\left(H_n(\varnothing)\,\middle|\,\mathcal{T}_n^{p_n}\right)\right)_{n\in\mathbb{N}^*}}$ is tight.
\end{itemize}
\end{thm}

The reason why the localisation in Theorem \ref{thm:main} occurs under the more restrictive assumption that $\nu$ is the uniform distribution on $\{-1,1\}$ is that we use a monotonicity property for the absolute value of $\nu$-integer-valued random height functions, which is not available for general distributions $\nu$ (see Subsection \ref{subsec:locivrhf} for more details).

\begin{rem}
At the very final stages of our work, we learned about the work of Heller and Spinka \cite{hellerspinka}, whose results (combined with our Theorem \ref{thm:dgff} below) yield a different proof of the localisation/delocalisation phase transition in Theorem \ref{thm:main}, in the case where $\nu$ is the uniform distribution on $\{-1,1\}$, or the uniform distribution on $\llbracket-L,L\rrbracket$ for some $L\in\mathbb{N}^*$.
They give up-to-constant comparisons between the variance of the height at the root of the tree for the integer-valued model and the effective resistance between the root and the leaves of the tree (which in fact corresponds to the variance of the height at the root of the tree for the discrete Gaussian free field, see below), but they do not consider the Gaussian scaling limit \eqref{eq:clt}, which is our main result.
\end{rem}

\paragraph{The Gaussian case as a toy model.}
The location of the phase transition can be predicted using the following as a toy model, namely the discrete Gaussian free field.
More precisely, for fixed $T_n$ and $B_n$ as above, let $\mathbb{G}^0(T_n,B_n)$ denote the set of functions ${g=(g_u)_{u\in T_n}\in\mathbb{R}^{T_n}}$ such that $g_u=0$ for all $u\in B_n$.
The discrete Gaussian free field on $T_n$ with zero boundary conditions on $B_n$ can be defined as a random variable ${G_n=(G_n(u))_{u\in T_n}}$ with values in $\mathbb{G}^0(T_n,B_n)$, whose density with respect to the Lebesgue measure on $\mathbb{G}^0(T_n,B_n)$ is proportional to
\begin{equation}\label{eq:dgff}
(g_u)_{u\in T_n}\in\mathbb{G}^0(T_n,B_n)\longmapsto\prod_{(u,v)\in\vec{E}(T_n)}\exp\left(-\frac{(g_v-g_u)^2}{2}\right).
\end{equation}
In this model, the random vector $(G_n(u))_{u\in T_n\setminus B_n}$ is a non-degenerate centred Gaussian vector, so that understanding the distribution of the height at the root of the tree $G_n(\varnothing)$ boils down to estimating its variance $\mathbb{V}(G_n(\varnothing))$, which can be calculated (in principle) by solving a recursion on the tree (see Section \ref{sec:dgff} for more details).

\quad

Consider now again the same random tree $\mathcal{T}_n^{p_n}$ as above, and conditionally on $\mathcal{T}_n^{p_n}$, let $G_n=(G_n(u))_{u\in\mathcal{T}_n^{p_n}}$ be the discrete Gaussian free field on $\mathcal{T}_n^{p_n}$ with zero boundary conditions on $\partial\mathcal{T}_n^{p_n}$.

\begin{thm}\label{thm:dgff}
The following localisation/delocalisation phase transition occurs:
\begin{itemize}

\item \emph{Delocalisation:} If $(m^n\cdot p_n)/n\rightarrow0$ as $n\to\infty$, then $\left.\#\partial\mathcal{T}_n^{p_n}\middle/n\right.\rightarrow0$ in probability as $n\to\infty$, and we have
\[\frac{\#\partial\mathcal{T}_n^{p_n}}{n}\cdot\mathbb{V}\left(G_n(\varnothing)\,\middle|\,\mathcal{T}_n^{p_n}\right)\longrightarrow1\quad\text{in probability as $n\to\infty$.}\]
In particular, we get that
\[\mathbb{P}\left(\frac{\sqrt{\#\partial\mathcal{T}_n^{p_n}}}{\sqrt{n}}\cdot G_n(\varnothing)\in\cdot\,\middle|\,\mathcal{T}_n^{p_n}\right)\longrightarrow\mathcal{N}(0,1)\quad\text{in probability as $n\to\infty$}\]
in the Polish space of Borel probability measures on $\mathbb{R}$ endowed with the topology of weak convergence.
\item \emph{Localisation:} If $\inf_{n\geq1}(m^n\cdot p_n)/n>0$, then the sequence of random variables
\[\left(\mathbb{V}\left(G_n(\varnothing)\,\middle|\,\mathcal{T}_n^{p_n}\right)\right)_{n\in\mathbb{N}^*}\]
is tight.
\end{itemize}
\end{thm}


\paragraph{The starting point for proofs.}
The random tree $\mathcal{T}_n^{p_n}$ turns out to have the same distribution as the genealogical tree associated with an explicit inhomogenous branching process (see Proposition \ref{prop:pruned} below).
This key observation has been made by Ayuso Ventura and Berger \cite[Proposition 5.1]{irenequentinising} in their study of the Ising model on supercritical Bienaymé trees with a sparse random external field.
As mentioned earlier, the variance of the height at the root of the tree for the discrete Gaussian free field on $\mathcal{T}_n^{p_n}$ with zero boundary conditions on $\partial\mathcal{T}_n^{p_n}$ can be calculated (in principle) by solving a recursion on the tree.
Auyso Ventura and Berger \cite{irenequentinising,irenequentinnonlinear} have shown in great generality how sharp quenched estimates can be derived from such recursions on branching random trees (see also \cite{chenhulin} by Chen, Hu and Lin), and the proof of Theorem \ref{thm:dgff} is another instance of the power of such branching processes arguments.

The proof of our main result for integer-valued random height functions (Theorem \ref{thm:main}) proceeds differently, and a key ingredient is the following quenched concentration result for supercritical branching random walks (see Theorem \ref{thm:concentration2} below).

\paragraph{Quenched concentration results for supercritical branching random walks.}
For each $n\in\mathbb{N}^*$, fix a finite rooted tree $T_n$ with height $n$ and with no internal leaves.
Then, let $(X_u)_{u\in T_n^*}$ be a collection of independent random variables with distribution $\nu$, where ${T_n^*=T_n\setminus\{\varnothing\}}$, and set
\[S_v=\sum_{\substack{u\in T_n^*\\u\preceq v}}X_u\quad\text{for all $v\in T_n$.}\]
Here, we do not assume that $\nu$ is discrete: we assume that $\nu$ is a Borel probability measure on $\mathbb{R}$, which is symmetric and has compact support.
Let ${\sigma^2=\int_{\mathbb{R}}x^2\cdot\nu(\mathrm{d}x)}$ be the variance of $\nu$, and set $L=\inf\{\ell\in\mathbb{R}_+:\nu([-\ell,\ell])=1\}$.
We derive concentration bounds for the following random variables:
\[\langle S\rangle_n:=\frac{1}{N_n}\cdot\sum_{v\in\partial T_n}S_v\quad\text{and}\quad\left\langle S^2\right\rangle_n:=\frac{1}{N_n}\cdot\sum_{v\in\partial T_n}S_v^2,\]
where $N_n=\#\partial T_n$.
Note that $\mathbb{E}[\langle S\rangle_n]=0$ and $\mathbb{E}\left[\left\langle S^2\right\rangle_n\right]=\sigma^2\cdot n$.
Our bounds involve the following tree-dependent quantity:
\begin{equation}\label{eq:gamman}
\gamma_n:=1\vee\left(\frac{1}{N_n}\cdot\sum_{k=1}^nk^3\cdot\max_{\substack{u\in T_n\\|u|=k}}N_n(u)\right),
\end{equation}
where $N_n(u)=\#\{v\in\partial T_n:v\succeq u\}$ for all $u\in T_n^*$.
A direct calculation shows that $\langle S\rangle_n$ is $\left(L^2\gamma_n\right)$-sub-Gaussian (see Proposition \ref{prop:concentration1} below), but the concentration of $\left\langle S^2\right\rangle_n$ is more delicate.
We obtain the following result:

\begin{thm}\label{thm:concentration2}
We have
\begin{equation}\label{eq:concentration2}
\mathbb{E}\left[\exp\left(\lambda\left(\left\langle S^2\right\rangle_n-\sigma^2n\right)\right)\right]\leq\exp\left(\frac{\lambda^2}{\delta_n^2}\right)\quad\text{for all $\lambda\in\left[-\delta_n,\delta_n\right]$,}
\end{equation}
where
\[\delta_n=\frac{1}{32L^2\sqrt{2\gamma_n}}.\]
\end{thm}

Thus, if $\gamma_n$ is bounded, then $\delta_n$ is bounded away from zero, and \eqref{eq:concentration2} entails that the random variables $\left\langle S^2\right\rangle_n$ are concentrated around their mean $\sigma^2n$ with uniform sub-exponential tails.
Note that if $T_n$ consists of the first $(n+1)$ generations of the $m$-ary tree, where $m\geq2$, then
\[\gamma_n=1\vee\left(\frac{1}{m^n}\cdot\sum_{k=1}^nk^3\cdot m^{n-k}\right)\leq1\vee\left(\sum_{k\geq1}\frac{k^3}{m^k}\right).\]
Furthermore, for our supercritical Bienaymé tree $\mathcal{T}$ with offspring distribution $\mu$, the sequence of random variables
\[\left(\Gamma_n:=\frac{1}{\#\partial[\mathcal{T}]_n}\cdot\sum_{k=1}^nk^3\cdot\max_{\substack{u\in\mathcal{T}\\|u|=k}}\#\{v\in\partial[\mathcal{T}]_n:v\succeq u\}\right)_{n\in\mathbb{N}^*}\]
is easily seen to be tight (see Remark \ref{rem:cutheightbienaymé} below).
In this respect, Theorem \ref{thm:concentration2} is a quenched concentration result for supercritical branching random walks.

\paragraph{Some context.}
To finish this introduction, let us put our results into context by discussing some examples which sparked the current article.
In particular, the model of supercritical Bienaymé trees with sparse boundary conditions considered in this paper ``interpolates'' in a natural way between the following two examples of regular and starlike trees.

\subparagraph{Regular trees.}
Let $T_n$ consist of the first $(n+1)$ generations of the $m$-ary tree, where $m\geq2$ (i.e, every node $u\in T_n$ with $|u|\leq n-1$ has $m$ children).
Random height functions on $T_n$ with zero boundary conditions on $B_n=\partial T_n$ are localised, in the following sense:
\begin{itemize}
\item For the discrete Gaussian free field $G_n=(G_n(u))_{u\in T_n}$ on $T_n$ with zero boundary conditions on $\partial T_n$, the variance of the height at the root of the tree can be calculated to be (see Remark \ref{rem:vardgfftrees} below)
\[\mathbb{V}(G_n(\varnothing))=\frac{1-m^{-n}}{m-1}\underset{n\to\infty}{\longrightarrow}\frac{1}{m-1}.\]

\item For the $\nu$-integer-valued random height function $H_n=(H_n(u))_{u\in T_n}$ on $T_n$ with zero boundary conditions on $\partial T_n$, it has been shown by Peled, Samotij and Yehudayoff \cite[Theorem 1.2]{peledsamotijyehudayoff} that in the case where $\nu$ is the uniform distribution on $\llbracket-L,L\rrbracket$ for some $L\in\mathbb{N}^*$, there exists a constant ${c=c(m,L)\in\mathbb{R}_+^*}$ such that for each $n\in\mathbb{N}^*$, we have
\[\mathbb{P}\left(H_n(\varnothing)=y\right)\leq\exp\left(-c\cdot m^{|y|/L}\right)\quad\text{for all $y\in\mathbb{Z}$.}\]
Moreover, it follows from a result of Lammers and Toninelli \cite[Lemma 3.7]{lammerstoninelli} that in the case where $\nu$ is the uniform distribution on $\{-1,1\}$, there exists an absolute constant $C\in\mathbb{R}_+^*$ such that $\mathbb{V}(H_n(\varnothing))\leq C$ for all $n\in\mathbb{N}^*$.\footnote{To be perfectly accurate, their result \cite[Lemma 3.7]{lammerstoninelli} applies to trees with minimal degree $3$, which is not the case of $T_n$ when $m=2$ due to the single root node.}
Let us also mention that recently, Butler, Krishnan, Ray and Spinka \cite[Theorem 1.1]{localconv} have studied the convergence in distribution of $(H_n(\varnothing))_{n\in\mathbb{N}^*}$, in the case where $\nu$ is the uniform distribution on $\llbracket-1,1\rrbracket$: it turns out that whether convergence occurs or not depends on $m$.
\end{itemize}

\subparagraph{Starlike trees.}
Let $T_n$ be the rooted tree in which the root has $N_n\in\mathbb{N}^*$ children, and every other node $u\in T_n$ with $|u|\leq n-1$ has one child.
Random height functions on $T_n$ with zero boundary conditions on $B_n=\partial T_n$ undergo the following localisation/delocalisation phase transition:
\begin{itemize}
\item For the discrete Gaussian free field $G_n=(G_n(u))_{u\in T_n}$ on $T_n$ with zero boundary conditions on $\partial T_n$, the variance of the height at the root of the tree can be calculated to be  (see Remark \ref{rem:vardgfftrees} below)
\[\mathbb{V}(G_n(\varnothing))=\frac{n}{N_n}.\]
Therefore, the height at the root of the tree is delocalised if $N_n/n\rightarrow0$ as $n\to\infty$, and localised if $\inf_{n\geq1}(N_n/n)>0$.

\item For the $\nu$-integer-valued random height function $H_n=(H_n(u))_{u\in T_n}$ on $T_n$ with zero boundary conditions on $\partial T_n$, it has been shown by Benjamini, Häggström and Mossel \cite[Proposition 3.1 and Proposition 3.2]{randomgraphhom} that in the case where $\nu$ is the uniform distribution on $\{-1,1\}$, the following holds:
\begin{itemize}
\item If $N_n/n\rightarrow0$ as $n\to\infty$, then
\[\frac{\sqrt{N_n}}{\sqrt{n}}\cdot H_n(\varnothing)\longrightarrow\mathcal{N}(0,1)\quad\text{in distribution as $n\to\infty$.}\]

\item If $\inf_{n\geq1}(N_n/n)>0$, then the sequence of random variables $(H_n(\varnothing))_{n\in\mathbb{N}^*}$ is tight.
Moreover, it has been shown that if $N_n/n\rightarrow\infty$ as $n\to\infty$, then we have
\[\mathbb{P}(H_{2n}(\varnothing)=0)\longrightarrow1\quad\text{and}\quad\mathbb{P}(H_{2n+1}(\varnothing)=\pm1)\longrightarrow\frac{1}{2}\quad\text{as $n\to\infty$.}\]
See \cite[Proposition 3.3]{randomgraphhom}.

In fact, in this case one can even go further and observe that when $N_n/n\rightarrow\rho\in\mathbb{R}_+^*$ as $n\to\infty$, the height at the root of the tree converges in distribution towards a discrete Gaussian random variable!\footnote{More precisley, for all $k\in\mathbb{Z}$, we have
\[\mathbb{P}(H_{2n}(\varnothing)=2k)\longrightarrow\frac{\exp\left(-2\rho k^2\right)}{\sum_{\ell=-\infty}^\infty\exp\left(-2\rho\ell^2\right)}\quad\text{and}\quad\mathbb{P}(H_{2n+1}(\varnothing)=2k+1)\longrightarrow\frac{\exp\left(-2\rho k(k+1)\right)}{\sum_{\ell=-\infty}^\infty\exp\left(-2\rho\ell(\ell+1)\right)}\]
as $n\to\infty$.}
Understanding the generality of such critical behaviour would of course be of interest too.
\end{itemize}
\end{itemize}

\paragraph{Acknowledgements.}
G.B.~is extremely grateful to Quentin Berger for enlightening and encouraging discussions around the recursion on the pruned tree, to Pascal Massart for his prompt advice on the concentration problem, and to Sébastien Ott for numerous helpful discussions.

\paragraph{AI use statement.}
The mathematical content of this article was complete before summer 2026, and we have not used generative AI tools in the preparation of the paper.

\tableofcontents

\section{Preliminary estimates on the pruned tree}\label{sec:pruned}

In this section, we collect some preliminary estimates on the ``pruned tree'' $\mathcal{T}_n^p$, whose definition we recall below.
First, the following proposition justifies the fact that for our models of random height functions, it suffices to study the distribution of the height at the root of this pruned tree.

\begin{prop}\label{prop:prune}
Let $T_n$ be a finite rooted tree with height $n\in\mathbb{N}^*$ and with no internal leaves, and let $B_n\subset\partial T_n$ be a non-empty subset of leaves.
Then, set
\[T_n'=\{u\in T_n:\text{there exists $v\in B_n$ such that $v\succeq u$}\}.\]
See Figure \ref{fig:tree} for an illustration.
We have the following identities in distribution:
\begin{itemize}
\item Let $G_n'=(G_n'(u))_{u\in T_n'}$ be the discrete Gaussian free field on $T_n'$ with zero boundary conditions on ${\partial T_n'=B_n}$, and let $(X_{u,v})_{(u,v)\in\left.\vec{E}(T_n)\middle\backslash\vec{E}(T_n')\right.}$ be a collection of independent random variables with standard Gaussian distribution, independent of $G_n'$.
Then, set
\[G_n(w)=G_n'(\varnothing)+\sum_{\substack{(u,v)\in\vec{E}(T_n)\\v\preceq w}}Y_{u,v}\quad\text{for all $w\in T_n$,}\]
where
\[Y_{u,v}=\begin{cases}
G_n'(v)-G_n'(u)&\text{if $(u,v)\in\vec{E}(T_n')$,}\\
X_{u,v}&\text{otherwise,}
\end{cases}\quad\text{for all $(u,v)\in\vec{E}(T_n)$.}\]
The random height function $G_n=(G_n(u))_{u\in T_n}\in\mathbb{G}^0(T_n,B_n)$ has the same distribution as the discrete Gaussian free field on $T_n$ with zero boundary conditions on $B_n$.

\item Let $H_n'=(H_n'(u))_{u\in T_n'}$ be a $\nu$-integer-valued random height function on $T_n'$ with zero boundary conditions on ${\partial T_n'=B_n}$, and let $(X_{u,v})_{(u,v)\in\left.\vec{E}(T_n)\middle\backslash\vec{E}(T_n')\right.}$ be a collection of independent random variables with distribution $\nu$, independent of $H_n'$.
Then, set
\[H_n(w)=H_n'(\varnothing)+\sum_{\substack{(u,v)\in\vec{E}(T_n)\\v\preceq w}}Y_{u,v}\quad\text{for all $w\in T_n$,}\]
where
\[Y_{u,v}=\begin{cases}
H_n'(v)-H_n'(u)&\text{if $(u,v)\in\vec{E}(T_n')$,}\\
X_{u,v}&\text{otherwise,}
\end{cases}\quad\text{for all $(u,v)\in\vec{E}(T_n)$.}\]
The random height function $H_n=(H_n(u))_{u\in T_n}\in\mathbb{H}^0(T_n,B_n)$ has the same distribution as a $\nu$-integer-valued random height function on $T_n$ with zero boundary conditions on $B_n$.
\end{itemize}
\end{prop}
\begin{proof}
\begin{itemize}
\item By \eqref{eq:dgff}, the random vector $G':=(G_n'(u))_{u\in T_n'\setminus B_n}$ has a density $\varphi$ with respect to the Lebesgue measure on $\mathbb{R}^{T_n'\setminus B_n}$ which is given by
\[\varphi(g')=\frac{1}{Z}\cdot\prod_{(u,v)\in\vec{E}(T_n')}\exp\left(-\frac{(g'_v-g'_u)^2}{2}\right)\quad\text{for all $g'=(g'_u)_{u\in T_n'\setminus B_n}\in\mathbb{R}^{T_n'\setminus B_n}$,}\]
where $g'_v=0$ for all $v\in B_n$, and where
\[Z=\int_{\mathbb{R}^{T_n'\setminus B_n}}\prod_{(u,v)\in\vec{E}(T_n')}\exp\left(-\frac{(g'_v-g'_u)^2}{2}\right)\mathrm{d}g'.\]
Since the random vector $X:=(X_{u,v})_{(u,v)\in\left.\vec{E}(T_n)\middle\backslash\vec{E}(T_n')\right.}$ is independent of $G'$, and since its coordinates are independent and have standard Gaussian distribution, we get that the random vector $(G',X)$ has a density $\chi$ with respect to the Lebesgue measure on ${\mathbb{R}^{T_n'\setminus B_n}\times\mathbb{R}^{\left.\vec{E}(T_n)\middle\backslash\vec{E}(T_n')\right.}}$ which is given by
\[\chi(g',x)=\varphi(g')\cdot\prod_{(u,v)\in\left.\vec{E}(T_n)\middle\backslash\vec{E}(T_n')\right.}\left(\frac{1}{\sqrt{2\pi}}\cdot\exp\left(-\frac{x_{u,v}^2}{2}\right)\right)\]
for all $g'=(g'_u)_{u\in T_n'}\in\mathbb{R}^{T_n'\setminus B_n}$ and $x=(x_{u,v})_{(u,v)\in\left.\vec{E}(T_n)\middle\backslash\vec{E}(T_n')\right.}\in\mathbb{R}^{\left.\vec{E}(T_n)\middle\backslash\vec{E}(T_n')\right.}$.
It follows that the random vector $G:=(G_n(u))_{u\in T_n\setminus B_n}$ has a density $\psi$ with respect to the Lebesgue measure on $\mathbb{R}^{T_n\setminus B_n}$ which is given by
\begin{eqnarray*}
\lefteqn{\psi(g)}\\
&=&\frac{1}{Z}\cdot\prod_{(u,v)\in\vec{E}(T_n')}\exp\left(-\frac{(g_v-g_u)^2)}{2}\right)\cdot\prod_{(u,v)\in\left.\vec{E}(T_n)\middle\backslash\vec{E}(T_n')\right.}\left(\frac{1}{\sqrt{2\pi}}\cdot\exp\left(-\frac{(g_v-g_u)^2}{2}\right)\right)\\
&=&\frac{1}{Z}\cdot\left(\frac{1}{\sqrt{2\pi}}\right)^{\#\left(\vec{E}(T_n)\middle\backslash\vec{E}(T_n')\right)}\cdot\prod_{(u,v)\in\vec{E}(T_n)}\exp\left(-\frac{(g_v-g_u)^2}{2}\right)
\end{eqnarray*}
for all $g=(g_u)_{u\in T_n\setminus B_n}\in\mathbb{R}^{T_n\setminus B_n}$.
To justify the previous change-of-variables properly, observe that $G=\Phi(G',X)$, where
\[\begin{matrix}
 \Phi:&\mathbb{R}^{T_n'\setminus B_n}\times\mathbb{R}^{\left.\vec{E}(T_n)\middle\backslash\vec{E}(T_n')\right.}&\longrightarrow&\mathbb{R}^{T_n\setminus B_n}\\
&(g',x)&\longmapsto&g
\end{matrix}\]
is given by
\[g_w=g'_\varnothing+\sum_{\substack{(u,v)\in\vec{E}(T_n)\\v\preceq w}}y_{u,v}\quad\text{for all $w\in T_n$,}\]
with
\[y_{u,v}=\begin{cases}
g'_v-g'_u&\text{if $(u,v)\in\vec{E}(T_n')$,}\\
x_{u,v}&\text{otherwise,}
\end{cases}\quad\text{for all $(u,v)\in\vec{E}(T_n)$.}\]
Notice that $\Phi$ is an invertible linear map with determinant $1$.
One way to see this is to enumerate the nodes of $T_n$ in depth-first order,\footnote{Breadth-first order would also work.} and consider the natural orders induced on $T_n\setminus B_n$ and on ${(T_n'\setminus B_n)\sqcup\left(\vec{E}(T_n)\middle\backslash\vec{E}(T_n')\right)}$, respectively.
With respect to these orders, the matrix associated with $\Phi$ is upper-triangular and has diagonal entries $1$.
See Figure \ref{fig:prop1}.
This shows that $\det\Phi=1$, which justifies the above change-of-variables.
Finally, by \eqref{eq:dgff}, we conclude that the random height function $(G_n(u))_{u\in T_n}$ has the same distribution as the discrete Gaussian free field on $T_n$ with zero boundary conditions on $B_n$, as desired.\footnote{Moreover, we get that
\[\int_{\mathbb{R}^{T_n\setminus B_n}}\prod_{(u,v)\in\vec{E}(T_n)}\exp\left(-\frac{(g_v-g_u)^2}{2}\right)\mathrm{d}g=Z\cdot\left(\sqrt{2\pi}\right)^{\left(\vec{E}(T_n)\middle\backslash\vec{E}(T_n')\right)}.\]}

\begin{figure}[h]
\centering
\begin{tabular}{lcr}
\iffigures\includegraphics[width=0.3\linewidth]{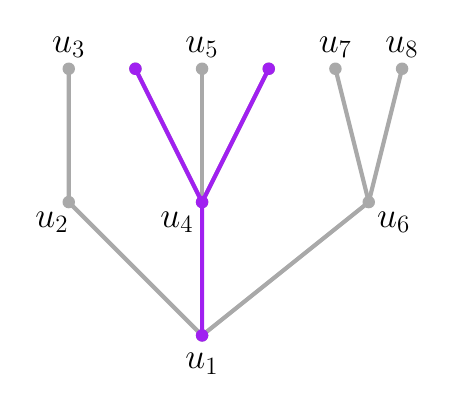}&&\includegraphics[width=0.3\linewidth]{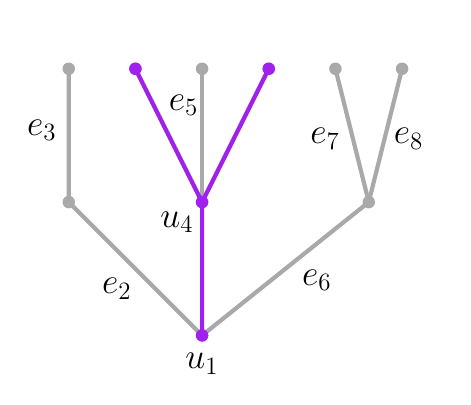}\fi
\end{tabular}
\caption{Illustration of the proof of Propostion \ref{prop:prune}.
In gray, the tree $T_n$.
In purple, the pruned tree $T_n'$.
The depth-first order on the nodes of $T_n$ induces a natural order on $T_n\setminus B_n$ (left) and on $(T_n'\setminus B_n)\sqcup\left(\vec{E}(T_n)\middle\backslash\vec{E}(T_n')\right)$ (right).
With respect to these orders, the matrix of $\Phi$ is upper-triangular and has diagonal entries $1$.
Indeed, we have
\begin{eqnarray*}
g_{u_1}&=&g'_{u_1},\\
g_{u_2}&=&g'_{u_1}+x_{e_2},\\
g_{u_3}&=&g'_{u_1}+x_{e_2}+x_{e_3},\\
g_{u_4}&=&g'_{u_4},\\
g_{u_5}&=&g'_{u_4}+x_{e_5},\\
g_{u_6}&=&g'_{u_1}+x_{e_6},\\
g_{u_7}&=&g'_{u_1}+x_{e_6}+x_{e_7},\\
g_{u_8}&=&g'_{u_4}+x_{e_6}+x_{e_8}.
\end{eqnarray*}}\label{fig:prop1}
\end{figure}

\item We omit the proof, which is based on the same ideas as in the previous point.
\end{itemize}
\end{proof}

Let us now dive into the study of the pruned tree $\mathcal{T}_n^p$.
Recall that $\mathcal{T}$ is a Bienaymé tree with offspring distribution $\mu$.
We assume that ${\mu(\{0\})=0}$ and ${\mu(\{1\})<1}$.
Moreover, we assume that $\sum_{k\geq0}k^2\cdot\mu(\{k\})<\infty$, and we denote by ${m=\sum_{k\geq0}k\cdot\mu(k)\in{]1,\infty[}}$ the mean of $\mu$.
Set ${Z_n=\#\partial[\mathcal{T}]_n}$ for all $n\in\mathbb{N}$.
As $n\to\infty$, we have
\begin{equation}\label{eq:W}
\frac{Z_n}{m^n}\longrightarrow W\quad\text{almost surely and in $L^2$.}
\end{equation}
In particular, there exists a constant $C\in\mathbb{R}_+^*$ such that
\begin{equation}\label{eq:Zn2}
\mathbb{E}\left[Z_n^2\right]\leq C\cdot m^{2n}\quad\text{for all $n\in\mathbb{N}$.}
\end{equation}
Moreover, since $\mu(\{0\})=0$, we have $W>0$ almost surely.\footnote{Recall that $W>0$ almost surely on the non-extinction event.}
It follows that the sequence $\left(Z_n\middle/m^n\right)_{n\in\mathbb{N}}$ is tight in $\mathbb{R}_+^*$: for each $\varepsilon\in{]0,1[}$, there exists $\delta<\Delta\in\mathbb{R}_+^*$ such that
\begin{equation}\label{eq:Znmntight}
\mathbb{P}\left(\delta\leq\frac{Z_n}{m^n}\leq\Delta\right)\leq\varepsilon\quad\text{for all $n\in\mathbb{N}$.}
\end{equation}

Now, fix $n\in\mathbb{N}^*$ and $p\in{]0,1]}$, and let $(X_u)_{u\in\partial[\mathbb{U}]_n}$ be a collection of independent Bernoulli random variables with parameter $p$, independent of $\mathcal{T}$, where $\mathbb{U}$ is the Ulam--Harris tree.\footnote{The Ulam--Harris tree is the infinite plane tree on the set of words on the alphabet $\mathbb{N}^*=\{1,2,\ldots\}$, that is
\[\mathbb{U}=\left\{u=i_1\ldots i_k\,;\,\text{$k\in\mathbb{N}$ and $i_1,\ldots,i_k\in\mathbb{N}^*$}\right\},\]
in which the root is the empty word $\varnothing$, and for all $u\in\mathbb{U}$ and $i\in\mathbb{N}^*$, the word $ui$ is the $i$-th children of $u$.
We denote by $|u|$ the length of the word $u\in\mathbb{U}$, which corresponds to the height of $u$ in the Ulam--Harris tree, and we set
\[[\mathbb{U}]_n=\{u\in\mathbb{U}:|u|\leq n\}\quad\text{and}\quad\partial[\mathbb{U}]_n=\{u\in\mathbb{U}:|u|=n\}.\]}
Then, set ${\mathcal{B}_n=\{u\in\partial[\mathcal{T}]_n:X_u=1\}}$, and let
\[\mathcal{T}_n'=\{u\in[\mathcal{T}]_n:\text{there exists $v\in\mathcal{B}_n$ such that $v\succeq u$}\}.\]

\begin{defn}[Pruned tree]
We denote by $\mathcal{T}_n^p$ the random tree $\mathcal{T}_n'$ conditionally on the event \eqref{eq:boundary}.
\end{defn}

As observed by Ayuso Ventura and Berger \cite[Proposition 5.1]{irenequentinising}, the random tree $\mathcal{T}_n^p$ has the same distribution as the genealogical tree associated with an explicit inhomogenous branching process.
To describe it, let us introduce some notation.
We denote by $\xi$ a random variable with distribution $\mu$, and we denote by $G:x\in[0,1]\mapsto\sum_{k\geq0}\mu(\{k\})\cdot x^k$ the generating function of $\mu$.
Then, for all $n\in\mathbb{N}^*$ and ${p\in{]0,1]}}$, we denote by $\mu_n^p$ the distribution of a random variable with distribution ${\mathrm{Binomial}\left(\xi,1-G^{n-1}(1-p)\right)}$ conditionally on $\xi$, conditioned on being strictly positive:\footnote{Observe that
\[\sum_{k\geq1}\mathbb{E}\left[\binom{\xi}{k}\cdot\left(1-G^{n-1}(1-p)\right)^k\cdot\left(G^{n-1}(1-p)\right)^{\xi-k}\right]=\mathbb{E}\left[1-\left(G^{n-1}(1-p)\right)^\xi\right]=1-G^n(1-p).\]}
\begin{equation}\label{eq:munp}
\mu_n^p(\{k\})=\frac{\mathbb{E}\left[\binom{\xi}{k}\cdot\left(1-G^{n-1}(1-p)\right)^k\cdot\left(G^{n-1}(1-p)\right)^{\xi-k}\right]}{1-G^n(1-p)}\quad\text{for all $k\in\mathbb{N}^*$.}
\end{equation}
Now, let $(\xi_n^p(u))_{u\in[\mathbb{U}]_{n-1}}$ be a collection of independent random variables such that for each $k\in\llbracket0,n-1\rrbracket$, the random variables $(\xi_n^p(u))_{u\in\partial[\mathbb{U}]_k}$ have distribution $\mu_{n-k}^p$.
Then, let $\mathbf{T}_n^p$ be the genealogical tree associated with the offspring random variables $\left(\xi_n^p(u)\right)_{u\in[\mathbb{U}]_{n-1}}$:
\begin{equation}\label{eq:genealogicalTnp}
\mathbf{T}_n^p=\left\{u=i_1\ldots i_k\in[\mathbb{U}]_n:\text{$i_{j+1}\leq\xi_n^p(i_1\ldots i_j)$ for all $j\in\llbracket0,k-1\rrbracket$}\right\}.
\end{equation}
The following proposition, due to Auyso Ventura and Berger, tells us that the pruned tree $\mathcal{T}_n^p$ has the same distribution as $\mathbf{T}_n^p$:

\begin{prop}[{\cite[Proposition 5.1]{irenequentinising}}]\label{prop:pruned}
For every finite plane tree $T$ with height $n$, we have
\[\mathbb{P}\left(\mathcal{T}_n^p=T\right)=\prod_{k=0}^{n-1}\prod_{\substack{u\in T\\|u|=k}}\mu_{n-k}^p(\#\{v\in T:\text{$v$ is a child of $u$}\})=\mathbb{P}(\mathbf{T}_n^p=T).\]
\end{prop}

\subsection{Basic estimates}

In this subsection, we collect some basic estimates on the pruned tree $\mathcal{T}_n^p$.
Depending on the context, it will be more convenient to use either the definition of $\mathcal{T}_n^p$ as the random tree $\mathcal{T}_n'$ conditionally on the event \eqref{eq:boundary}, or the fact that $\mathcal{T}_n^p$ has the same distribution as the genealogical tree $\mathbf{T}_n^p$ associated with the offspring random variables $(\xi_n^p(u))_{u\in[\mathbb{U}]_{n-1}}$.

\begin{lem}\label{lem:collection}
For all $n\in\mathbb{N}^*$ and $p\in{]0,1]}$, the following holds.
\begin{itemize}
\item We have
\begin{equation}\label{eq:meanTnp}
\mathbb{E}\left[\#\partial\mathcal{T}_n^p\right]=\frac{m^n\cdot p}{1-G^n(1-p)},
\end{equation}
and
\begin{equation}\label{eq:Tnp2}
\mathbb{E}\left[\left(\#\partial\mathcal{T}_n^p\right)^2\right]=\frac{\mathbb{E}[Z_n(Z_n-1)]\cdot p^2+m^n\cdot p}{1-G^n(1-p)}.
\end{equation}
Moreover, we have
\begin{equation}\label{eq:Tnpgeq2}
\mathbb{P}\left(\#\partial\mathcal{T}_n^p\geq2\right)\leq\frac{\mathbb{E}[Z_n(Z_n-1)]\cdot p^2}{1-G^n(1-p)}.
\end{equation}

\item We have
\begin{equation}\label{eq:meanxinp}
\mathbb{E}\left[\xi_n^p(\varnothing)\right]=\frac{m\cdot\left(1-G^{n-1}(1-p)\right)}{1-G^n(1-p)},
\end{equation}
and
\begin{equation}\label{eq:varxinp}
\mathbb{V}\left(\xi_n^p(\varnothing)\right)\leq\frac{\mathbb{E}[\xi(\xi-1)]\cdot\left(1-G^{n-1}(1-p)\right)^2}{1-G^n(1-p)}.
\end{equation}

\item For all $k\in\llbracket0,n\rrbracket$, we have
\begin{equation}\label{eq:Tnpgenk}
\mathbb{E}\left[\#\partial\left[\mathcal{T}_n^p\right]_k\right]=\frac{m^k\cdot\left(1-G^{n-k}(1-p)\right)}{1-G^n(1-p)}.
\end{equation}
\end{itemize}
\end{lem}
\begin{proof}
\begin{itemize}
\item By the definition of $\mathcal{T}_n^p$, we have
\[\mathbb{E}\left[\#\partial\mathcal{T}_n^p\right]=\frac{\mathbb{E}\left[\sum_{u\in\partial[\mathcal{T}]_n}X_u\right]}{1-G^n(1-p)}\quad\text{and}\quad\mathbb{E}\left[\left(\#\partial\mathcal{T}_n^p\right)^2\right]=\frac{\mathbb{E}\left[\left(\sum_{u\in\partial[\mathcal{T}]_n}X_u\right)^2\right]}{1-G^n(1-p)},\]
where $(X_u)_{u\in\partial[\mathbb{U}]_n}$ are independent Bernoulli random variables with parameter $p$, and are independent of $\mathcal{T}$.
Now, conditionally on $\mathcal{T}$, we have
\begin{equation}\label{eq:1stmomentcond}
\mathbb{E}\left[\sum_{u\in\partial[\mathcal{T}]_n}X_u\,\middle|\,\mathcal{T}\right]=Z_n\cdot p\quad\text{almost surely,}
\end{equation}
and
\begin{eqnarray}
\mathbb{E}\left[\left(\sum_{u\in\partial[\mathcal{T}]_n}X_u\right)^2\,\middle|\,\mathcal{T}\right]&=&\mathbb{E}\left[\sum_{\substack{u_1,u_2\in\partial[\mathcal{T}]_n\\u_1\neq u_2}}X_{u_1}X_{u_2}+\sum_{u\in\partial[\mathcal{T}]_n}X_u\,\middle|\,\mathcal{T}\right]\notag\\
&=&Z_n(Z_n-1)\cdot p_n^2+Z_n\cdot p_n\quad\text{almost surely.}\label{eq:2ndmomentcond}
\end{eqnarray}
This yields \eqref{eq:meanTnp} and \eqref{eq:Tnp2}.
Finally, we have
\begin{eqnarray*}
\mathbb{P}\left(\#\partial\mathcal{T}_n^p\geq2\right)&=&\frac{\mathbb{P}\left(\exists u_1\neq u_2\in\partial[\mathcal{T}]_n:\text{$X_{u_1}=1$ and $X_{u_2}=1$}\right)}{1-G^n(1-p)}\\
&\leq&\frac{\mathbb{E}\left[\sum_{\substack{u_1,u_2\in\partial[\mathcal{T}]_n\\u_1\neq u_2}}X_{u_1}X_{u_2}\right]}{1-G^n(1-p)}=\frac{\mathbb{E}[Z_n(Z_n-1)]\cdot p^2}{1-G^n(1-p)}.
\end{eqnarray*}

\item By the definition of $\mu_n^p$ (see \eqref{eq:munp}), we have
\[\mathbb{E}\left[\xi_n^p(\varnothing)\right]=\frac{\mathbb{E}\left[\xi\cdot\left(1-G^{n-1}(1-p)\right)\right]}{1-G^n(1-p)}=\frac{m\cdot\left(1-G^{n-1}(1-p)\right)}{1-G^n(1-p)},\]
and
\begin{eqnarray*}
\mathbb{E}\left[\xi_n^p(\varnothing)^2\right]&=&\frac{\mathbb{E}\left[\xi(\xi-1)\cdot\left(1-G^{n-1}(1-p)\right)^2+\xi\cdot\left(1-G^{n-1}(1-p)\right)\right]}{1-G^n(1-p)}\\
&=&\frac{\mathbb{E}[\xi(\xi-1)]\cdot\left(1-G^{n-1}(1-p)\right)^2}{1-G^n(1-p)}+\mathbb{E}\left[\xi_n^p(\varnothing)\right].
\end{eqnarray*}
Thus, since $\mathbb{E}[\xi_n^p(\varnothing)]\geq1$, we get
\begin{eqnarray*}
\mathbb{V}\left(\xi_n^p(\varnothing)\right)&=&\mathbb{E}\left[\xi_n^p(\varnothing)^2\right]-\mathbb{E}\left[\xi_n^p(\varnothing)\right]^2\notag\\
&\leq&\mathbb{E}\left[\xi_n^p(\varnothing)^2\right]-\mathbb{E}\left[\xi_n^p(\varnothing)\right]=\frac{\mathbb{E}[\xi(\xi-1)]\cdot\left(1-G^{n-1}(1-p)\right)^2}{1-G^n(1-p)}.
\end{eqnarray*}

\item By Proposition \ref{prop:pruned}, it suffices to show that for all ${k\in\llbracket0,n\rrbracket}$, we have
\begin{equation}\label{eq:recTnpgenk}
\mathbb{E}\left[\#\partial\left[\mathbf{T}_n^p\right]_k\right]=\frac{m^k\cdot\left(1-G^{n-k}(1-p)\right)}{1-G^n(1-p)}.
\end{equation}
We prove this by induction on $k$:
For $k=0$, we have $\mathbb{E}\left[\#\partial\left[\mathbf{T}_n^p\right]_k\right]=1$, hence \eqref{eq:recTnpgenk} holds.
Next, fix $k\in\llbracket0,n-1\rrbracket$ such that \eqref{eq:recTnpgenk} holds, and let us calculate the expectation of
\[\#\partial\left[\mathbf{T}_n^p\right]_{k+1}=\sum_{\substack{u\in\mathbf{T}_n^p\\|u|=k}}\xi_n^p(u).\]
Let $\mathcal{F}_k$ be the $\sigma$-algebra generated by the random variables $\left(\xi_n^p(u)\right)_{u\in[\mathbb{U}]_{k-1}}$.
Since $\left[\mathbf{T}_n^p\right]_k$ is measurable with respect to $\mathcal{F}_k$, and since the random variables $\left(\xi_n^p(u)\right)_{u\in\partial[\mathbb{U}]_k}$ are independent of $\mathcal{F}_k$, we have
\[\mathbb{E}\left[\#\partial\left[\mathbf{T}_n^p\right]_{k+1}\,\middle|\,\mathcal{F}_k\right]=\sum_{\substack{u\in\mathbf{T}_n^p\\|u|=k}}\mathbb{E}\left[\xi_n^p(u)\right]=\#\partial\left[\mathbf{T}_n^p\right]_k\cdot\mathbb{E}\left[\xi_{n-k}^p(\varnothing)\right].\]
Taking expectations, we get
\[\mathbb{E}\left[\#\partial\left[\mathbf{T}_n^p\right]_{k+1}\right]=\mathbb{E}\left[\#\partial\left[\mathbf{T}_n^p\right]_k\right]\cdot\mathbb{E}\left[\xi_{n-k}^p(\varnothing)\right].\]
Now, by the induction hypothesis \eqref{eq:recTnpgenk}, and by \eqref{eq:meanxinp}, we obtain
\begin{eqnarray*}
\mathbb{E}\left[\#\partial\left[\mathbf{T}_n^p\right]_{k+1}\right]&=&\frac{m^k\cdot\left(1-G^{n-k}(1-p)\right)}{1-G^n(1-p)}\cdot\frac{m\cdot\left(1-G^{n-k-1}(1-p)\right)}{1-G^{n-k}(1-p)}\\
&=&\frac{m^{k+1}\cdot\left(1-G^{n-(k+1)}(1-p)\right)}{1-G^n(1-p)},
\end{eqnarray*}
which concludes the proof by induction.

\emph{Nota bene.} Later on, we will shorten such arguments into: by the branching property of $\mathcal{T}_n^p$, we have
\[\mathbb{E}\left[\#\partial[\mathcal{T}_n^p]_{k+1}\right]=\mathbb{E}[\#\partial[\mathcal{T}_n^p]_k]\cdot\mathbb{E}\left[\xi_{n-k}^p(\varnothing)\right].\]
\end{itemize}
\end{proof}

In the next lemma, we provide some estimates on the non-empty boundary probability ${1-G^n(1-p)}$, which we will use repeatedly throughout the paper.

\begin{lem}\label{lem:asymptoticsurvival}
\begin{itemize}
\item There exists a constant $c\in{]0,1[}$ such that
\[c\cdot\left(\left(m^n\cdot p\right)\wedge1\right)\leq1-G^n(1-p)\leq\left(m^n\cdot p\right)\wedge1\quad\text{for all $n\in\mathbb{N}$ and $p\in[0,1]$.}\]

\item For each $n\in\mathbb{N}^*$, fix $p_n\in{]0,1]}$.
If $m^n\cdot p_n\rightarrow\infty$ as $n\to\infty$, then $1-G^n(1-p_n)\rightarrow1$ as $n\to\infty$.
\end{itemize}
\end{lem}
\begin{proof}
\begin{itemize}
\item Fix $p\in[0,1]$.
It is easy to check by induction on $n$ that for all $n\in\mathbb{N}$, we have ${1-G^n(1-p)=F^n(p)}$, where ${F:x\in[0,1]\mapsto1-G(1-x)}$.
Now, for all $x\in[0,1]$, we have
\begin{equation}\label{eq:F}
m\cdot x-\frac{\mathbb{E}[\xi(\xi-1)]}{2}\cdot x^2\leq F(x)\leq m\cdot x.
\end{equation}
In particular, we get that $F^n(p)\leq m^n\cdot p$ for all $n\in\mathbb{N}$, which yields the desired upper bound $F^n(p)\leq(m^n\cdot p)\wedge1$.
To prove the lower bound, set
\[c=\log_m\left(\frac{\mathbb{E}[\xi(\xi-1)]}{m(m-1)}\right)\in\mathbb{R}_+.\]
First, we claim that
\begin{equation}\label{eq:lowerboundFnp1}
F^n(p)\geq\frac{1}{2}\cdot m^n\cdot p\quad\text{for all $n\in\mathbb{N}$ such that $n\leq\log_m\left(\frac{1}{p}\right)-c$.}
\end{equation}
Indeed, fix $n\in\mathbb{N}$ such that $n\leq\log_m(1/p)-c$.
By \eqref{eq:F}, we have
\begin{equation}\label{eq:telescopeF}
\frac{F^k(p)}{m^k\cdot p}-\frac{F^{k+1}(p)}{m^{k+1}\cdot p}\leq\frac{\mathbb{E}[\xi(\xi-1)]}{2}\cdot\frac{F^k(p)^2}{m^{k+1}\cdot p}\leq\frac{\mathbb{E}[\xi(\xi-1)]}{2m}\cdot m^k\cdot p\quad\text{for all $k\in\llbracket0,n-1\rrbracket$.}
\end{equation}
Next, by the definition of $c$, we have
\[\sum_{k=0}^{n-1}\frac{\mathbb{E}[\xi(\xi-1)]}{2m}\cdot m^k\cdot p\leq\frac{\mathbb{E}[\xi(\xi-1)]}{2m}\cdot\frac{m^n}{m-1}\cdot p\leq\frac{\mathbb{E}[\xi(\xi-1)]}{2m(m-1)}\cdot m^{\log_m(1/p)-c}\cdot p=\frac{1}{2}.\]
Thus, summing \eqref{eq:telescopeF} over $k$, we obtain
\[1-\frac{F^n(p)}{m^n\cdot p}\leq\frac{1}{2},\]
which yields \eqref{eq:lowerboundFnp1}.
Next, we claim that
\begin{equation}\label{eq:lowerboundFnp2}
F^n(p)\geq\frac{1}{2}\cdot m^{-c-1}\quad\text{for all $n\in\mathbb{N}$ such that $n\geq\log_m\left(\frac{1}{p}\right)-c$.}
\end{equation}
Indeed:
\begin{itemize}
\item If $\log_m(1/p)-c<0$, then $p\geq m^{-c}$, hence for all $n\in\mathbb{N}$ such that $n\geq\log_m(1/p)-c$, we have $F^n(p)\geq F^0(p)=p\geq m^{-c}$.

\item If $\log_m(1/p)-c\geq0$, then $k:=\lfloor\log_m(1/p)-c\rfloor\in\mathbb{N}$, hence for all $n\in\mathbb{N}$ such that $n\geq\log_m(1/p)-c$, we have
\[F^n(p)\geq F^k(p)\geq\frac{1}{2}\cdot m^k\cdot p\geq\frac{1}{2}\cdot m^{\log_m(1/p)-c-1}\cdot p=\frac{1}{2}\cdot m^{-c-1}.\]
\end{itemize}
In any case, we obtain \eqref{eq:lowerboundFnp2}.
Combining \eqref{eq:lowerboundFnp1} and \eqref{eq:lowerboundFnp2}, we deduce that
\[F^n(p)\geq\left(\frac{1}{2}\cdot m^n\cdot p\right)\wedge\left(\frac{1}{2}\cdot m^{-c-1}\right)\quad\text{for all $n\in\mathbb{N}$.}\]
It follows that there exists a constant $c'\in{]0,1[}$ such that
\[F^n(p)\geq c'\cdot((m^n\cdot p)\wedge1)\quad\text{for all $n\in\mathbb{N}$.}\]

\item First, note that for each $x\in{]0,1]}$, we have $F^k(x)\rightarrow1$ as $k\to\infty$.
Indeed, since $\mu(\{0\})=0$, we have $G(x)\leq x$ for all $x\in[0,1]$, hence $F(x)\geq x$ for all $x\in[0,1]$.
It follows that for each $x\in{]0,1]}$, the sequence $\left(F^k(x)\right)_{k\in\mathbb{N}}$ is non-decreasing.
Since it is bounded by $1$, it converges, and by the continuity of $F$, its limit $y\in[0,1]$ must be a fixed point of $F$, i.e, $1-y$ must be a fixed point of $G$.
Since $\mu(\{0\})=0$ and $m>1$, the set of fixed points of $G$ is $\{0,1\}$, hence $y\in\{0,1\}$.
Since $y\geq x>0$, we conclude that $y=1$, as desired.

Now, for each $n\in\mathbb{N}^*$, let $\ell_n=\lceil\log_m(1/p_n)\rceil$ be the smallest integer $k\in\mathbb{N}$ such that $m^k\cdot p_n\geq1$.
By the previous point, there exists a constant $c\in{]0,1[}$ such that $F^{\ell_n}(p_n)\geq c$ for all $n\in\mathbb{N}^*$.
Then, since
\[n-\ell_n\geq n-\log_m\left(\frac{1}{p_n}\right)-1=\log_m\left(m^n\cdot p_n\right)-1\longrightarrow\infty\quad\text{as $n\to\infty$,}\]
we get that
\[1-G^n(1-p_n)=F^n(p_n)=F^{n-\ell_n}\left(F^{\ell_n}(p_n)\right)\geq F^{n-\ell_n}(c)\longrightarrow1\quad\text{as $n\to\infty$.}\]
\end{itemize}
\end{proof}

The following lemma shows that $\#\partial\mathcal{T}_n^p\asymp1\vee(m^n\cdot p)$ in probability, which will be useful throughout the paper.

\begin{lem}\label{lem:leavesTnp}
The sequence of random variables
\[\left(\frac{\#\partial\mathcal{T}_n^p}{1\vee\left(m^n\cdot p\right)}\right)_{\substack{n\in\mathbb{N}^*\\p\in{]0,1]}}}\]
is tight in $\mathbb{R}_+^*$: for each $\varepsilon\in{]0,1[}$, there exists $\delta<\Delta\in\mathbb{R}_+^*$ such that
\[\mathbb{P}\left(\delta\leq\frac{\#\partial\mathcal{T}_n^p}{1\vee\left(m^n\cdot p\right)}\leq\Delta\right)\geq1-\varepsilon\quad\text{for all $n\in\mathbb{N}^*$ and $p\in{]0,1]}$.}\]
\end{lem}
\begin{proof}
By \eqref{eq:meanTnp} and by Lemma \ref{lem:asymptoticsurvival}, there exists a constant $c\in{]0,1[}$ such that for all $n\in\mathbb{N}^*$ and $p\in{]0,1]}$, we have
\[\mathbb{E}\left[\#\partial\mathcal{T}_n^p\right]=\frac{m^n\cdot p}{1-G^n(1-p)}\leq\frac{m^n\cdot p}{c\cdot\left(\left(m^n\cdot p\right)\wedge 1\right)}=\frac{1}{c}\cdot\left(1\vee\left(m^n\cdot p\right)\right).\]
By Markov's inequality, it follows that for each $\varepsilon\in{]0,1[}$, we have
\[\mathbb{P}\left(\frac{\#\partial\mathcal{T}_n^p}{1\vee\left(m^n\cdot p\right)}\geq\frac{1}{c\varepsilon}\right)\leq\varepsilon.\]
Now, let us show that
\[\sup_{\substack{n\in\mathbb{N}^*\\p\in{]0,1]}}}\mathbb{P}\left(\frac{\#\partial\mathcal{T}_n^p}{1\vee\left(m^n\cdot p\right)}\leq\delta\right)\longrightarrow0\quad\text{as $\delta\to0^+$.}\]
To this end, fix $\delta\in{]0,1[}$, and let us upper bound the probability
\[\pi_n^p(\delta):=\mathbb{P}\left(\frac{\#\partial\mathcal{T}_n^p}{1\vee\left(m^n\cdot p\right)}\leq\delta\right)\]
uniformly in $n\in\mathbb{N}^*$ and $p\in{]0,1]}$.
First, if $m^n\cdot p<1/\delta$, then $\delta\cdot(1\vee\left(m^n\cdot p\right))<1$, and since $\#\partial\mathcal{T}_n^p\geq1$, we get $\pi_n^p(\delta)=0$.
Now, assume that $m^n\cdot p\geq1/\delta$.
Then, by the definition of $\mathcal{T}_n^p$, we have
\[\pi_n^p(\delta)=\mathbb{P}\left(\#\partial\mathcal{T}_n^p\leq\delta\cdot m^n\cdot p\right)=\frac{\mathbb{P}\left(1\leq\sum_{u\in\partial[\mathcal{T}]_n}X_u\leq\delta\cdot m^n\cdot p\right)}{1-G^n(1-p)},\]
where $(X_u)_{u\in\partial[\mathbb{U}]_n}$ are independent Bernoulli random variables with parameter $p$, and are independent of $\mathcal{T}$.
By Lemma \ref{lem:asymptoticsurvival}, there exists a constant $c\in{]0,1[}$ such that $1-G^n(1-p)\geq c$, which settles the denominator.
To handle the numerator, we write
\begin{eqnarray*}
\lefteqn{\mathbb{P}\left(1\leq\sum_{u\in\partial[\mathcal{T}]_n}X_u\leq\delta\cdot m^n\cdot p\right)}\\
&\leq&\mathbb{P}\left(Z_n\leq\sqrt{\delta}\cdot m^n\right)+\mathbb{P}\left(\mathrm{Binomial}\left(\left\lceil\sqrt{\delta}\cdot m^n\right\rceil,p\right)\leq\delta\cdot m^n\cdot p\right).
\end{eqnarray*}
Now, on the one hand, by \eqref{eq:Znmntight}, we have ${\sup_{n\geq0}\mathbb{P}\left(Z_n\leq\sqrt{\delta}\cdot m^n\right)\rightarrow0}$ as $\delta\to0$.
On the other hand, assuming further that $\delta\leq1/4$, so that
\[\delta\cdot m^n\cdot p\leq\frac{\sqrt{\delta}\cdot m^n\cdot p}{2}\leq\frac{\left\lceil\sqrt{\delta}\cdot m^n\right\rceil\cdot p}{2},\]
we get by a routine exponential Markov inequality argument (recall that $m^n\cdot p\geq1/\delta$):
\begin{eqnarray*}
\lefteqn{\mathbb{P}\left(\mathrm{Binomial}\left(\left\lceil\sqrt{\delta}\cdot m^n\right\rceil,p\right)\leq\delta\cdot m^n\cdot p\right)}\\
&\leq&\mathbb{P}\left(\mathrm{Binomial}\left(\left\lceil\sqrt{\delta}\cdot m^n\right\rceil,p\right)\leq\frac{\left\lceil\sqrt{\delta}\cdot m^n\right\rceil\cdot p}{2}\right)\\
&\leq&\exp\left(-\eta\cdot\left\lceil\sqrt{\delta}\cdot m^n\right\rceil\cdot p\right)\\
&\leq&\exp\left(-\eta\cdot\sqrt{\delta}\cdot m^n\cdot p\right)\leq\exp\left(-\frac{\eta}{\sqrt\delta}\right),
\end{eqnarray*}
where $\eta=(1/2)\ln(1/2)-1/2+1>0$.
Altogether, we obtain
\[\pi_n^p(\delta)\leq\frac{1}{c}\cdot\left(\sup_{n\geq0}\mathbb{P}\left(Z_n\leq\sqrt{\delta}\cdot m^n\right)+\exp\left(-\frac{\eta}{\sqrt{\delta}}\right)\right).\]
This completes the proof, since the right-hand side does not depend on $n$ or $p$, and goes to zero as $\delta\to0$.
\end{proof}

The next proposition complements the previous result:

\begin{prop}\label{prop:leavesTnp'}
For each $n\in\mathbb{N}^*$, fix $p_n\in{]0,1]}$.
\begin{itemize}
\item If $(m^n\cdot p_n)/n\rightarrow0$ as $n\to\infty$, then
\[\frac{\#\partial\mathcal{T}_n^{p_n}}{n}\longrightarrow0\quad\text{in probability as $n\to\infty$.}\]

\item If $m^n\cdot p_n\rightarrow\infty$ as $n\to\infty$, then we have
\[\frac{\#\partial\mathcal{T}_n^{p_n}}{m^n\cdot p_n}\longrightarrow W\quad\text{in distribution as $n\to\infty$,}\]
where $W$ is given by \eqref{eq:W}.
\end{itemize}
\end{prop}
\begin{proof}
\begin{itemize}
\item By Lemma \ref{lem:leavesTnp}, the sequence of random variables
\[\left(\frac{\#\partial\mathcal{T}_n^{p_n}}{1\vee(m^n\cdot p_n)}\right)_{n\in\mathbb{N}^*}\]
is tight.
Therefore, if $(m^n\cdot p_n)/n\rightarrow0$ as $n\to\infty$, then $(1\vee(m^n\cdot p_n))/n\rightarrow0$ as $n\to\infty$, and thus
\[\frac{\#\partial\mathcal{T}_n^{p_n}}{n}=\frac{1\vee(m^n\cdot p_n)}{n}\cdot\frac{\#\partial\mathcal{T}_n^{p_n}}{1\vee(m^n\cdot p_n)}\longrightarrow0\quad\text{in probability as $n\to\infty$.}\]

\item For each $n\in\mathbb{N}^*$, by the definition of the pruned tree $\mathcal{T}_n^{p_n}$, the following holds: for any bounded continuous function $\varphi:\mathbb{R}\rightarrow\mathbb{R}$, we have
\[\mathbb{E}\left[\varphi\left(\frac{\#\partial\mathcal{T}_n^{p_n}}{m^n\cdot p_n}\right)\right]=\mathbb{E}\left[\varphi\left(\frac{\sum_{u\in\partial[\mathcal{T}]_n}X_u}{m^n\cdot p_n}\right)\right]\cdot\frac{1}{1-G^n(1-p_n)},\]
where $(X_u)_{u\in\partial[\mathbb{U}]_n}$ are independent Bernoulli random variables with parameter $p_n$, and are independent of $\mathcal{T}$.
Since $m^n\cdot p_n\rightarrow\infty$ as $n\to\infty$, by Lemma \ref{lem:asymptoticsurvival}, we have ${1-G^n(1-p_n)\rightarrow1}$ as $n\to\infty$.
Therefore, it suffices to show that
\[\frac{\sum_{u\in\partial[\mathcal{T}]_n}X_u}{m^n\cdot p_n}\longrightarrow W\quad\text{in probability as $n\to\infty$.}\]
To prove this, let us show that (see, e.g, \cite[Exercise 11.7 (i)]{legall})
\[\mathbb{E}\left[\left(\frac{\sum_{u\in\partial[\mathcal{T}]_n}X_u}{m^n\cdot p_n}-W\right)^2\,\middle|\,\mathcal{T}\right]\longrightarrow0\quad\text{in probability as $n\to\infty$.}\]
Since $\sum_{u\in\partial[\mathcal{T}]_n}X_u$ and $W$ are square-integrable, and since $W$ is measurable with respect to $\mathcal{T}$, we have
\begin{eqnarray*}
\lefteqn{\mathbb{E}\left[\left(\frac{\sum_{u\in\partial[\mathcal{T}]_n}X_u}{m^n\cdot p_n}-W\right)^2\,\middle|\,\mathcal{T}\right]}\\
&=&\frac{\mathbb{E}\left[\left(\sum_{u\in\partial[\mathcal{T}]_n}X_u\right)^2\,\middle|\,\mathcal{T}\right]}{(m^n\cdot p_n)^2}-2\cdot\frac{\mathbb{E}\left[\sum_{u\in\partial[\mathcal{T}]_n}X_u\,\middle|\,\mathcal{T}\right]}{m^n\cdot p_n}\cdot W+W^2\quad\text{almost surely.}
\end{eqnarray*}
By \eqref{eq:1stmomentcond} and by \eqref{eq:2ndmomentcond}, we get that
\begin{eqnarray*}
\lefteqn{\mathbb{E}\left[\left(\frac{\sum_{u\in\partial[\mathcal{T}]_n}X_u}{m^n\cdot p_n}-W\right)^2\,\middle|\,\mathcal{T}\right]}\\
&=&\frac{Z_n(Z_n-1)\cdot p_n^2+Z_n\cdot p_n}{(m^n\cdot p_n)^2}-2\cdot\frac{Z_n\cdot p_n}{m^n\cdot p_n}\cdot W+W^2\\
&=&\frac{Z_n}{m^n}\cdot\left(\frac{Z_n}{m^n}-\frac{1}{m^n}\right)+\frac{Z_n}{m^n}\cdot\frac{1}{m^n\cdot p_n}-2\cdot\frac{Z_n}{m^n}\cdot W+W^2.
\end{eqnarray*}
Thus, since $m^n\cdot p_n\rightarrow\infty$ as $n\to\infty$, we obtain that
\[\mathbb{E}\left[\left(\frac{\sum_{u\in\partial[\mathcal{T}]_n}X_u}{m^n\cdot p_n}-W\right)^2\,\middle|\,\mathcal{T}\right]\longrightarrow W\cdot W+W\cdot0-2\cdot W\cdot W+W^2=0\]
almost surely as $n\to\infty$, which yields the desired result.
\end{itemize}
\end{proof}

\begin{rem}
Furthermore, it can be shown that if $m^n\cdot p_n\rightarrow c\in\mathbb{R}_+$ as $n\to\infty$, then $\#\partial\mathcal{T}_n^{p_n}$ converges in distribution towards a ${\mathrm{Poisson}(W\cdot c)}$ random variable conditioned on being strictly positive: for each $k\in\mathbb{N}^*$, we have
\[\mathbb{P}\left(\#\partial\mathcal{T}_n^{p_n}=k\right)\longrightarrow\frac{\mathbb{E}\left[\frac{(W\cdot c)^k}{k!}\cdot e^{-W\cdot c}\right]}{\mathbb{E}\left[1-e^{-W\cdot c}\right]}\quad\text{as $n\to\infty$.}\]
In particular, if $c=0$, then we get that $\#\partial\mathcal{T}_n^{p_n}\rightarrow1$ in probability as $n\to\infty$.
Since this is not directly useful for our purposes, we omit the proof.
\end{rem}

\subsection{Key lemma}

The basic estimates of the previous subsection culminate in the following key lemma, which we isolate in this subsection:

\begin{lem}\label{lem:cut-height}
Let $K_n^p$ be the smallest integer $k\in\llbracket0,n\rrbracket$ such that every node $u\in\mathcal{T}_n^p$ with ${k\leq|u|\leq n-1}$ has exactly one child,\footnote{Note that this property holds vacuously for $k=n$.} and set, in view of \eqref{eq:gamman}:
\[\Gamma_n^p:=\frac{1}{\#\mathcal{T}_n^p}\cdot\sum_{k=1}^{K_n^p}k^3\cdot\max_{\substack{u\in\mathcal{T}_n^p\\|u|=k}}\#\left\{v\in\partial\mathcal{T}_n^p:v\succeq u\right\}.\]
\begin{itemize}
\item There exists a constant $C\in\mathbb{R}_+^*$ such that for all $n\in\mathbb{N}^*$ and $p\in{]0,1]}$, we have
\[\mathbb{P}\left(K_n^p>x\right)\leq C\cdot\frac{\left(m^n\cdot p\right)\vee\left(m^n\cdot p\right)^2}{m^x}\quad\text{for all $x\in\mathbb{R}_+$.}\]

\item The collection of random variables $\left(\Gamma_n^p\right)_{\substack{n\in\mathbb{N}^*\\p\in{]0,1]}}}$ is tight.
\end{itemize}
\end{lem}
\begin{proof}
\begin{itemize}
\item Fix $n\in\mathbb{N}^*$ and $p\in{]0,1]}$, and let $k\in\llbracket0,n-1\rrbracket$.
By the branching property of $\mathcal{T}_n^p$, by \eqref{eq:Tnpgenk}, and by \eqref{eq:Tnpgeq2}, we have
\begin{eqnarray*}
\mathbb{P}\left(K_n^p>k\right)&=&\mathbb{P}\left(\exists u\in\partial[\mathcal{T}_n^p]_k:\#\{v\in\partial\mathcal{T}_n^p:v\succeq u\}\geq2\right)\\
&\leq&\mathbb{E}\left[\sum_{\substack{u\in\mathcal{T}_n^p\\|u|=k}}\mathbf{1}\left(\#\left\{v\in\partial\mathcal{T}_n^p:v\succeq u\right\}\geq2\right)\right]\\
&=&\mathbb{E}\left[\#\partial\left[\mathcal{T}_n^p\right]_k\right]\cdot\mathbb{P}\left(\#\partial\mathcal{T}_{n-k}^p\geq2\right)\\
&\leq&\frac{m^k\cdot\left(1-G^{n-k}(1-p)\right)}{1-G^n(1-p)}\cdot\frac{\mathbb{E}[Z_{n-k}(Z_{n-k}-1)]\cdot p^2}{1-G^{n-k}(1-p)}\\
&=&\frac{m^k\cdot\mathbb{E}[Z_{n-k}(Z_{n-k}-1)]\cdot p^2}{1-G^n(1-p)}.
\end{eqnarray*}
Now, by \eqref{eq:Zn2} and by Lemma \ref{lem:asymptoticsurvival}, there exists constants $C\in\mathbb{R}_+^*$ and $c\in{]0,1[}$ such that
\[\mathbb{P}\left(K_n^p>k\right)\leq\frac{m^k\cdot C\cdot m^{2(n-k)}\cdot p^2}{c\cdot\left(\left(m^n\cdot p\right)\wedge1\right)}=\frac{C}{c}\cdot\frac{\left(m^n\cdot p\right)^2}{\left(m^n\cdot p\right)\wedge1}\cdot m^{-k}.\]
The desired result follows.

\item Fix $A\in\mathbb{R}_+^*$, and let us upper bound $\mathbb{P}\left(\Gamma_n^p>A\right)$ uniformly in $n\in\mathbb{N}^*$ and $p\in{]0,1]}$.
By the previous point, there exists a constant $C\in\mathbb{R}_+^*$ such that
\begin{equation}\label{eq:cuttight}
\mathbb{P}\left(K_n^p>x\right)\leq C\cdot\frac{(m^n\cdot p)\vee(m^n\cdot p)^2}{m^x}\quad\text{for all $x\in\mathbb{R}_+$.}
\end{equation}
Now, if $m^n\cdot p\leq1$, then we simply bound $\Gamma_n^p\leq\sum_{k=1}^{K_n^p}k^3\leq\left(K_n^p\right)^4$.
By \eqref{eq:cuttight}, it follows that
\[\mathbb{P}\left(\Gamma_n^p>A\right)\leq\mathbb{P}\left(K_n^p>A^{1/4}\right)\leq C\cdot m^{-A^{1/4}}.\]
Next, assume that $m^n\cdot p\geq1$.
Let $B\in\mathbb{R}_+^*$ be a parameter to be adjusted, and let $\ell_n^p$ be the largest integer $k\in\llbracket0,n\rrbracket$ such that $k\leq2\log_m\left(m^n\cdot p\right)+B$.
Moreover, let $\delta\in\mathbb{R}_+^*$ be a parameter to be adjusted.
We write
\begin{eqnarray*}
\mathbb{P}\left(\Gamma_n^p>A\right)&\leq&\mathbb{P}\left(K_n^p>2\log_m(m^n\cdot p)+B\right)+\mathbb{P}\left(\#\partial\mathcal{T}_n^p\leq\delta\cdot m^n\cdot p\right)\\
&&+\mathbb{P}\left(\sum_{k=1}^{\ell_n^p}k^3\cdot\max_{\substack{u\in\mathcal{T}_n^p\\|u|=k}}\#\left\{v\in\partial\mathcal{T}_n^p:v\succeq u\right\}>A\delta\cdot m^n\cdot p\right).
\end{eqnarray*}
Then, fix $\varepsilon\in{]0,1[}$, and let $B$ be large enough so that $C\cdot m^{-B}\leq\varepsilon$.
By \eqref{eq:cuttight}, we have
\[\mathbb{P}\left(K_n^p>2\log_m(m^n\cdot p)+B\right)\leq C\cdot\frac{(m^n\cdot p)^2}{m^{2\log_m(m^n\cdot p)+B}}=C\cdot m^{-B}\leq\varepsilon.\]
Moreover, by Lemma \ref{lem:leavesTnp}, we can fix $\delta$ small enough so that ${\mathbb{P}\left(\#\partial\mathcal{T}_n^p\leq\delta\cdot m^n\cdot p\right)\leq\varepsilon}$ for all $n\in\mathbb{N}^*$ and $p\in{]0,1]}$ such that $m^n\cdot p\geq1$.
Now that $B$ and $\delta$ are fixed, let us show that for $A$ large enough, we have
\[\mathbb{P}\left(\sum_{k=1}^{\ell_n^p}k^3\cdot\max_{\substack{u\in\mathcal{T}_n^p\\|u|=k}}\#\left\{v\in\partial\mathcal{T}_n^p:v\succeq u\right\}>A\delta\cdot m^n\cdot p\right)\leq\varepsilon\]
uniformly in $n$ and $p$.
To this end, we introduce, in hindsight:
\[a_k:=k\cdot\sqrt{m^k\wedge(m^n\cdot p)}\cdot\left(1\vee\left(m^{n-k}\cdot p\right)\right)\quad\text{for all $k\in\llbracket1,n\rrbracket$.}\]
First, we claim that there exists a constant $D\in\mathbb{R}_+^*$ such that
\begin{equation}\label{eq:D}
\sum_{k=1}^{\ell_n^p}k^3\cdot a_k\leq D\cdot m^n\cdot p.
\end{equation}
Indeed, let $\kappa_n^p$ be the largest integer $k\in\llbracket0,n\rrbracket$ such that $k\leq\log_m(m^n\cdot p)$ (note that $\kappa_n^p\leq\ell_n^p$).
By definition, we have
\begin{eqnarray*}
\sum_{k=1}^{\ell_n^p}k^3\cdot a_k&=&\sum_{k=1}^{\kappa_n^p}k^4\cdot m^{k/2}\cdot m^{n-k}\cdot p+\sum_{k=\kappa_n^p+1}^{\ell_n^p}k^4\cdot\sqrt{m^n\cdot p}\\
&=&\left(\sum_{k=1}^{\kappa_n^p}\frac{k^4}{m^{k/2}}\right)\cdot m^n\cdot p+\left(\sum_{k=\kappa_n^p+1}^{\ell_n^p}k^4\right)\cdot\sqrt{m^n\cdot p}\\
&\leq&\left(\sum_{k\geq1}\frac{k^4}{m^{k/2}}\right)\cdot m^n\cdot p+\left(\ell_n^p\right)^5\cdot\sqrt{m^n\cdot p}\\
&\leq&\left(\sum_{k\geq1}\frac{k^4}{m^{k/2}}\right)\cdot m^n\cdot p+\left(2\log_m(m^n\cdot p)+B\right)^5\cdot\sqrt{m^n\cdot p}.
\end{eqnarray*}
Since $m^n\cdot p\geq1$, this yields \eqref{eq:D}.
With \eqref{eq:D} at hand, we can upper bound
\begin{eqnarray*}
\lefteqn{\mathbb{P}\left(\sum_{k=1}^{\ell_n^p}k^3\cdot\max_{\substack{u\in\mathcal{T}_n^p\\|u|=k}}\#\left\{v\in\partial\mathcal{T}_n^p:v\succeq u\right\}>A\delta\cdot m^n\cdot p\right)}\\
&\leq&\mathbb{P}\left(\exists k\in\llbracket1,\ell_n^p\rrbracket:\max_{\substack{u\in\mathcal{T}_n^p\\|u|=k}}\#\{v\in\partial\mathcal{T}_n^p:v\succeq u\}>\frac{A\delta}{D}\cdot a_k\right).
\end{eqnarray*}
Now, let us first impose $A>D/\delta$, so that
\[\frac{A\delta}{D}\cdot a_n=\frac{A\delta}{D}\cdot n\cdot\sqrt{m^n\cdot p}>1.\]
Then, by the union bound, we get
\begin{eqnarray*}
\lefteqn{\mathbb{P}\left(\exists k\in\llbracket1,\ell_n^p\rrbracket:\max_{\substack{u\in\mathcal{T}_n^p\\|u|=k}}\#\{v\in\partial\mathcal{T}_n^p:v\succeq u\}>\frac{A\delta}{D}\cdot a_k\right)}\\
&=&\mathbb{P}\left(\exists k\in\llbracket1,\ell_n^p\wedge(n-1)\rrbracket:\max_{\substack{u\in\mathcal{T}_n^p\\|u|=k}}\#\{v\in\partial\mathcal{T}_n^p:v\succeq u\}>\frac{A\delta}{D}\cdot a_k\right)\\
&\leq&\sum_{k=1}^{\ell_n^p\wedge(n-1)}\mathbb{P}\left(\exists u\in\partial\left[\mathcal{T}_n^p\right]_k:\#\left\{v\in\partial\mathcal{T}_n^p:v\succeq u\right\}>\frac{A\delta}{D}\cdot a_k\right).
\end{eqnarray*}
Now, for each $k\in\llbracket1,\ell_n^p\wedge(n-1)\rrbracket$, by the union bound, by the branching property of $\mathcal{T}_n^p$, by \eqref{eq:Tnpgenk}, and by \eqref{eq:Tnp2}, we have
\begin{eqnarray*}
\lefteqn{\mathbb{P}\left(\exists u\in\partial[\mathcal{T}_n^p]_k:\#\{v\in\partial\mathcal{T}_n^p:v\succeq u\}>\frac{A\delta}{D}\cdot a_k\right)}\\
&\leq&\mathbb{E}\left[\sum_{\substack{u\in\mathcal{T}_n^p\\|u|=k}}\mathbf{1}\left(\#\left\{v\in\partial\mathcal{T}_n^p:v\succeq u\right\}>\frac{A\delta}{D}\cdot a_k\right)\right]\\
&=&\mathbb{E}\left[\#\partial\left[\mathcal{T}_n^p\right]_k\right]\cdot\mathbb{P}\left(\#\partial\mathcal{T}_{n-k}^p>\frac{A\delta}{D}\cdot a_k\right)\\
&\leq&\mathbb{E}\left[\#\partial\left[\mathcal{T}_n^p\right]_k\right]\cdot\frac{\mathbb{E}\left[\left(\#\partial\mathcal{T}_{n-k}^p\right)^2\right]}{((A\delta/D)\cdot a_k)^2}\\
&=&\frac{m^k\cdot\left(1-G^{n-k}(1-p)\right)}{1-G^n(1-p)}\cdot\frac{\mathbb{E}[Z_{n-k}(Z_{n-k}-1)]\cdot p^2+m^{n-k}\cdot p}{1-G^{n-k}(1-p)}\cdot\left(\frac{A\delta}{D}\cdot a_k\right)^{-2}\\
&=&\frac{m^k\cdot\left(\mathbb{E}[Z_{n-k}(Z_{n-k}-1)]\cdot p^2+m^{n-k}\cdot p\right)}{1-G^n(1-p)}\cdot\left(\frac{A\delta}{D}\cdot a_k\right)^{-2}.
\end{eqnarray*}
Thus, by \eqref{eq:Zn2} and by Lemma \ref{lem:asymptoticsurvival}, there exists constants $C'\in\mathbb{R}_+^*$ and $c\in{]0,1[}$ such that
\begin{eqnarray*}
\lefteqn{\mathbb{P}\left(\exists u\in\partial[\mathcal{T}_n^p]_k:\#\{v\in\partial\mathcal{T}_n^p:v\succeq u\}>\frac{A\delta}{D}\cdot a_k\right)}\\
&\leq&\frac{m^k\cdot C'\cdot\left(\left(m^{n-k}\cdot p\right)\vee\left(m^{n-k}\cdot p\right)^2\right)}{c}\cdot\left(\frac{A\delta}{D}\cdot a_k\right)^{-2}.
\end{eqnarray*}
By the definition of $\kappa_n^p$, it follows that
\begin{eqnarray*}
\lefteqn{\sum_{k=1}^{\ell_n^p\wedge(n-1)}\mathbb{P}\left(\exists u\in\partial[\mathcal{T}_n^p]_k:\#\{v\in\partial\mathcal{T}_n^p:v\succeq u\}>\frac{A\delta}{D}\cdot a_k\right)}\\
&\leq&\sum_{k=1}^{\kappa_n^p\wedge(n-1)}\frac{m^k\cdot C'\cdot\left(m^{n-k}\cdot p\right)^2}{c}\cdot\left(\frac{A\delta}{D}\cdot k\cdot m^{k/2}\cdot m^{n-k}\cdot p\right)^{-2}\\
&&+\sum_{k=\kappa_n^p+1}^{\ell_n^p\wedge(n-1)}\frac{m^k\cdot C'\cdot m^{n-k}\cdot p}{c}\cdot\left(\frac{A\delta}{D}\cdot k\cdot\sqrt{m^n\cdot p}\right)^{-2}\\
&=&\frac{C'D^2}{c(A\delta)^2}\cdot\sum_{k=1}^{\kappa_n^p\wedge(n-1)}\frac{1}{k^2}+\frac{C'D^2}{c(A\delta)^2}\cdot\sum_{k=\kappa_n^p+1}^{\ell_n^p\wedge(n-1)}\frac{1}{k^2}\leq\frac{C'D^2}{c(A\delta)^2}\cdot\sum_{k\geq1}\frac{1}{k^2}.
\end{eqnarray*}
This completes the proof, since the right-hand side does not depend on $n$ or $p$, and can be made arbitrarily small by choosing $A$ large enough.
\end{itemize}
\end{proof}

\begin{rem}\label{rem:cutheightbienaymé}
The previous proof can be easily adapted to show that the sequence of random variables
\[\left(\Gamma_n:=\frac{1}{\#\partial[\mathcal{T}]_n}\cdot\sum_{k=1}^nk^3\cdot\max_{\substack{u\in\mathcal{T}\\|u|=k}}\#\{v\in\partial[\mathcal{T}]_n:v\succeq u\}\right)_{n\in\mathbb{N}^*}\]
is tight, as claimed in the introduction.
Again, this is not directly useful for our purposes, so we omit the details.
\end{rem}

\section{The Gaussian case}\label{sec:dgff}

In this section, we consider the discrete Gaussian free field on $\mathcal{T}_n^{p_n}$, and prove Theorem \ref{thm:dgff}.

\subsection{Preliminaries}

We first collect some standard results on the discrete Gaussian free field.
Throughout this subsection, we fix a finite rooted tree $T_n$ with height $n\in\mathbb{N}^*$ and with no internal leaves, and we let ${G_n=(G_n(u))_{u\in T_n}}$ be the discrete Gaussian free field on $T_n$ with zero boundary conditions on $B_n=\partial T_n$.
In view of \eqref{eq:dgff}, the random vector $(G_n(u))_{u\in T_n\setminus\partial T_n}$ is a non-degenerate centred Gaussian vector.
Therefore, understanding the distribution of the height at the root of the tree $G_n(\varnothing)$ boils down to estimating its variance $\mathbb{V}(G_n(\varnothing))$, which can be calculated (in principle) by solving a recursion on the tree:

\begin{prop}\label{prop:green}
We have
\[\mathbb{V}(G_n(\varnothing))=\frac{1}{A_n(\varnothing)},\]
where $A_n:T_n\setminus\partial T_n\rightarrow\mathbb{R}_+^*$ satisfies the following recursion: we have
\[A_n(u)=\#\{v\in T_n:\text{$v$ is a child of $u$}\}\quad\text{for all $u\in T_n$ such that $|u|=n-1$,}\]
and
\[A_n(u)=\sum_{\substack{v\in T_n\\\text{$v$ child of $u$}}}\frac{A_n(v)}{1+A_n(v)}\quad\text{for all $u\in T_n$ such that $|u|\leq n-2$.}\]
\end{prop}

This result is a simple consequence of the fact that the function ${C_n:u\in T_n\mapsto\mathbb{E}[G_n(\varnothing)G_n(u)]}$ satisfies\footnote{We refer to \cite[Subsection 1.1.3]{wernerpowellgff} for a derivation of \eqref{eq:green}, which stems from the fact that for each $u\in T_n\setminus\partial T_n$, the random variable
\[G_n(u)-\frac{1}{\deg(u)}\cdot\sum_{\substack{v\in T_n\\v\sim u}}G_n(v)\sim\mathcal{N}\left(0,\frac{1}{\mathrm{deg}(u)}\right)\]
and the random vector $(G_n(v))_{v\in T_n\setminus\{u\}}$ are independent, which in turn follows from \eqref{eq:dgff}.}
\begin{equation}\label{eq:green}
\sum_{\substack{v\in T_n\\v\sim u}}(C_n(v)-C_n(u))=\begin{cases}
-1&\text{if $u=\emptyset$,}\\
0&\text{otherwise,}
\end{cases}\quad\text{for all ${u\in T_n\setminus\partial T_n}$.}
\end{equation}

\begin{proof}
With \eqref{eq:green} at hand, let us show that
\begin{equation}\label{eq:vardgffres}
C_n(\varnothing)=\frac{1}{A_n(\varnothing)}.
\end{equation}
If $n=1$, then we directly deduce from \eqref{eq:green} that
\[\sum_{\substack{v\in T_n\\\text{$v$ child of $\varnothing$}}}(0-C_n(\varnothing))=-1,\]
which yields \eqref{eq:vardgffres}, since $A_n(\varnothing)=\#\{v\in T_n:\text{$v$ is a child of $\varnothing$}\}$.
For $n\geq2$, let us first prove by induction on $k\in\{n-2,\ldots,0\}$ that
\begin{equation}\label{eq:reck}
C_n(u)=(1+A_n(v))\cdot C_n(v)\quad\text{for all $(u,v)\in\vec{E}(T_n)$ such that $k\leq|u|\leq n-2$.}
\end{equation}
To prove \eqref{eq:reck} for $k=n-2$, let $(u,v)\in\vec{E}(T_n)$ be such that $|u|=n-2$.
By \eqref{eq:green}, we have
\[(C_n(u)-C_n(v))+\sum_{\substack{w\in T_n\\\text{$w$ child of $v$}}}(0-C_n(v))=0.\]
Thus, we get
\[C_n(u)=C_n(v)+\#\{w\in T_n:\text{$w$ is a child of $v$}\}\cdot C_n(v)=(1+A_n(v))\cdot C_n(v),\]
as desired.
Next, fix $k\in\{n-2,\ldots,1\}$ such that \eqref{eq:reck} holds, and let us show that for all ${(u,v)\in\vec{E}(T_n)}$ with $|u|=k-1$, we have $C_n(u)=(1+A_n(v))\cdot C_n(v)$.
By \eqref{eq:green}, we have
\[(C_n(u)-C_n(v))+\sum_{\substack{w\in T_n\\\text{$w$ child of $v$}}}(C_n(w)-C_n(v))=0.\]
Now, by the induction hypothesis, for all children $w\in T_n$ of $v$, we have
\[C_n(w)-C_n(v)=\frac{C_n(v)}{1+A_n(w)}-C_n(v)=-\frac{A_n(w)}{1+A_n(w)}\cdot C_n(v).\]
It follows that
\[C_n(u)=C_n(v)+\sum_{\substack{w\in T_n\\\text{$w$ child of $v$}}}\frac{A_n(w)}{1+A_n(w)}\cdot C_n(v)=(1+A_n(v))\cdot C_n(v),\]
which proves \eqref{eq:reck} by induction.
To conclude, we specify \eqref{eq:green} to $u=\varnothing$: we get
\[\sum_{\substack{v\in T_n\\\text{$v$ child of $\varnothing$}}}(C_n(v)-C_n(\varnothing))=-1.\]
Now, by \eqref{eq:reck}, for all children $v\in T_n$ of $\varnothing$, we have
\[C_n(v)-C_n(\varnothing)=\frac{C_n(\varnothing)}{1+A_n(v)}-C_n(\varnothing)=-\frac{A_n(v)}{1+A_n(v)}\cdot C_n(\varnothing).\]
This finally yields \eqref{eq:vardgffres}, since
\[A_n(\varnothing)=\sum_{\substack{v\in T_n\\\text{$v$ child of $\varnothing$}}}\frac{A_n(v)}{1+A_n(v)}.\]
\end{proof}

\begin{rem}\label{rem:vardgfftrees}
\begin{itemize}
\item Using Proposition \ref{prop:green}, it is straightforward to calculate the variance of the height at the root of the tree for the discrete Gaussian free field on regular and starlike trees:
\begin{itemize}
\item If $T_n$ consists of the first $(n+1)$ generations of the $m$-ary tree, where ${m\geq2}$, then as claimed in the introduction, we have
\[\mathbb{V}(G_n(\varnothing))=\frac{1-m^{-n}}{m-1}.\]
Indeed, with the notation of Proposition \ref{prop:green}: by the spherical symmetry of $T_n$, we have ${A_n(u)=a_{n-|u|}}$ for all $u\in T_n\setminus\partial T_n$, where $(a_k)_{k\in\llbracket1,n\rrbracket}$ is given by $a_1=m$, and
\[a_{k+1}=m\cdot\frac{a_k}{1+a_k}\quad\text{for all $k\in\llbracket1,n-1\rrbracket$.}\]
The previous equation can be rearranged as follows:
\[\frac{m^{k+1}}{a_{k+1}}=m^k\cdot\left(\frac{1}{a_k}+1\right)\quad\text{for all $k\in\llbracket1,n-1\rrbracket$.}\]
Then, telescoping the above equality, we get
\[\frac{m^n}{a_n}=\sum_{k=1}^{n-1}\left(\frac{m^{k+1}}{a_{k+1}}-\frac{m^k}{a_k}\right)+\frac{m}{a_1}=\sum_{k=1}^{n-1}m^k+1=\sum_{k=0}^{n-1}m^k=\frac{m^n-1}{m-1},\]
which yields the desired result:
\[\frac{1}{a_n}=\frac{1-m^{-n}}{m-1}.\]

\item If $T_n$ is the tree with height $n$ in which the root has $N_n\in\mathbb{N}^*$ children, and every other node $u\in T_n$ with $|u|\leq n-1$ has one child, then as claimed in the introduction, we have
\[\mathbb{V}(G_n(\varnothing))=\frac{n}{N_n}.\]
Indeed, with the notation of Proposition \ref{prop:green}: by the spherical symmetry of $T_n$, we have ${A_n(u)=a_{n-|u|}}$ for all $u\in T_n\setminus\partial T_n$, where $(a_k)_{k\in\llbracket1,n\rrbracket}$ is given by $a_1=1$,
\[a_{k+1}=\frac{a_k}{1+a_k}\quad\text{for all $k\in\llbracket1,n-2\rrbracket$,}\]
and
\[a_n=N_n\cdot\frac{a_{n-1}}{1+a_{n-1}}.\]
The second equation can be rearranged as follows:
\[\frac{1}{a_{k+1}}=\frac{1}{a_k}+1\quad\text{for all $k\in\llbracket1,n-2\rrbracket$.}\]
Then, telescoping the above equality, we get
\[\frac{1}{a_{n-1}}=\sum_{k=1}^{n-2}\left(\frac{1}{a_{k+1}}-\frac{1}{a_k}\right)+\frac{1}{a_1}=\sum_{k=1}^{n-2}1+1=n-1,\]
which yields the desired result:
\[\frac{1}{a_n}=\frac{1}{N_n}\cdot\left(\frac{1}{a_{n-1}}+1\right)=\frac{n}{N_n}.\]
\end{itemize}

\item To emphasize that the order of magnitude of the variance of the height at the root of the tree really depends on the shape of the underlying tree, and not only, e.g, on the ratio of the number of leaves to the height of the tree, let us consider the following two examples: Fix $m\geq2$, fix $\kappa_n\in\llbracket0,n-1\rrbracket$, set $N_n:=m^{\kappa_n}$, and assume that $\kappa_n/n\rightarrow0$ as $n\to\infty$.

\begin{figure}[h]
\centering
\begin{tabular}{lcr}
\includegraphics[width=0.3\linewidth]{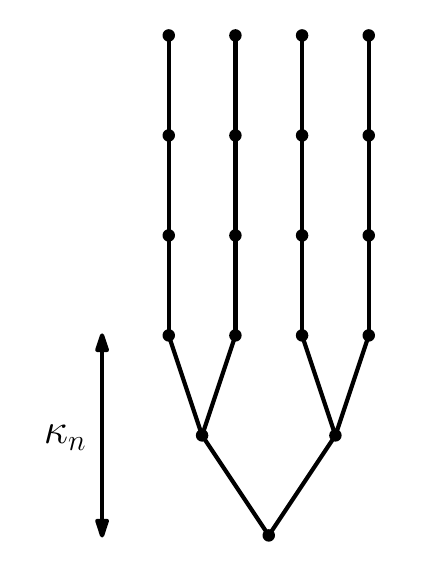}&&\includegraphics[width=0.3\linewidth]{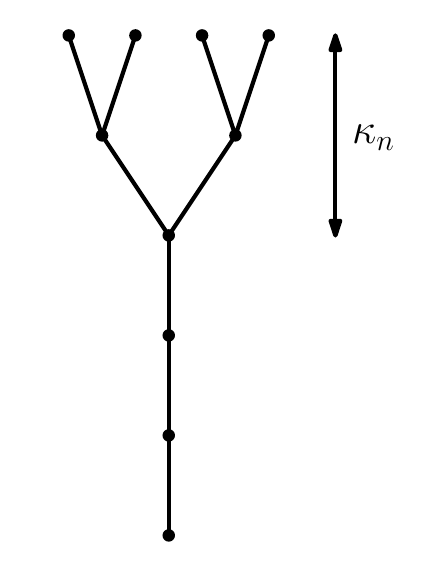}\\
\end{tabular}
\caption{Two examples of trees with the same ratio of number of leaves to height, for which the variance of the height at the root of the tree for the discrete Gaussian free field has a very different behaviour.}\label{fig:trees}
\end{figure}

\begin{enumerate}[label=(\arabic*)]
\item Let $T_n$ be the tree with height $n$ in which every node $u\in T_n$ with $|u|\leq\kappa_n-1$ has $m$ children, and every node $u\in T_n$ with $\kappa_n\leq|u|\leq n-1$ has one child (see Figure \ref{fig:trees} for an illustration).
Observe that $\#\partial T_n=N_n$.
We claim that
\[\mathbb{V}(G_n(\varnothing))\sim\frac{n}{N_n}\quad\text{as $n\to\infty$.}\]
Indeed, with the notation of Proposition \ref{prop:green}: by the spherical symmetry of $T_n$, we have $A_n(u)=a_{n-|u|}$ for all $u\in T_n\setminus\partial T_n$, where $(a_k)_{k\in\llbracket1,n\rrbracket}$ is given by $a_1=1$,
\[a_{k+1}=\frac{a_k}{1+a_k}\quad\text{for all $k\in\llbracket1,n-\kappa_n-1\rrbracket$,}\]
and
\[a_{k+1}=m\cdot\frac{a_k}{1+a_k}\quad\text{for all $k\in\llbracket n-\kappa_n,n-1\rrbracket$.}\]
With the same manipulations as above, we get that
\[\frac{1}{a_{n-\kappa_n}}=\sum_{k=1}^{n-\kappa_n-1}1+1=n-\kappa_n,\]
and
\[\frac{m^n}{a_n}=\sum_{k=n-\kappa_n}^{n-1}m^k+\frac{m^{n-\kappa_n}}{a_{n-\kappa_n}}=\frac{m^n-m^{n-\kappa_n}}{m-1}+m^{n-\kappa_n}\cdot(n-\kappa_n),\]
hence
\[\frac{1}{a_n}=\frac{1-1/N_n}{m-1}+\frac{n-\kappa_n}{N_n}\underset{n\to\infty}{\sim}\frac{n}{N_n}.\]

\item Let $T_n$ be the tree with height $n$ in which every node $u\in T_n$ with $|u|\leq n-\kappa_n-1$ has one child, and every node $u\in T_n$ with $n-\kappa_n\leq|u|\leq n-1$ has $m$ children (see Figure \ref{fig:trees} for an illustration).
Observe that $\#\partial T_n=N_n$.
We claim that
\[\mathbb{V}(G_n(\varnothing))\sim n\quad\text{as $n\to\infty$.}\]
Indeed, with the notation of Proposition \ref{prop:green}: by the spherical symmetry of $T_n$, we have $A_n(u)=a_{n-|u|}$ for all $u\in T_n\setminus\partial T_n$, where $(a_k)_{k\in\llbracket1,n\rrbracket}$ is given by $a_1=m$,
\[a_{k+1}=m\cdot\frac{a_k}{1+a_k}\quad\text{for all $k\in\llbracket1,\kappa_n-1\rrbracket$,}\]
and
\[a_{k+1}=\frac{a_k}{1+a_k}\quad\text{for all $k\in\llbracket\kappa_n,n-1\rrbracket$.}\]
With the same manipulations as above, we get that
\[\frac{m^{\kappa_n}}{a_{\kappa_n}}=\sum_{k=1}^{\kappa_n-1}m^k+\frac{m}{a_1}=\sum_{k=0}^{\kappa_n-1}m^k=\frac{m^{\kappa_n}-1}{m-1},\]
and
\[\frac{1}{a_n}=\sum_{k=\kappa_n}^{n-1}1+\frac{1}{a_{\kappa_n}}=(n-\kappa_n)+\frac{1-m^{-\kappa_n}}{m-1},\]
hence
\[\frac{1}{a_n}=(n-\kappa_n)+\frac{1-1/N_n}{m-1}\underset{n\to\infty}{\sim}n.\]
\end{enumerate}
\end{itemize}
\end{rem}

The result of Proposition \ref{prop:green} can also be seen as a manifestation of the link between the discrete Gaussian free field and electric networks: in fact, the variance of $G_n(\varnothing)$ corresponds to the effective resistance between $\varnothing$ and $\partial T_n$, which can be calculated (in principle) by using standard network reduction techniques (namely, the series and parallel laws).
We refer the interested reader to the authoritative reference \cite[Chapter 2]{lyonsperes}.
Another tool to estimate the variance of $G_n(\varnothing)$ is given by the following proposition:

\begin{prop}\label{prop:thomson}
We have the following variational representation:
\[\mathbb{V}(G_n(\varnothing))=\inf_{\varphi\in\Phi(T_n)}\sum_{(u,v)\in\vec{E}(T_n)}\varphi(u,v)^2,\]
where $\Phi(T_n)$ denotes the set of unit flows $\varphi$ on $T_n$ from $\varnothing$ to $\partial T_n$, i.e, functions $\varphi:\vec{E}(T_n)\rightarrow\mathbb{R}$ such that
\begin{equation}\label{eq:flow}
\varphi(u,v)=\sum_{\substack{w\in T_n\\\text{$w$ child of $v$}}}\varphi(v,w)\quad\text{for all $(u,v)\in\vec{E}(T_n)$ with $|v|\leq n-1$,}
\end{equation}
and
\[\sum_{\substack{v\in T_n\\\text{$v$ child of $\varnothing$}}}\varphi(\varnothing,v)=1.\]
\end{prop}

Given that $\mathbb{V}(G_n(\varnothing))$ corresponds to the effective resistance between $\varnothing$ and $\partial T_n$, this result is a direct consequence of Thomson's principle for electric networks.
For the convenience of the reader, we provide a self-contained proof of Proposition \ref{prop:thomson} in Appendix \ref{sec:thomson}.

\subsection{Localisation}

Now that we have all the necessary tools at hand, let us dive into the proof of Theorem \ref{thm:dgff}.
For each $n\in\mathbb{N}^*$, fix $p_n\in{]0,1]}$.
Then, conditionally on the random tree $\mathcal{T}_n^{p_n}$, let $G_n=(G_n(u))_{u\in\mathcal{T}_n^{p_n}}$ be the discrete Gaussian free field on $\mathcal{T}_n^{p_n}$ with zero boundary conditions on $\partial\mathcal{T}_n^{p_n}$.
In this subsection, we prove the localisation in Theorem \ref{thm:dgff}: we show that if
\begin{equation}\label{eq:asslocdgff}
\inf_{n\geq1}\frac{m^n\cdot p_n}{n}>0,
\end{equation}
then the sequence of random variables $\left(\mathbb{V}\left(G_n(\varnothing)\,\middle|\,\mathcal{T}_n^{p_n}\right)\right)_{n\in\mathbb{N}^*}$ is tight.

\begin{proof}[Proof of the localisation in Theorem \ref{thm:dgff}]
By Proposition \ref{prop:thomson}, we have
\[\mathbb{V}\left(G_n(\varnothing)\,\middle|\,\mathcal{T}_n^{p_n}\right)\leq\sum_{k=1}^n\sum_{\substack{u\in\mathcal{T}_n^{p_n}\\|u|=k}}\left(\frac{\#\left\{v\in\partial\mathcal{T}_n^{p_n}:v\succeq u\right\}}{\#\partial\mathcal{T}_n^{p_n}}\right)^2.\]
Thus, for all $B\in\mathbb{R}_+$, we have
\[\mathbb{P}\left(\mathbb{V}\left(G_n(\varnothing)\,\middle|\,\mathcal{T}_n^{p_n}\right)\geq B\right)\leq\mathbb{P}\left(\frac{1}{\left(\#\partial\mathcal{T}_n^{p_n}\right)^2}\cdot\sum_{k=1}^n\sum_{\substack{u\in\mathcal{T}_n^{p_n}\\|u|=k}}\#\left\{v\in\partial\mathcal{T}_n^{p_n}:v\succeq u\right\}^2\geq B\right).\]
Now, fix $\varepsilon\in{]0,1[}$.
By Lemma \ref{lem:leavesTnp}, and since $\inf_{n\geq1}(m^n\cdot p_n)>0$ by \eqref{eq:asslocdgff}, there exists $\delta\in\mathbb{R}_+^*$ such that $\mathbb{P}\left(\#\partial\mathcal{T}_n^{p_n}\leq\delta\cdot m^n\cdot p_n\right)\leq\varepsilon$ for all $n\in\mathbb{N}^*$.
It follows that
\begin{eqnarray*}
\lefteqn{\mathbb{P}\left(\mathbb{V}\left(G_n(\varnothing)\,\middle|\,\mathcal{T}_n^{p_n}\right)\geq B\right)}\\
&\leq&\mathbb{P}\left(\#\partial\mathcal{T}_n^{p_n}\leq\delta\cdot m^n\cdot p_n\right)+\mathbb{P}\left(\sum_{k=1}^n\sum_{\substack{u\in\mathcal{T}_n^{p_n}\\|u|=k}}\#\left\{v\in\mathcal{T}_n^{p_n}:v\succeq u\right\}^2\geq B\cdot(\delta\cdot m^n\cdot p_n)^2\right)\\
&\leq&\varepsilon+\frac{\mathbb{E}\left[\sum_{k=1}^n\sum_{\substack{u\in\mathcal{T}_n^{p_n}\\|u|=k}}\#\left\{v\in\mathcal{T}_n^{p_n}:v\succeq u\right\}^2\right]}{B\cdot(\delta\cdot m^n\cdot p_n)^2}.
\end{eqnarray*}
To complete the proof, let us show that there exists a constant $C\in\mathbb{R}_+^*$ such that
\begin{equation}\label{eq:upperboundthomson}
\mathbb{E}\left[\sum_{k=1}^n\sum_{\substack{u\in\mathcal{T}_n^{p_n}\\|u|=k}}\#\left\{v\in\partial\mathcal{T}_n^{p_n}:v\succeq u\right\}^2\right]\leq C\cdot(m^n\cdot p_n)^2\quad\text{for all $n\in\mathbb{N}^*$.}
\end{equation}
First, for $k=n$, by \eqref{eq:meanTnp}, we have
\[\mathbb{E}\left[\sum_{\substack{u\in\mathcal{T}_n^{p_n}\\|u|=n}}\#\left\{v\in\partial\mathcal{T}_n^{p_n}:v\succeq u\right\}^2\right]=\mathbb{E}\left[\#\partial\mathcal{T}_n^{p_n}\right]=\frac{m^n\cdot p_n}{1-G^n(1-p_n)}.\]
Next, for each $k\in\llbracket1,n-1\rrbracket$, by the branching property of $\mathcal{T}_n^{p_n}$, by \eqref{eq:Tnpgenk}, and by \eqref{eq:Tnp2}, we have
\begin{eqnarray*}
\lefteqn{\mathbb{E}\left[\sum_{\substack{u\in\mathcal{T}_n^{p_n}\\|u|=k}}\#\left\{v\in\partial\mathcal{T}_n^{p_n}:v\succeq u\right\}^2\right]}\\
&=&\mathbb{E}\left[\#\partial\left[\mathcal{T}_n^{p_n}\right]_k\right]\cdot\mathbb{E}\left[\left(\#\partial\mathcal{T}_{n-k}^{p_n}\right)^2\right]\\
&=&\frac{m^k\cdot\left(1-G^{n-k}(1-p_n)\right)}{1-G^n(1-p_n)}\cdot\frac{\mathbb{E}[Z_{n-k}(Z_{n-k}-1)]\cdot p_n^2+m^{n-k}\cdot p_n}{1-G^{n-k}(1-p)}\\
&=&\frac{m^k\cdot\mathbb{E}[Z_{n-k}(Z_{n-k}-1)]\cdot p_n^2+m^n\cdot p_n}{1-G^n(1-p_n)}.
\end{eqnarray*}
Thus, we get
\begin{eqnarray*}
\lefteqn{\mathbb{E}\left[\sum_{k=1}^n\sum_{\substack{u\in\mathcal{T}_n^{p_n}\\|u|=k}}\#\left\{v\in\partial\mathcal{T}_n^{p_n}:v\succeq u\right\}^2\right]}\\
&=&\sum_{k=1}^{n-1}\frac{m^k\cdot\mathbb{E}[Z_{n-k}(Z_{n-k}-1)]\cdot p_n^2+m^n\cdot p_n}{1-G^n(1-p_n)}+\frac{m^n\cdot p_n}{1-G^n(1-p_n)}\\
&=&\sum_{k=1}^{n-1}\frac{m^k\cdot\mathbb{E}[Z_{n-k}(Z_{n-k}-1)]\cdot p_n^2}{1-G^n(1-p_n)}+\frac{m^n\cdot p_n}{1-G^n(1-p_n)}\cdot n
\end{eqnarray*}
By \eqref{eq:Zn2} and by Lemma \ref{lem:asymptoticsurvival}, and since $\inf_{n\geq1}(m^n\cdot p_n)>0$ by \eqref{eq:asslocdgff}, it follows that there exists constants $C\in\mathbb{R}_+^*$ and $c\in{]0,1[}$ such that
\begin{eqnarray*}
\mathbb{E}\left[\sum_{k=1}^n\sum_{\substack{u\in\mathcal{T}_n^{p_n}\\|u|=k}}\#\left\{v\in\partial\mathcal{T}_n^{p_n}:v\succeq u\right\}^2\right]&\leq&\sum_{k=1}^{n-1}\frac{m^k\cdot C\cdot m^{2(n-k)}\cdot p_n^2}{c}+\frac{m^n\cdot p_n}{c}\cdot n\\
&=&\frac{C}{c}\cdot\left(\sum_{k=1}^{n-1}m^{-k}\right)\cdot(m^n\cdot p_n)^2+\frac{m^n\cdot p_n}{c}\cdot n\\
&\leq&\left(\frac{C}{c}\cdot\frac{1}{m-1}+\frac{1}{c}\cdot\frac{n}{m^n\cdot p_n}\right)\cdot(m^n\cdot p_n)^2.
\end{eqnarray*}
By \eqref{eq:asslocdgff}, we conclude that \eqref{eq:upperboundthomson} holds, which completes the proof.
\end{proof}

\subsection{Quantitative delocalisation}

In this subsection, we prove the quantitative delocalisation in Theorem \ref{thm:dgff}: we show that if
\begin{equation}\label{eq:assdelocdgff}
\frac{m^n\cdot p_n}{n}\longrightarrow0\quad\text{as $n\to\infty$,}
\end{equation}
then
\begin{equation}\label{eq:resdelocdgff}
\frac{\#\partial\mathcal{T}_n^{p_n}}{n}\cdot\mathbb{V}(G_n(\varnothing)\,|\,\mathcal{T}_n^{p_n})\longrightarrow1\quad\text{in probability as $n\to\infty$.}
\end{equation}
Since
\[\mathbb{P}\left(\frac{\sqrt{\partial\mathcal{T}_n^{p_n}}}{\sqrt{n}}\cdot G_n(\varnothing)\in\cdot\,\middle|\,\mathcal{T}_n^{p_n}\right)=\mathcal{N}\left(0,\frac{\#\partial\mathcal{T}_n^{p_n}}{n}\cdot\mathbb{V}(G_n(\varnothing)\,|\,\mathcal{T}_n^{p_n})\right),\]
this entails that
\[\mathbb{P}\left(\frac{\sqrt{\partial\mathcal{T}_n^{p_n}}}{\sqrt{n}}\cdot G_n(\varnothing)\in\cdot\,\middle|\,\mathcal{T}_n^{p_n}\right)\longrightarrow\mathcal{N}(0,1)\quad\text{in probability as $n\to\infty$}\]
in the Polish space of Borel probability measures on $\mathbb{R}$ endowed with the topology of weak convergence (see Remark \ref{rem:cvM1R} in the devoted Appendix \ref{sec:cvM1R}).
Recall from Proposition \ref{prop:leavesTnp'} that \eqref{eq:assdelocdgff} entails that $\left.\#\partial\mathcal{T}_n^{p_n}\middle/n\right.\rightarrow0$ in probability as $n\to\infty$.

Recall from Section \ref{sec:pruned} the definition of $\mathbf{T}_n^p$ as the genealogical tree associated with the offspring random variables $(\xi_n^p(u))_{u\in[\mathbb{U}]_{n-1}}$ (see \eqref{eq:genealogicalTnp}): recall that $\mathbf{T}_n^p$ has the same distribution as $\mathcal{T}_n^p$ (see Proposition \ref{prop:pruned}).
Now, by Proposition \ref{prop:green}, the variance of the height at the root of the tree for the discrete Gaussian free field on $\mathbf{T}_n^p$ with zero boundary conditions on $\partial\mathbf{T}_n^p$ is equal to $\left.1\middle/A_n^p(\varnothing)\right.$, where $\left(A_n^p(u)\right)_{u\in[\mathbb{U}]_{n-1}}$ is the collection of random variables constructed out of the offspring random variables $(\xi_n^p(u))_{u\in[\mathbb{U}]_{n-1}}$ by induction as follows: we set
\[A_n^p(u)=\xi_n^p(u)\quad\text{for all $u\in\mathbb{U}$ such that $|u|=n-1$,}\]
and we let
\[A_n^p(u)=\sum_{i=1}^{\xi_n^p(u)}\frac{A_n^p(ui)}{1+A_n^p(ui)}\quad\text{for all $u\in\mathbb{U}$ such that $|u|\leq n-2$.}\]
Therefore, proving \eqref{eq:resdelocdgff} is equivalent to showing that
\begin{equation}\label{eq:vardgffpruned}
\frac{n}{\#\partial\mathbf{T}_n^{p_n}}\cdot A_n^{p_n}(\varnothing)\longrightarrow1\quad\text{in probability as $n\to\infty$.}
\end{equation}
To prove this, we will exploit the following property, which we will refer to as the ``branching property'' of the random variables $(A_n^p(u))_{u\in[\mathbb{U}]_{n-1}}$: for each ${k\in\llbracket0,n-1\rrbracket}$, the random variables $(A_n^p(u))_{u\in\partial[\mathbb{U}]_k}$ are independent, have the same distribution as $A_{n-k}^p(\varnothing)$, and are independent of the random variables $(\xi_n^p(u))_{u\in[\mathbb{U}]_{k-1}}$.
We first exploit this branching property to establish the following rough estimate:

\begin{lem}\label{lem:estimateanp}
There exists a constant $C\in\mathbb{R}_+^*$ such that
\[\mathbb{E}\left[A_n^p(\varnothing)^2\right]\leq C\cdot\left(\frac{1\vee(m^n\cdot p)}{n}\right)^2\quad\text{for all $n\in\mathbb{N}^*$ and $p\in{]0,1]}$.}\]
\end{lem}
\begin{proof}
We prove this in several steps.

\quad

\emph{Step 1.} For all $n\in\mathbb{N}^*$ and $p\in{]0,1]}$, set $a_n^p=\mathbb{E}\left[A_n^p(\varnothing)\right]$, and let us show that
\begin{equation}\label{eq:anp}
a_n^p\leq\frac{1}{1-G^n(1-p)}\cdot\frac{m^n\cdot p}{n}.
\end{equation}
For $n=1$, we have
\[a_1^p=\mathbb{E}\left[\xi_1^p(\varnothing)\right]=\frac{m\cdot p}{1-G(1-p)}\]
by \eqref{eq:meanxinp}, hence \eqref{eq:anp} holds.
Next, fix $n\geq2$, and let us prove \eqref{eq:anp}.
By the branching property of the random variables $(A_n^p(u))_{u\in[\mathbb{U}]_{n-1}}$, and by \eqref{eq:meanxinp}, we have
\begin{eqnarray}
a_n^p&=&\mathbb{E}\left[\xi_n^p(\varnothing)\right]\cdot\mathbb{E}\left[\frac{A_{n-1}^p(\varnothing)}{1+A_{n-1}^p(\varnothing)}\right]\label{eq:anpbranching}\\
&=&\frac{m\cdot\left(1-G^{n-1}(1-p)\right)}{1-G^n(1-p)}\cdot\mathbb{E}\left[\frac{A_{n-1}^p(\varnothing)}{1+A_{n-1}^p(\varnothing)}\right].\notag
\end{eqnarray}
Then, since $x\in\mathbb{R}_+\mapsto x/(1+x)$ is concave, by Jensen's inequality, we get
\[a_n^p\leq\frac{m\cdot\left(1-G^{n-1}(1-p)\right)}{1-G^n(1-p)}\cdot\frac{a_{n-1}^p}{1+a_{n-1}^p}.\]
We rearrange this as follows:
\[\frac{m^{n-1}}{1-G^{n-1}(1-p)}\cdot\left(\frac{1}{a_{n-1}^p}+1\right)\leq\frac{m^n}{1-G^n(1-p)}\cdot\frac{1}{a_n^p}.\]
Telescoping the above inequality, we obtain
\[\begin{split}
\frac{m^n}{1-G^n(1-p)}\cdot\frac{1}{a_n^p}&\geq\frac{m}{1-G(1-p)}\cdot\frac{1}{a_1^p}+\sum_{k=2}^n\frac{m^{k-1}}{1-G^{k-1}(1-p)}\\
&=\frac{m}{1-G(1-p)}\cdot\frac{1-G(1-p)}{m\cdot p}+\sum_{l=1}^{n-1}\frac{m^l}{1-G^l(1-p)}\\
&=\sum_{k=0}^{n-1}\frac{m^k}{1-G^k(1-p)}.
\end{split}\]
Finally, since $1-G^k(1-p)\leq m^k\cdot p$ for all $k\in\mathbb{N}$ (see the proof of Lemma \ref{lem:asymptoticsurvival}), we get
\[\frac{m^n}{1-G^n(1-p)}\cdot\frac{1}{a_n^p}\geq\sum_{k=0}^{n-1}\frac{m^k}{m^k\cdot p}=\frac{n}{p},\]
which yields \eqref{eq:anp}.

\quad

\emph{Step 2.} For all $n\in\mathbb{N}^*$ and $p\in{]0,1]}$, set $b_n^p=\mathbb{V}\left(A_n^p(\varnothing)\right)=\mathbb{E}\left[A_n^p(\varnothing)^2\right]-\left(a_n^p\right)^2$, and let us show that
\begin{equation}\label{eq:bnp}
b_n^p\leq\frac{\mathbb{E}[\xi(\xi-1)]}{m^2}\cdot\frac{m^n}{1-G^n(1-p)}\cdot\sum_{k=1}^n\frac{m^k\cdot p^2}{k^2}.
\end{equation}
For $n=1$, we have
\[b_1^p=\mathbb{V}\left(\xi_1^p(\varnothing)\right)\leq\frac{\mathbb{E}[\xi(\xi-1)]\cdot p^2}{1-G(1-p)}\]
by \eqref{eq:varxinp}, hence \eqref{eq:bnp} holds.
Next, fix $n\geq2$, and let us prove \eqref{eq:bnp}.
By the branching property of the random variables $(A_n^p(u))_{u\in[\mathbb{U}]_{n-1}}$, and by the ``law of total variance'', we have
\begin{eqnarray*}
b_n^p&=&\mathbb{E}\left[\xi_n^p(\varnothing)\cdot\mathbb{V}\left(\frac{A_{n-1}^p(\varnothing)}{1+A_{n-1}^p(\varnothing)}\right)\right]+\mathbb{V}\left(\xi_n^p(\varnothing)\cdot\mathbb{E}\left[\frac{A_{n-1}^p(\varnothing)}{1+A_{n-1}^p(\varnothing)}\right]\right)\\
&=&\mathbb{E}\left[\xi_n^p(\varnothing)\right]\cdot\mathbb{V}\left(\frac{A_{n-1}^p(\varnothing)}{1+A_{n-1}^p(\varnothing)}\right)+\mathbb{V}\left(\xi_n^p(\varnothing)\right)\cdot\mathbb{E}\left[\frac{A_{n-1}^p(\varnothing)}{1+A_{n-1}^p(\varnothing)}\right]^2.
\end{eqnarray*}
Then, since the function $\varphi:x\in\mathbb{R}_+\mapsto x/(1+x)$ is $1$-Lipschitz, we have
\[\begin{split}
\mathbb{V}\left(\frac{A_{n-1}^p(\varnothing)}{1+A_{n-1}^p(\varnothing)}\right)&=\mathbb{V}\left(\varphi\left(A_{n-1}^p(\varnothing)\right)\right)\\
&\leq\mathbb{E}\left[\left(\varphi\left(A_{n-1}^p(\varnothing)\right)-\varphi\left(a_{n-1}^p\right)\right)^2\right]\\
&\leq\mathbb{E}\left[\left(A_{n-1}^p(\varnothing)-a_{n-1}^p\right)^2\right]=b_{n-1}^p.
\end{split}\]
Moreover, by \eqref{eq:anpbranching}, we have
\[\mathbb{E}\left[\frac{A_{n-1}^p(\varnothing)}{1+A_{n-1}^p(\varnothing)}\right]=\frac{1-G^n(1-p)}{m\cdot\left(1-G^{n-1}(1-p)\right)}\cdot a_n^p.\]
Thus, by \eqref{eq:meanxinp} and by \eqref{eq:varxinp}, we get
\begin{eqnarray*}
b_n^p&\leq&\frac{m\cdot\left(1-G^{n-1}(1-p)\right)}{1-G^n(1-p)}\cdot b_{n-1}^p\\
&&+\frac{\mathbb{E}[\xi(\xi-1)]\cdot\left(1-G^{n-1}(1-p)\right)^2}{1-G^n(1-p)}\cdot\left(\frac{1-G^n(1-p)}{m\cdot\left(1-G^{n-1}(1-p)\right)}\cdot a_n^p\right)^2\\
&=&\frac{m\cdot\left(1-G^{n-1}(1-p)\right)}{1-G^n(1-p)}\cdot b_{n-1}^p+\frac{\mathbb{E}[\xi(\xi-1)]}{m^2}\cdot\left(1-G^n(1-p)\right)\cdot\left(a_n^p\right)^2.
\end{eqnarray*}
Plugging \eqref{eq:anp} into the previous inequality, we obtain
\[b_n^p\leq\frac{m\cdot\left(1-G^{n-1}(1-p)\right)}{1-G^n(1-p)}\cdot b_{n-1}^p+\frac{\mathbb{E}[\xi(\xi-1)]}{m^2}\cdot\frac{1}{1-G^n(1-p)}\cdot\left(\frac{m^n\cdot p}{n}\right)^2.\]
We rearrange this as follows:
\[\frac{1-G^n(1-p)}{m^n}\cdot b_n^p\leq\frac{1-G^{n-1}(1-p)}{m^{n-1}}\cdot b_{n-1}^p+\frac{\mathbb{E}[\xi(\xi-1)]}{m^2}\cdot\frac{m^n\cdot p^2}{n^2}.\]
Telescoping the above inequality, we obtain
\begin{eqnarray*}
\frac{1-G^n(1-p)}{m^n}\cdot b_n^p&\leq&\frac{1-G(1-p)}{m}\cdot b_1^p+\frac{\mathbb{E}[\xi(\xi-1)]}{m^2}\cdot\sum_{k=2}^n\frac{m^k\cdot p^2}{k^2}\\
&\leq&\frac{1-G(1-p)}{m}\cdot\frac{\mathbb{E}[\xi(\xi-1)]\cdot p^2}{1-G(1-p)}+\frac{\mathbb{E}[\xi(\xi-1)]}{m^2}\cdot\sum_{k=2}^n\frac{m^k\cdot p^2}{k^2}\\
&=&\frac{\mathbb{E}[\xi(\xi-1)]}{m^2}\cdot\sum_{k=1}^n\frac{m^k\cdot p^2}{k^2},
\end{eqnarray*}
which yields \eqref{eq:bnp}.

\quad

\emph{Step 3.} Observe that there exists a constant $C\in\mathbb{R}_+^*$ such that\footnote{Indeed, we have
\[\frac{n^2}{m^n}\cdot\sum_{k=1}^n\frac{m^k}{k^2}=\sum_{k=1}^n\frac{m^{-(n-k)}}{(k/n)^2}=\sum_{l=0}^{n-1}\frac{m^{-l}}{((n-l)/n)^2}\leq\sum_{l=0}^{\lfloor n/2\rfloor}\frac{m^{-l}}{(1/2)^2}+\sum_{l=\lfloor n/2\rfloor+1}^{n-1}\frac{m^{-n/2}}{(1/n)^2}\leq4\cdot\sum_{l\geq0}m^{-l}+n^3\cdot m^{-n/2}.\]}
\[\sum_{k=1}^n\frac{m^k}{k^2}\leq C\cdot\frac{m^n}{n^2}\quad\text{for all $n\in\mathbb{N}^*$.}\]
Thus, by \eqref{eq:bnp}, we obtain
\[b_n^p\leq\frac{\mathbb{E}[\xi(\xi-1)]}{m^2}\cdot\frac{m^n}{1-G^n(1-p)}\cdot C\cdot\frac{m^n}{n^2}\cdot p^2=\frac{C\cdot\mathbb{E}[\xi(\xi-1)]}{m^2\cdot\left(1-G^n(1-p)\right)}\cdot\left(\frac{m^n\cdot p}{n}\right)^2.\]
Since $\mathbb{E}\left[A_n^p(\varnothing)^2\right]=b_n^p+\left(a_n^p\right)^2$, by \eqref{eq:anp}, we get (since $1-G^n(1-p)\leq1$):
\begin{eqnarray*}
\mathbb{E}\left[A_n^p(\varnothing)^2\right]&\leq&\frac{C\cdot\mathbb{E}[\xi(\xi-1)]}{m^2\cdot\left(1-G^n(1-p)\right)}\cdot\left(\frac{m^n\cdot p}{n}\right)^2+\left(\frac{1}{1-G^n(1-p)}\cdot\frac{m^n\cdot p}{n}\right)^2\\
&\leq&\left(\frac{C\cdot\mathbb{E}[\xi(\xi-1)]}{m^2}+1\right)\cdot\frac{1}{\left(1-G^n(1-p)\right)^2}\cdot\left(\frac{m^n\cdot p}{n}\right)^2.
\end{eqnarray*}
By Lemma \ref{lem:asymptoticsurvival}, it follows that there exists a constant $c\in{]0,1[}$ such that
\begin{eqnarray*}
\mathbb{E}\left[A_n^p(\varnothing)^2\right]&\leq&\left(\frac{C\cdot\mathbb{E}[\xi(\xi-1)]}{m^2}+1\right)\cdot\frac{1}{\left(c\cdot((m^n\cdot p)\wedge1)\right)^2}\cdot\left(\frac{m^n\cdot p}{n}\right)^2\\
&=&\left(\frac{C\cdot\mathbb{E}[\xi(\xi-1)]}{m^2}+1\right)\cdot\frac{1}{c^2}\cdot\left(\frac{1\vee(m^n\cdot p)}{n}\right)^2,
\end{eqnarray*}
which concludes the proof.
\end{proof}

We are now ready to prove the more precise estimate \eqref{eq:vardgffpruned}:

\begin{proof}[Proof of the quantitative delocalisation in Theorem \ref{thm:dgff}]
For all $u\in\mathbf{T}_n^{p_n}$ such that $|u|\leq n-2$, we have
\[A_n^{p_n}(u)=\sum_{\substack{v\in\mathbf{T}_n^{p_n}\\\text{$v$ child of $u$}}}\frac{A_n^{p_n}(v)}{1+A_n^{p_n}(v)}=\sum_{\substack{v\in\mathbf{T}_n^{p_n}\\\text{$v$ child of $u$}}}A_n^{p_n}(v)-\sum_{\substack{v\in\mathbf{T}_n^{p_n}\\\text{$v$ child of $u$}}}\frac{A_n^{p_n}(v)^2}{1+A_n^{p_n}(v)}.\]
Thus, for all $\ell\in\llbracket0,n-1\rrbracket$, we have
\begin{equation}\label{eq:decomp}
A_n^{p_n}(\varnothing)=\sum_{\substack{u\in\mathbf{T}_n^{p_n}\\|u|=\ell}}A_n^{p_n}(u)-\sum_{k=1}^\ell\sum_{\substack{u\in\mathbf{T}_n^{p_n}\\|u|=k}}\frac{A_n^{p_n}(u)^2}{1+A_n^{p_n}(u)}.
\end{equation}
Now, let $K_n$ be the smallest integer $k\in\llbracket0,n\rrbracket$ such that every node $u\in\mathbf{T}_n^{p_n}$ with ${k\leq|u|\leq n-1}$ has exactly one child.\footnote{Note that this property holds vacuously for $k=n$.}
By the same calculation as in Remark \ref{rem:vardgfftrees}, we have
\[A_n^{p_n}(u)=\frac{1}{n-|u|}\quad\text{for all $u\in\mathbf{T}_n^{p_n}$ such that $K_n\leq|u|\leq n-1$.}\]
Thus, on the good event $(K_n\leq n-1)$, by \eqref{eq:decomp}, we have
\[\frac{n}{\#\partial\mathbf{T}_n^{p_n}}\cdot A_n^{p_n}(\varnothing)=\frac{n}{n-K_n}-\frac{n}{\#\partial\mathbf{T}_n^{p_n}}\cdot\sum_{k=1}^{K_n}\sum_{\substack{u\in\mathbf{T}_n^{p_n}\\|u|=k}}\frac{A_n(u)^2}{1+A_n(u)}.\]
Indeed, by the definition of $K_n$, we have $\#\partial[\mathbf{T}_n^{p_n}]_{K_n}=\#\partial\mathbf{T}_n^{p_n}$.
Therefore, to prove \eqref{eq:vardgffpruned}, it suffices to show that
\begin{equation}\label{eq:Kn}
\frac{K_n}{n}\longrightarrow0\quad\text{in probability as $n\to\infty$,}
\end{equation}
and
\begin{equation}\label{eq:remainder}
\frac{n}{\#\partial\mathbf{T}_n^{p_n}}\cdot\sum_{k=1}^{K_n}\sum_{\substack{u\in\mathbf{T}_n^{p_n}\\|u|=k}}\frac{A_n^{p_n}(u)^2}{1+A_n^{p_n}(u)}\longrightarrow0\quad\text{in probability as $n\to\infty$.}
\end{equation}
We start with \eqref{eq:Kn}.
Set $\ell_n=\lceil2\log_m(n)\rceil$.
By Proposition \ref{prop:pruned} and by Lemma \ref{lem:cut-height}, there exists a constant $C\in\mathbb{R}_+^*$ such that
\[\mathbb{P}(K_n>\ell_n)\leq C\cdot\frac{(m^n\cdot p_n)\vee(m^n\cdot p_n)^2}{m^{\ell_n}}\leq C\cdot\frac{(m^n\cdot p_n)\vee(m^n\cdot p_n)^2}{n^2}.\]
Thus, by \eqref{eq:assdelocdgff}, we get that $\mathbb{P}(K_n>\ell_n)\rightarrow0$ as $n\to\infty$.
This proves \eqref{eq:Kn}, since $\ell_n/n\rightarrow0$ as $n\to\infty$.
Finally, let us prove \eqref{eq:remainder}.
Since $\mathbb{P}(K_n>\ell_n)\rightarrow0$ as $n\to\infty$, it suffices to show that
\[\frac{n}{\#\partial\mathbf{T}_n^{p_n}}\cdot\sum_{k=1}^{\ell_n}\sum_{\substack{u\in\mathbf{T}_n^{p_n}\\|u|=k}}\frac{A_n^{p_n}(u)^2}{1+A_n^{p_n}(u)}\longrightarrow0\quad\text{in probability as $n\to\infty$.}\]
Now, by Proposition \ref{prop:pruned} and by Lemma \ref{lem:leavesTnp}, the sequence of random variables
\[\left(\frac{\#\partial\mathbf{T}_n^{p_n}}{1\vee(m^n\cdot p_n)}\right)_{n\in\mathbb{N}^*}\]
is tight in $\mathbb{R}_+^*$.
Therefore, it suffices to show that
\[\frac{n}{1\vee(m^n\cdot p_n)}\cdot\sum_{k=1}^{\ell_n}\sum_{\substack{u\in\mathbf{T}_n^{p_n}\\|u|=k}}\frac{A_n^{p_n}(u)^2}{1+A_n^{p_n}(u)}\longrightarrow0\quad\text{in probability as $n\to\infty$.}\]
To this end, let us show that
\begin{equation}\label{eq:error0}
\frac{n}{1\vee(m^n\cdot p_n)}\cdot\mathbb{E}\left[\sum_{k=1}^{\ell_n}\sum_{\substack{u\in\mathbf{T}_n^{p_n}\\|u|=k}}A_n^{p_n}(u)^2\right]\longrightarrow0\quad\text{as $n\to\infty$.}
\end{equation}
For all $n\in\mathbb{N}^*$ sufficiently large so that $\ell_n\leq n-1$, the following holds: for each $k\in\llbracket1,\ell_n\rrbracket$, by the branching property of the random variables $\left(A_n^{p_n}(u)\right)_{u\in[\mathbb{U}]_{n-1}}$, by Proposition \ref{prop:pruned} and by \eqref{eq:Tnpgenk}, we have
\begin{eqnarray*}
\mathbb{E}\left[\sum_{\substack{u\in\mathbf{T}_n^{p_n}\\|u|=k}}A_n^{p_n}(u)^2\right]&=&\mathbb{E}\left[\#\partial\left[\mathbf{T}_n^{p_n}\right]_k\right]\cdot\mathbb{E}\left[A_{n-k}^{p_n}(\varnothing)^2\right]\\
&=&\frac{m^k\cdot\left(1-G^{n-k}(1-p_n)\right)}{1-G^n(1-p_n)}\cdot\mathbb{E}\left[A_{n-k}^{p_n}(\varnothing)^2\right].
\end{eqnarray*}
Thus, by Lemma \ref{lem:asymptoticsurvival} and by Lemma \ref{lem:estimateanp}, there exists constants $c\in{]0,1[}$ and $C'\in\mathbb{R}_+^*$ such that
\begin{eqnarray*}
\mathbb{E}\left[\sum_{\substack{u\in\mathbf{T}_n^{p_n}\\|u|=k}}A_n^{p_n}(u)^2\right]&\leq&\frac{m^k\cdot\left(\left(m^{n-k}\cdot p_n\right)\wedge1\right)}{c\cdot(m^n\cdot p_n)\wedge1)}\cdot C'\cdot\left(\frac{1\vee\left(m^{n-k}\cdot p_n\right)}{n-k}\right)^2\\
&=&\frac{C'}{c}\cdot\frac{m^n\cdot p_n}{(m^n\cdot p_n)\wedge1}\cdot\frac{1\vee\left(m^{n-k}\cdot p_n\right)}{(n-k)^2}.
\end{eqnarray*}
It follows that
\[\sum_{k=1}^{\ell_n}\mathbb{E}\left[\sum_{\substack{u\in\mathbf{T}_n^{p_n}\\|u|=k}}A_n^{p_n}(u)^2\right]\leq\frac{C'}{c}\cdot(1\vee(m^n\cdot p_n))\cdot\sum_{k=1}^{\ell_n}\frac{1\vee\left(m^{n-k}\cdot p_n\right)}{(n-k)^2},\]
hence
\begin{eqnarray*}
\frac{n}{1\vee(m^n\cdot p_n)}\cdot\sum_{k=1}^{\ell_n}\mathbb{E}\left[\sum_{\substack{u\in\mathbf{T}_n^{p_n}\\|u|=k}}A_n^{p_n}(u)^2\right]&\leq&\frac{C'}{c}\cdot n\cdot\sum_{k=1}^{\ell_n}\frac{1\vee\left(m^{n-k}\cdot p_n\right)}{(n-k)^2}.
\end{eqnarray*}
Finally, we have
\[\sum_{k=1}^{\ell_n}\frac{1\vee\left(m^{n-k}\cdot p_n\right)}{(n-k)^2}\leq\sum_{k=1}^{\ell_n}\frac{1+m^{n-k}\cdot p_n}{(n-\ell_n)^2}\leq\frac{1}{(n-\ell_n)^2}\cdot\left(\ell_n+\frac{m^n\cdot p_n}{m-1}\right).\]
Altogether, we obtain that
\begin{eqnarray*}
\frac{n}{1\vee(m^n\cdot p_n)}\cdot\sum_{k=1}^{\ell_n}\mathbb{E}\left[\sum_{\substack{u\in\mathbf{T}_n^{p_n}\\|u|=k}}A_n^{p_n}(u)^2\right]&\leq&\frac{C'}{c}\cdot \frac{n}{(n-\ell_n)^2}\cdot\left(\ell_n+\frac{m^n\cdot p_n}{m-1}\right).
\end{eqnarray*}
Since $\ell_n/n\rightarrow0$ as $n\to\infty$, and by virtue of \eqref{eq:assdelocdgff}, the above right-hand side vanishes as ${n\to\infty}$, which proves \eqref{eq:error0}.
The proof of (the quantitative delocalisation in) Theorem \ref{thm:dgff} is now complete.
\end{proof}

\section{Quenched concentration results for supercritical branching random walks}\label{sec:concentration}

In this interlude, which can be read independently of the rest of the paper, we prove quenched concentration results for supercritical branching random walks, which will be instrumental in the proof of Theorem \ref{thm:main}.
For each $n\in\mathbb{N}^*$, fix a finite rooted tree $T_n$ with height $n$ and with no internal leaves (i.e, such that each node $u\in T_n$ with $|u|\leq n-1$ has at least one child).
Then, let $(X_u)_{u\in T_n^*}$ be a collection of independent random variables with distribution $\nu$, where $T_n^*=T_n\setminus\{\varnothing\}$, and set
\[S_v=\sum_{\substack{u\in T_n^*\\u\preceq v}}X_u\quad\text{for all $v\in T_n$.}\]
Here, we do not assume that $\nu$ is discrete: we assume that $\nu$ is a Borel probability measure on $\mathbb{R}$, which is symmetric and has compact support.
Let ${\sigma^2=\int_{\mathbb{R}}x^2\cdot\nu(\mathrm{d}x)}$ be the variance of $\nu$, and set $L=\inf\{\ell\in\mathbb{R}_+:\nu([-\ell,\ell])=1\}$.
We assume that $\nu(\{0\})<1$, so that $\sigma^2,L>0$ (otherwise, the statement of Theorem \ref{thm:concentration2} is trivial, understanding $\delta_n$ as $\infty$).
For each $n\in\mathbb{N}^*$, we set $N_n=\#\partial T_n$, and we let
\[\langle S\rangle_n=\frac{1}{N_n}\cdot\sum_{v\in\partial T_n}S_v\quad\text{and}\quad\left\langle S^2\right\rangle_n=\frac{1}{N_n}\cdot\sum_{v\in\partial T_n}S_v^2.\]
We want to derive concentration bounds for the random variables $\langle S\rangle_n$ and $\left\langle S^2\right\rangle_n$.
To begin with, note that $\mathbb{E}[\langle S\rangle_n]=0$ and $\mathbb{E}\left[\left\langle S^2\right\rangle_n\right]=\sigma^2\cdot n$.
For simplicity, we will express all our bounds in terms of the following single tree-dependent quantity:
\[\gamma_n:=1\vee\left(\frac{1}{N_n}\cdot\sum_{k=1}^nk^3\cdot\max_{\substack{u\in T_n\\|u|=k}}N_n(u)\right),\]
where $N_n(u)=\#\{v\in\partial T_n:v\succeq u\}$ for all $u\in T_n^*$.
We start with the concentration of $\langle S\rangle_n$, which is straightforward.
It turns out that $\langle S\rangle_n$ is $\left(L^2\gamma_n\right)$-sub-Gaussian:

\begin{prop}\label{prop:concentration1}
For all $\lambda\in\mathbb{R}$, we have
\[\mathbb{E}\left[e^{\lambda\langle S\rangle_n}\right]\leq\exp\left(\frac{L^2\gamma_n\lambda^2}{2}\right).\]
In particular, we get that
\[\mathbb{E}\left[\exp\left(\lambda\langle S\rangle_n^2\right)\right]\leq\exp\left(8L^2\gamma_n\lambda\right)\quad\text{for all $\lambda\in\left[0,\frac{1}{4L^2\gamma_n}\right]$.}\]
\end{prop}
\begin{proof}
We have
\[\langle S\rangle_n=\frac{1}{N_n}\cdot\sum_{v\in\partial T_n}\sum_{\substack{u\in T_n^*\\u\preceq v}}X_u=\frac{1}{N_n}\cdot\sum_{u\in T_n^*}\sum_{\substack{v\in\partial T_n\\v\succeq u}}X_u=\frac{1}{N_n}\cdot\sum_{u\in T_n^*}X_u\cdot N_n(u).\]
Thus, for all $\lambda\in\mathbb{R}$, we have
\[\mathbb{E}\left[e^{\lambda\langle S\rangle_n}\right]=\mathbb{E}\left[\prod_{u\in T_n^*}\exp\left(\lambda\cdot\frac{N_n(u)}{N_n}\cdot X_u\right)\right]=\prod_{u\in T_n^*}\mathbb{E}\left[\exp\left(\lambda\cdot\frac{N_n(u)}{N_n}\cdot X_u\right)\right].\]
Now, by Hoeffding's lemma (see \cite[Lemma 2.2]{concentration}), the distribution $\nu$ is $\left(L^2\right)$-sub-Gaussian: we get
\begin{eqnarray*}
\mathbb{E}\left[e^{\lambda\langle S\rangle_n}\right]&\leq&\prod_{u\in T_n^*}\exp\left(\frac{L^2}{2}\cdot\left(\lambda\cdot\frac{N_n(u)}{N_n}\right)^2\right)\\
&=&\exp\left(\frac{L^2\lambda^2}{2}\cdot\sum_{u\in T_n^*}\left(\frac{N_n(u)}{N_n}\right)^2\right)\leq\exp\left(\frac{L^2\lambda^2}{2}\cdot\gamma_n\right).
\end{eqnarray*}
Indeed, we have
\begin{eqnarray}
\sum_{u\in T_n^*}\left(\frac{N_n(u)}{N_n}\right)^2&=&\frac{1}{N_n^2}\cdot\sum_{k=1}^n\sum_{\substack{u\in T_n\\|u|=k}}N_n(u)^2\notag\\
&\leq&\frac{1}{N_n^2}\cdot\sum_{k=1}^n\max_{\substack{u\in T_n\\|u|=k}}N_n(u)\cdot\sum_{\substack{u\in T_n\\|u|=k}}N_n(u)\notag\\
&=&\frac{1}{N_n}\cdot\sum_{k=1}^n\max_{\substack{u\in T_n\\|u|=k}}N_n(u)\notag\\
&\leq&\frac{1}{N_n}\cdot\sum_{k=1}^nk^3\cdot\max_{\substack{u\in T_n\\|u|=k}}N_n(u)=\gamma_n.\label{eq:boundgamman}
\end{eqnarray}
By the Chernoff bound, we deduce that
\[\mathbb{P}(|\langle S\rangle_n|\geq y)\leq2\exp\left(-\frac{y^2}{2L^2\gamma_n}\right)\quad\text{for all $y\in\mathbb{R}_+$.}\]
It follows that for all $\lambda\in\left[0,\left.1\middle/\left(2L^2\gamma_n\right)\right.\right[$, we have
\begin{eqnarray*}
\mathbb{E}\left[\exp\left(\lambda\langle S\rangle_n^2\right)\right]&=&1+\int_0^\infty\lambda e^{\lambda x}\cdot\mathbb{P}\left(\langle S\rangle_n^2\geq x\right)\mathrm{d}x\\
&\leq&1+\int_0^\infty\lambda e^{\lambda x}\cdot2\exp\left(-\frac{x}{2L^2\gamma_n}\right)\mathrm{d}x\\
&=&1+\frac{2\lambda}{\left.1\middle/\left(2L^2\gamma_n\right)\right.-\lambda}=1+\frac{4L^2\gamma_n\lambda}{1-2L^2\gamma_n\lambda}.
\end{eqnarray*}
Thus, for all $\lambda\in\left[0,\left.1\middle/\left(4L^2\gamma_n\right)\right.\right]$, we get
\[\mathbb{E}\left[\exp\left(\lambda\langle S\rangle_n^2\right)\right]\leq1+\frac{4L^2\gamma_n\lambda}{1/2}\leq\exp\left(8L^2\gamma_n\lambda\right).\]
\end{proof}

Next, we turn to the concentration of $\left\langle S^2\right\rangle_n$, which will occupy us for the rest of this section.

\subsection{The entropy method}

The proof of Theorem \ref{thm:concentration2} is based on the so-called ``entropy method''.
We refer the reader to the authoritative reference \cite{concentration} for a detailed account.
To implement the entropy method, let $X'=\left(X'_u\right)_{u\in T_n^*}$ be an independent copy of $X=(X_u)_{u\in T_n^*}$.
For each ${u\in T_n^*}$, we denote by $\overline{X}^u$ the random vector with values in $\mathbb{R}^{T_n^*}$ given by
\[\overline{X}^u_v=\begin{cases}
X'_v&\text{if $v=u$,}\\
X_v&\text{otherwise,}
\end{cases}\quad\text{for all $v\in T_n^*$.}\]
The first ingredient of the entropy method is a logarithmic Sobolev inequality.
In our case, it takes the following form: for any function $f:\mathbb{R}^{T_n^*}\rightarrow\mathbb{R}$, we have\footnote{This inequality is easily deduced from \cite[Theorem 6.15]{concentration} by using the inequality ${\tau(-z)\leq z^2}$ for all $z\in\mathbb{R}_+$.}
\begin{equation}\label{eq:logsobolev}
\mathrm{Ent}\left(e^{\lambda f(X)}\right)\leq\lambda^2\cdot\mathbb{E}\left[e^{\lambda f(X)}\cdot\sum_{u\in T_n^*}\left(f(X)-f\left(\overline{X}^u\right)\right)^2\right]\quad\text{for all $\lambda\in\mathbb{R}$.}
\end{equation}
The second ingredient of the entropy method is an instance of Herbst's argument.
For us, it will be given by the following lemma:

\begin{lem}\label{lem:herbst}
Let $Y$ be a real random variable such that $\mathbb{E}\left[e^{\lambda Y}\right]<\infty$ for all $\lambda\in\mathbb{R}$, let $Z$ be a non-negative random variable such that $\mathbb{E}\left[e^{\lambda Z}\right]<\infty$ for some $\lambda\in\mathbb{R}_+^*$, and suppose that
\[\mathrm{Ent}\left(e^{\lambda Y}\right)\leq\lambda^2\cdot\mathbb{E}\left[e^{\lambda Y}\cdot Z\right]\quad\text{for all $\lambda\in\mathbb{R}$.}\]
Then, fix $\lambda_0\in\mathbb{R}_+^*$ such that $\mathbb{E}\left[e^{\lambda_0Z}\right]\leq e$.
We have
\[\mathbb{E}\left[e^{\lambda(Y-\mathbb{E}[Y])}\right]\leq\exp\left(\frac{\lambda^2}{\delta^2}\right)\quad\text{for all $\lambda\in[-\delta,\delta]$,}\]
where $\delta=\sqrt{\lambda_0/2}$.
\end{lem}

The proof of this lemma is similar to \cite[proof of Theorem 6.16]{concentration}.

\begin{proof}[Proof of Lemma \ref{lem:herbst}]
Set
\[\psi(\lambda)=\ln\left(\mathbb{E}\left[e^{\lambda(Y-\mathbb{E}[Y])}\right]\right)=\ln\left(\mathbb{E}\left[e^{\lambda Y}\right]\right)-\lambda\mathbb{E}[Y]\quad\text{for all $\lambda\in\mathbb{R}$,}\]
and notice that for all $\lambda\in\mathbb{R}$, we have
\[\begin{split}
\lambda\cdot\psi'(\lambda)-\psi(\lambda)&=\lambda\cdot\frac{\mathbb{E}\left[Ye^{\lambda Y}\right]}{\mathbb{E}\left[e^{\lambda Y}\right]}-\ln\left(\mathbb{E}\left[e^{\lambda Y}\right]\right)\\
&=\frac{\mathbb{E}\left[e^{\lambda Y}\cdot\ln\left(e^{\lambda Y}\right)\right]-\mathbb{E}\left[e^{\lambda Y}\right]\cdot\ln\left(\mathbb{E}\left[e^{\lambda Y}\right]\right)}{\mathbb{E}\left[e^{\lambda Y}\right]}=\frac{\mathrm{Ent}\left(e^{\lambda Y}\right)}{\mathbb{E}\left[e^{\lambda Y}\right]}.
\end{split}\]
Dividing through by $\lambda^2$, we get that
\[\frac{\psi'(\lambda)}{\lambda}-\frac{\psi(\lambda)}{\lambda^2}\leq\frac{1}{\lambda^2}\cdot\frac{\mathrm{Ent}\left(e^{\lambda Y}\right)}{\mathbb{E}\left[e^{\lambda Y}\right]}\quad\text{for all $\lambda\in\mathbb{R}^*$,}\]
where in the left-hand side, we recognise the derivative of $\lambda\in\mathbb{R^*}\mapsto\psi(\lambda)/\lambda$.
Now, let us upper bound the right-hand side.
First, by assumption, we have
\[\mathrm{Ent}\left(e^{\lambda Y}\right)\leq\lambda^2\cdot\mathbb{E}\left[e^{\lambda Y}\cdot Z\right]=\frac{\lambda^2}{\lambda_0}\cdot\mathbb{E}\left[e^{\lambda Y}\cdot(\lambda_0Z)\right]\]
for all $\lambda\in\mathbb{R}$.
Then, by the duality formula for entropy (see \cite[Theorem 4.13 and Remark 4.4]{concentration}), we have
\[\mathbb{E}\left[e^{\lambda Y}\cdot(\lambda_0Z)\right]\leq\mathbb{E}\left[e^{\lambda Y}\right]\cdot\ln\left(\mathbb{E}\left[e^{\lambda_0Z}\right]\right)+\mathrm{Ent}\left(e^{\lambda Y}\right)\leq\mathbb{E}\left[e^{\lambda Y}\right]+\mathrm{Ent}\left(e^{\lambda Y}\right).\]
Now, set $\delta=\sqrt{\lambda_0/2}$.
For all $\lambda\in[-\delta,\delta]$, we get
\[\mathrm{Ent}\left(e^{\lambda Y}\right)\leq\frac{\lambda^2}{\lambda_0}\cdot\left(\mathbb{E}\left[e^{\lambda Y}\right]+\mathrm{Ent}\left(e^{\lambda Y}\right)\right)\leq\frac{\lambda^2}{\lambda_0}\cdot\mathbb{E}\left[e^{\lambda Y}\right]+\frac{1}{2}\cdot\mathrm{Ent}\left(e^{\lambda Y}\right),\]
hence
\[\mathrm{Ent}\left(e^{\lambda Y}\right)\leq\frac{2\lambda^2}{\lambda_0}\cdot\mathbb{E}\left[e^{\lambda Y}\right]=\frac{\lambda^2}{\delta^2}\cdot\mathbb{E}\left[e^{\lambda Y}\right].\]
We deduce that
\begin{equation}\label{eq:herbst}
\frac{\psi'(\lambda)}{\lambda}-\frac{\psi(\lambda)}{\lambda^2}\leq\frac{1}{\delta^2}\quad\text{for all $\lambda\in[-\delta,\delta]\setminus\{0\}$.}
\end{equation}
Integrating \eqref{eq:herbst} between $0^+$ and $\lambda$, we get:\footnote{Since $\psi(0)=0$ and $\psi'(0)=0$, we have $\psi(\lambda)/\lambda\rightarrow0$ as $\lambda\to0$.}
\[\frac{\psi(\lambda)}{\lambda}\leq\frac{\lambda}{\delta^2}\quad\text{for all $\lambda\in{]0,\delta]}$,}\]
hence
\[\psi(\lambda)\leq\frac{\lambda^2}{\delta^2}\quad\text{for all $\lambda\in{]0,\delta]}$.}\]
Similarly, integrating \eqref{eq:herbst} between $\lambda$ and $0^-$, we get
\[-\frac{\psi(\lambda)}{\lambda}\leq-\frac{\lambda}{\delta^2}\quad\text{for all $\lambda\in{[-\delta,0[}$,}\]
hence (multiplying through by $-\lambda>0$):
\[\psi(\lambda)\leq\frac{\lambda^2}{\delta^2}\quad\text{for all $\lambda\in{[-\delta,0[}$.}\]
Altogether, we obtain that
\[\psi(\lambda)\leq\frac{\lambda^2}{\delta^2}\quad\text{for all $\lambda\in[-\delta,\delta]$,}\]
which yields the desired result.
\end{proof}

Now, with these two ingredients at hand, let us dive into the proof of Theorem \ref{thm:concentration2}.
A first use of the logarithmic Sobolev inequality \eqref{eq:logsobolev} yields the following inequality:

\begin{lem}\label{lem:In}
For all $\lambda\in\mathbb{R}$, we have
\begin{equation}\label{eq:In}
\mathrm{Ent}\left(\exp\left(\lambda\left\langle S^2\right\rangle_n\right)\right)\leq\lambda^2\cdot\mathbb{E}\left[\exp\left(\lambda\left\langle S^2\right\rangle_n\right)\cdot\left(32L^2Y_n+32L^4\gamma_n\right)\right],
\end{equation}
where
\[Y_n=\frac{1}{N_n^2}\cdot\sum_{v,w\in\partial T_n}S_vS_w\cdot|v\wedge w|.\]
\end{lem}

Let us defer the proof of this lemma (which is just a calculation) to continue unrolling the proof of Theorem \ref{thm:concentration2}.
With \eqref{eq:In} at hand, in view of Lemma \ref{lem:herbst}, we now want some control on $Y_n$.
Set
\begin{equation}\label{eq:defIn}
Z_n=\frac{1}{N_n^2}\cdot\sum_{v_1,v_2\in\partial T_n}S_{v_1}S_{v_2}\cdot\zeta_n(v_1,v_2),
\end{equation}
where
\[\zeta_n(v_1,v_2)=\frac{1}{N_n^2}\cdot\sum_{w_1,w_2\in\partial T_n}|v_1\wedge w_1|\cdot|v_2\wedge w_2|\cdot|w_1\wedge w_2|\quad\text{for all $v_1,v_2\in\partial T_n$.}\]
Another use of the logarithmic Sobolev inequality \eqref{eq:logsobolev} yields the following inequality:

\begin{lem}\label{lem:Jn}
For all $\lambda\in\mathbb{R}$, we have
\begin{equation}\label{eq:Jn}
\mathrm{Ent}\left(e^{\lambda Y_n}\right)\leq\lambda^2\cdot\mathbb{E}\left[e^{\lambda Y_n}\cdot\left(32L^2Z_n+32L^4\gamma_n\right)\right].
\end{equation}
\end{lem}

Again, let us defer the proof of this lemma (which is just a calculation) to continue unrolling the proof of Theorem \ref{thm:concentration2}.
With \eqref{eq:Jn} at hand, in view of Lemma \ref{lem:herbst}, we now want some control on $Z_n$.
Luckily, this entropy method procedure terminates, since we have the following direct control on $Z_n$:

\begin{lem}\label{lem:Jn'}
For all $\lambda\in\left[0,\left.1\middle/\left(4L^2\gamma_n^2\right)\right.\right]$, we have
\[\mathbb{E}\left[e^{\lambda Z_n}\right]\leq\exp\left(8L^2\gamma_n^2\lambda\right).\]
\end{lem}

Again, let us defer the proof of this lemma to finish unrolling the proof of Theorem \ref{thm:concentration2}:

\begin{proof}[Proof of Theorem \ref{thm:concentration2} assuming Lemma \ref{lem:In}, Lemma \ref{lem:Jn} and Lemma \ref{lem:Jn'}]
By Lemma \ref{lem:Jn'}, the following holds: for all $\lambda\in\left[0,\left.1\middle/\left(128L^4\gamma_n^2\right)\right.\right]$, we have
\begin{eqnarray*}
\mathbb{E}\left[\exp\left(\lambda\left(32L^2Z_n+32L^4\gamma_n\right)\right)\right]&=&\mathbb{E}\left[\exp\left(32L^2\lambda Z_n\right)\right]\cdot \exp\left(32L^4\gamma_n\lambda\right)\\
&\leq&\exp\left(8L^2\gamma_n^2\cdot32L^2\lambda\right)\cdot\exp\left(32L^4\gamma_n\lambda\right)\\
&=&\exp\left(128L^4\gamma_n^2\lambda+32L^4\gamma_n\lambda\right)\leq\exp\left(256L^4\gamma_n^2\lambda\right).
\end{eqnarray*}
Now, we set $\lambda_0=\left.1\middle/\left(256L^4\gamma_n^2\right)\right.$, so that $\mathbb{E}\left[\exp\left(\lambda_0\left(32L^2Z_n+32L^4\gamma_n\right)\right)\right]\leq e$.
In view of \eqref{eq:Jn}, by Lemma \ref{lem:herbst}, we get that
\[\mathbb{E}\left[e^{\lambda(Y_n-\mathbb{E}[Y_n])}\right]\leq\exp\left(\frac{\lambda^2}{\delta^2}\right)\quad\text{for all $\lambda\in[-\delta,\delta]$,}\]
where $\delta=\sqrt{\lambda_0/2}$.
Now, observe that
\[\mathbb{E}\left[Y_n\right]=\frac{1}{N_n^2}\cdot\sum_{v,w\in\partial T_n}\mathbb{E}[S_vS_w]\cdot|v\wedge w|=\frac{1}{N_n^2}\cdot\sum_{v,w\in\partial T_n}\sigma^2\cdot|v\wedge w|^2,\]
hence $\mathbb{E}[Y_n]\leq\sigma^2\cdot\gamma_n\leq L^2\cdot\gamma_n$.
Indeed, the variance of a random variable with values in $[-L,L]$ is at most $L^2$, and for each $v\in\partial T_n$, we have
\begin{eqnarray}
\frac{1}{N_n}\cdot\sum_{w\in\partial T_n}|v\wedge w|^2&=&\frac{1}{N_n}\cdot\sum_{\substack{u\in T_n\\u\preceq v}}\sum_{\substack{w\in\partial T_n\\v\wedge w=u}}|v\wedge w|^2\notag\\
&=&\frac{1}{N_n}\cdot\sum_{\substack{u\in T_n\\u\preceq v}}|u|^2\cdot\#\{w\in\partial T_n:v\wedge w=u\}\notag\\
&\leq&\frac{1}{N_n}\cdot\sum_{\substack{u\in T_n\\u\preceq v}}|u|^2\cdot\#\{w\in\partial T_n:w\succeq u\}\notag\\
&=&\frac{1}{N_n}\cdot\sum_{\substack{u\in T_n^*\\u\preceq v}}|u|^2\cdot N_n(u)\notag\\
&\leq&\frac{1}{N_n}\cdot\sum_{k=1}^nk^2\cdot\max_{\substack{u\in T_n\\|u|=k}}N_n(u)\notag\\
&\leq&\frac{1}{N_n}\cdot\sum_{k=1}^nk^3\cdot\max_{\substack{u\in T_n\\|u|=k}}N_n(u)=\gamma_n,\label{eq:calculation}
\end{eqnarray}
hence
\[\frac{1}{N_n^2}\cdot\sum_{v,w\in\partial T_n}|v\wedge w|^2=\frac{1}{N_n}\cdot\sum_{v\in\partial T_n}\left(\frac{1}{N_n}\cdot\sum_{w\in\partial T_n}|v\wedge w|^2\right)\leq\frac{1}{N_n}\cdot\sum_{v\in\partial T_n}\gamma_n=\gamma_n.\]
It follows that for all $\lambda\in[0,\delta]$, we have
\[\mathbb{E}\left[e^{\lambda Y_n}\right]=\mathbb{E}\left[e^{\lambda(Y_n-\mathbb{E}[Y_n])}\right]\cdot e^{\lambda\mathbb{E}[Y_n]}\leq\exp\left(\frac{\lambda^2}{\delta^2}\right)\cdot\exp\left(\lambda\cdot L^2\gamma_n\right)=\exp\left(\frac{\lambda^2}{\delta^2}+L^2\gamma_n\lambda\right).\]
Thus, for all $\lambda\in\left[0,\left.\delta\middle/\left(32L^2\right)\right.\right]$, we get
\begin{eqnarray*}
\mathbb{E}\left[\exp\left(\lambda\left(32L^2Y_n+32L^4\gamma_n\right)\right)\right]&=&\mathbb{E}\left[\exp\left(32L^2\lambda Y_n\right)\right]\cdot\exp\left(32L^4\gamma_n\lambda\right)\\
&\leq&\exp\left(\frac{\left(32L^2\lambda\right)^2}{\delta^2}+L^2\gamma_n\cdot32L^2\lambda\right)\cdot\exp\left(32L^4\gamma_n\lambda\right)\\
&\leq&\exp\left(\frac{\left(32L^2\lambda\right)^2}{\delta^2}+64L^4\gamma_n\lambda\right).
\end{eqnarray*}
Now, we set
\[\lambda_1=\frac{\left.\delta\middle/\left(32L^2\right)\right.}{\sqrt{2}}=\frac{\sqrt{\lambda_0}}{64L^2}=\frac{1}{1024L^4\gamma_n}\leq\frac{1}{128L^4\gamma_n},\]
so that
\begin{eqnarray*}
\mathbb{E}\left[\exp\left(\lambda_1\left(32L^2Y_n+32L^4\gamma_n\right)\right)\right]&\leq&\exp\left(\frac{\left(\delta\middle/\sqrt{2}\right)^2}{\delta^2}+64L^4\gamma_n\cdot\frac{1}{128L^4\gamma_n}\right)\\
&\leq&\exp\left(\frac{1}{2}+\frac{1}{2}\right)=e.
\end{eqnarray*}
By a final application of Lemma \ref{lem:herbst}, we obtain that
\[\mathbb{E}\left[\exp\left(\lambda\left(\left\langle S^2\right\rangle_n-\sigma^2n\right)\right)\right]\leq\exp\left(\frac{\lambda^2}{\delta_n^2}\right)\quad\text{for all $\lambda\in[-\delta_n,\delta_n]$,}\]
where
\[\delta_n=\sqrt{\frac{\lambda_1}{2}}=\sqrt{\frac{1}{2048L^4\gamma_n}}=\frac{1}{32L^2\sqrt{2\gamma_n}}.\]
This is the desired result.
\end{proof}

Let us conclude this subsection by stating, for future use, the following corollary of Proposition \ref{prop:concentration1} and Theorem \ref{thm:concentration2}:

\begin{cor}\label{cor:zeta}
For all $\lambda_1,\lambda_2\in\mathbb{R}$ such that $|\lambda_2|\leq\delta_n/2$, we have
\[1\leq\mathbb{E}\left[\exp\left(\lambda_1\langle S\rangle_n+\lambda_2\left(\left\langle S^2\right\rangle_n-\sigma^2n\right)\right)\right]\leq\exp\left(L^2\gamma_n\lambda_1^2+\frac{2\lambda_2^2}{\delta_n^2}\right).\]
\end{cor}
\begin{proof}
The lower bound is just Jensen's inequality, and in fact holds for all $\lambda_2\in\mathbb{R}$.
We now turn to the upper bound.
By the Cauchy--Schwarz inequality, we have
\begin{eqnarray*}
\lefteqn{\mathbb{E}\left[\exp\left(\lambda_1\langle S\rangle_n+\lambda_2\left(\left\langle S^2\right\rangle_n-\sigma^2n\right)\right)\right]}\\
&\leq&\mathbb{E}\left[\exp\left(2\lambda_1\langle S\rangle_n\right)\right]^{1/2}\cdot\mathbb{E}\left[\exp\left(2\lambda_2\left(\left\langle S^2\right\rangle_n-\sigma^2n\right)\right)\right]^{1/2}.
\end{eqnarray*}
Then, on the one hand, by Proposition \ref{prop:concentration1}, we have
\[\mathbb{E}\left[\exp\left(2\lambda_1\langle S\rangle_n\right)\right]^{1/2}\leq\exp\left(\frac{L^2\gamma_n\cdot(2\lambda_1)^2}{2}\right)^{1/2}=\exp\left(L^2\gamma_n\lambda_1^2\right).\]
On the other hand, by Theorem \ref{thm:concentration2}, since $|\lambda_2|\leq\delta_n/2$, we get
\[\mathbb{E}\left[\exp\left(2\lambda_2\left(\left\langle S^2\right\rangle_n-\sigma^2n\right)\right)\right]^{1/2}\leq\exp\left(\frac{(2\lambda_2)^2}{\delta_n^2}\right)^{1/2}=\exp\left(\frac{2\lambda_2^2}{\delta_n^2}\right).\]
This yields the desired upper bound.
\end{proof}

\subsection{Deferred calculations}

To conclude this section, let us present the calculations that we have deferred in the previous subsection, namely the proofs of Lemma \ref{lem:In}, Lemma \ref{lem:Jn} and Lemma \ref{lem:Jn'}.

\begin{proof}[Proof of Lemma \ref{lem:In}]
Consider the function ${f:\mathbb{R}^{T_n^*}\rightarrow\mathbb{R}_+}$ given by
\[f(x)=\frac{1}{N_n}\cdot\sum_{v\in\partial T_n}\left(\sum_{\substack{u\in T_n^*\\u\preceq v}}x_u\right)^2\quad\text{for all $x=(x_u)_{u\in T_n^*}\in\mathbb{R}^{T_n^*}$.}\]
By definition, we have $\left\langle S^2\right\rangle_n=f(X)$.
Thus, by \eqref{eq:logsobolev}, we get
\[\mathrm{Ent}\left(\exp\left(\lambda\left\langle S^2\right\rangle_n\right)\right)\leq\lambda^2\cdot\mathbb{E}\left[\exp\left(\lambda\left\langle S^2\right\rangle_n\right)\cdot\sum_{u\in T_n^*}\left(f(X)-f\left(\overline{X}^u\right)\right)^2\right]\quad\text{for all $\lambda\in\mathbb{R}$.}\]
To complete the proof, let us show that
\begin{equation}\label{eq:logsobolevIn}
\sum_{u\in T_n^*}\left(f(X)-f\left(\overline{X}^u\right)\right)^2\leq32L^2Y_n+32L^4\gamma_n.
\end{equation}
First, fix $u\in T_n^*$, and let us calculate the difference
\begin{eqnarray*}
f(X)-f\left(\overline{X}^u\right)&=&\frac{1}{N_n}\cdot\sum_{\substack{v\in\partial T_n\\v\succeq u}}\left(S_v^2-\left(S_v-\left(X_u-X'_u\right)\right)^2\right)\\
&=&\frac{1}{N_n}\cdot\sum_{\substack{v\in\partial T_n\\v\succeq u}}\left(2\cdot S_v\cdot\left(X_u-X'_u\right)-\left(X_u-X'_u\right)^2\right)\\
&=&\frac{2}{N_n}\cdot\left(\sum_{\substack{v\in\partial T_n\\v\succeq u}}S_v\right)\cdot\left(X_u-X'_u\right)-\frac{N_n(u)}{N_n}\cdot\left(X_u-X'_u\right)^2.
\end{eqnarray*}
Using the inequality $(a+b)^2\leq2\left(a^2+b^2\right)$ for all $a,b\in\mathbb{R}$, we deduce that
\begin{eqnarray*}
\left(f(X)-f\left(\overline{X}^u\right)\right)^2&\leq&2\left(\frac{4}{N_n^2}\cdot\left(\sum_{\substack{v\in\partial T_n\\v\succeq u}}S_v\right)^2\cdot\left(X_u-X'_u\right)^2+\left(\frac{N_n(u)}{N_n}\right)^2\cdot\left(X_u-X'_u\right)^4\right)\\
&\leq&2\left(\frac{4}{N_n^2}\cdot\left(\sum_{\substack{v\in\partial T_n\\v\succeq u}}S_v\right)^2\cdot(2L)^2+\left(\frac{N_n(u)}{N_n}\right)^2\cdot(2L)^4\right).
\end{eqnarray*}
Then, summing over $u$, we obtain \eqref{eq:logsobolevIn}, as desired.
Indeed, on the one hand, we have
\begin{eqnarray*}
\frac{1}{N_n^2}\cdot\sum_{u\in T_n^*}\left(\sum_{\substack{v\in\partial T_n\\v\succeq u}}S_v\right)^2&=&\frac{1}{N_n^2}\cdot\sum_{u\in T_n^*}\sum_{\substack{v,w\in\partial T_n\\v,w\succeq u}}S_vS_w\\
&=&\frac{1}{N_n^2}\cdot\sum_{v,w\in\partial T_n}\sum_{\substack{u\in T_n^*\\u\preceq v,w}}S_vS_w\\
&=&\frac{1}{N_n^2}\cdot\sum_{v,w\in\partial T_n}S_vS_w\cdot|v\wedge w|=Y_n.
\end{eqnarray*}
On the other hand, by \eqref{eq:boundgamman}, we have
\[\sum_{u\in T_n^*}\left(\frac{N_n(u)}{N_n}\right)^2\leq\gamma_n.\]
The proof is now complete.
\end{proof}

\begin{proof}[Proof of Lemma \ref{lem:Jn}]
Consider the function $g:\mathbb{R}^{T_n^*}\rightarrow\mathbb{R}$ given by
\[g(x)=\frac{1}{N_n^2}\cdot\sum_{v,w\in\partial T_n}\left(\sum_{\substack{u\in T_n^*\\u\preceq v}}x_u\right)\left(\sum_{\substack{u\in T_n^*\\u\preceq w}}x_u\right)\cdot|v\wedge w|\quad\text{for all $x=(x_u)_{u\in T_n^*}\in\mathbb{R}^{T_n^*}$.}\]
By definition, we have $Y_n=g(X)$.
Thus, by \eqref{eq:logsobolev}, we get
\[\mathrm{Ent}\left(e^{\lambda Y_n}\right)\leq\lambda^2\cdot\mathbb{E}\left[e^{\lambda Y_n}\cdot\sum_{u\in T_n^*}\left(g(X)-g\left(\overline{X}^u\right)\right)^2\right]\quad\text{for all $\lambda\in\mathbb{R}$.}\]
To complete the proof, let us show that
\begin{equation}\label{eq:logsobolevJn}
\sum_{u\in T_n^*}\left(g(X)-g\left(\overline{X}^u\right)\right)^2\leq32L^2Z_n+32L^4\gamma_n.
\end{equation}
First, fix $u\in T_n^*$, and let us calculate, for all $v,w\in\partial T_n$:
\begin{eqnarray*}
\lefteqn{S_vS_w-\left(S_v-\left(X_u-X'_u\right)\cdot\mathbf{1}(v\succeq u)\right)\left(S_w-\left(X_u-X_u'\right)\cdot\mathbf{1}(w\succeq u)\right)}\\
&=&S_v\cdot\left(X_u-X'_u\right)\cdot\mathbf{1}(w\succeq u)+S_w\cdot\left(X_u-X'_u\right)\cdot\mathbf{1}(v\succeq u)-\left(X_u-X'_u\right)^2\cdot\mathbf{1}(v,w\succeq u).
\end{eqnarray*}
It follows that
\begin{eqnarray*}
g(X)-g\left(\overline{X}^u\right)&=&\frac{1}{N_n^2}\cdot\left(\sum_{\substack{v,w\in\partial T_n\\w\succeq u}}S_v\cdot|v\wedge w|\right)\cdot\left(X_u-X'_u\right)\\
&&+\frac{1}{N_n^2}\cdot\left(\sum_{\substack{v,w\in\partial T_n\\v\succeq u}}S_w\cdot|v\wedge w|\right)\cdot\left(X_u-X'_u\right)\\
&&-\frac{1}{N_n^2}\cdot\left(\sum_{\substack{v,w\in\partial T_n\\v,w\succeq u}}|v\wedge w|\right)\cdot\left(X_u-X'_u\right)^2\\
&=&\frac{2}{N_n^2}\cdot\left(\sum_{\substack{v,w\in\partial T_n\\w\succeq u}}S_v\cdot|v\wedge w|\right)\cdot\left(X_u-X'_u\right)\\
&&-\frac{1}{N_n^2}\left(\sum_{\substack{v,w\in\partial T_n\\v,w\succeq u}}|v\wedge w|\right)\cdot\left(X_u-X'_u\right)^2.
\end{eqnarray*}
Using the inequality $(a+b)^2\leq2a^2+2b^2$ for all $a,b\in\mathbb{R}$, we deduce that
\begin{eqnarray*}
\left(g(X)-g\left(\overline{X}^u\right)\right)^2&\leq&\frac{8}{N_n^4}\cdot\left(\sum_{\substack{v,w\in\partial T_n\\w\succeq u}}S_v\cdot|v\wedge w|\right)^2\cdot\left(X_u-X'_u\right)^2\\
&&+\frac{2}{N_n^4}\cdot\left(\sum_{\substack{v,w\in\partial T_n\\v,w\succeq u}}|v\wedge w|\right)^2\cdot\left(X_u-X'_u\right)^4\\
&\leq&\frac{8}{N_n^4}\cdot\left(\sum_{\substack{v,w\in\partial T_n\\w\succeq u}}S_v\cdot|v\wedge w|\right)^2\cdot(2L)^2\\
&&+\frac{2}{N_n^4}\cdot\left(\sum_{\substack{v,w\in\partial T_n\\v,w\succeq u}}|v\wedge w|\right)^2\cdot(2L)^4.
\end{eqnarray*}
Then, summing over $u$, we obtain \eqref{eq:logsobolevJn}, as desired.
Indeed, on the one hand, we have
\begin{eqnarray*}
\lefteqn{\frac{1}{N_n^4}\cdot\sum_{u\in T_n^*}\left(\sum_{\substack{v,w\in\partial T_n\\w\succeq u}}S_v\cdot|v\wedge w|\right)^2}\\
&=&\frac{1}{N_n^4}\cdot\sum_{u\in T_n^*}\sum_{\substack{v_1,v_2,w_1,w_2\in\partial T_n\\w_1,w_2\succeq u}}S_{v_1}S_{v_2}\cdot|v_1\wedge w_1|\cdot|v_2\wedge w_2|\\
&=&\frac{1}{N_n^4}\cdot\sum_{v_1,v_2,w_1,w_2\in\partial T_n}\sum_{\substack{u\in T_n^*\\u\preceq w_1,w_2}}S_{v_1}S_{v_2}\cdot|v_1\wedge w_1|\cdot|v_2\wedge w_2|\\
&=&\frac{1}{N_n^4}\cdot\sum_{v_1,v_2,w_1,w_2\in\partial T_n}S_{v_1}S_{v_2}\cdot|v_1\wedge w_1|\cdot|v_2\wedge w_2|\cdot|w_1\wedge w_2|\\
&=&\frac{1}{N_n^2}\cdot\sum_{v_1,v_2\in\partial T_n}S_{v_1}S_{v_2}\cdot\left(\frac{1}{N_n^2}\cdot\sum_{w_1,w_2\in\partial T_n}|v_1\wedge w_1|\cdot|v_2\wedge w_2|\cdot|w_1\wedge w_2|\right)\\
&=&\frac{1}{N_n^2}\cdot\sum_{v_1,v_2\in\partial T_n}S_{v_1}S_{v_2}\cdot\zeta_n(v_1,v_2)=Z_n.
\end{eqnarray*}
On the other hand, we have
\begin{eqnarray*}
\lefteqn{\frac{1}{N_n^4}\cdot\sum_{u\in T_n^*}\left(\sum_{\substack{v,w\in\partial T_n\\v,w\succeq u}}|v\wedge w|\right)^2}\\
&=&\frac{1}{N_n^4}\cdot\sum_{u\in T_n^*}\sum_{\substack{v_1,v_2,w_1,w_2\in\partial T_n\\v_1,v_2,w_1,w_2\succeq u}}|v_1\wedge w_1|\cdot|v_2\wedge w_2|\\
&=&\frac{1}{N_n^4}\cdot\sum_{v_1,v_2,w_1,w_2\in\partial T_n}\sum_{\substack{u\in T_n^*\\u\preceq v_1,v_2,w_1,w_2}}|v_1\wedge w_1|\cdot|v_2\wedge w_2|\\
&=&\frac{1}{N_n^4}\cdot\sum_{v_1,v_2,w_1,w_2\in\partial T_n}|v_1\wedge w_1|\cdot|v_2\wedge w_2|\cdot\#\{u\in T_n^*:u\preceq v_1,v_2,w_1,w_2\}\\
&\leq&\frac{1}{N_n^4}\cdot\sum_{v_1,v_2,w_1,w_2\in\partial T_n}|v_1\wedge w_1|\cdot|v_2\wedge w_2|\cdot|w_1\wedge w_2|\\
&=&\frac{1}{N_n^2}\cdot\sum_{v_1,v_2\in\partial T_n}\left(\frac{1}{N_n^2}\cdot\sum_{w_1,w_2\in\partial T_n}|v_1\wedge w_1|\cdot|v_2\wedge w_2|\cdot|w_1\wedge w_2|\right)\\
&=&\frac{1}{N_n^2}\cdot\sum_{v_1,v_2\in\partial T_n}\zeta_n(v_1,v_2)\leq\gamma_n.
\end{eqnarray*}
The last inequality relies on the upper bound
\begin{equation}\label{eq:thetan}
\zeta_n(v_1,v_2)\leq\gamma_n\quad\text{for all $v_1,v_2\in\partial T_n$,}
\end{equation}
which follows from the generalised Hölder inequality:
\begin{eqnarray*}
\lefteqn{\zeta_n(v_1,v_2)^3}\\
&\leq&\left(\frac{1}{N_n^2}\cdot\sum_{w_1,w_2\in\partial T_n}|v_1\wedge w_1|^3\right)\cdot\left(\frac{1}{N_n^2}\cdot\sum_{w_1,w_2\in\partial T_n}|v_2\wedge w_2|^3\right)\cdot\left(\frac{1}{N_n^2}\cdot\sum_{w_1,w_2\in\partial T_n}|w_1\wedge w_2|^3\right)\\
&=&\left(\frac{1}{N_n}\cdot\sum_{w_1\in\partial T_n}|v_1\wedge w_1|^3\right)\cdot\left(\frac{1}{N_n}\cdot\sum_{w_2\in\partial T_n}|v_2\wedge w_2|^3\right)\cdot\left(\frac{1}{N_n^2}\cdot\sum_{w_1,w_2\in\partial T_n}|w_1\wedge w_2|^3\right),
\end{eqnarray*}
together with the same calculation as in \eqref{eq:calculation}, which shows that for each $v\in\partial T_n$, we have
\[\frac{1}{N_n}\cdot\sum_{w\in\partial T_n}|v\wedge w|^3\leq\gamma_n.\]
Plugging the last inequality into the previous generalised Hölder inequality yields ${\zeta_n(v_1,v_2)^3\leq\gamma_n^3}$, hence \eqref{eq:thetan}.
The proof of the lemma is now complete.
\end{proof}

\begin{proof}[Proof of Lemma \ref{lem:Jn'}]
Recall that
\[Z_n=\frac{1}{N_n^4}\cdot\sum_{u\in T_n^*}\left(\sum_{\substack{v,w\in\partial T_n\\w\succeq u}}S_v\cdot|v\wedge w|\right)^2\geq0.\]
In particular, by the monotone convergence theorem, we have
\[\mathbb{E}\left[e^{\lambda Z_n}\right]=\sum_{k\geq0}\frac{\lambda^k\cdot\mathbb{E}\left[Z_n^k\right]}{k!}\quad\text{for all $\lambda\in\mathbb{R}_+$.}\]
Now, for each $k\in\mathbb{N}$, by \eqref{eq:defIn}, we have
\[\mathbb{E}\left[Z_n^k\right]=\frac{1}{N_n^{2k}}\cdot\sum_{v_1,\ldots,v_{2k}\in\partial T_n}\mathbb{E}\left[S_{v_1}\ldots S_{v_{2k}}\right]\cdot\zeta_n(v_1,v_2)\ldots\zeta_n(v_{2k-1},v_{2k}).\]
At this point, notice that $\mathbb{E}\left[S_{v_1}\ldots S_{v_{2k}}\right]\geq0$ for all $v_1,\ldots,v_{2k}\in\partial T_n$.
Indeed, we can write
\begin{eqnarray*}
\mathbb{E}\left[S_{v_1}\ldots S_{v_{2k}}\right]&=&\sum_{\substack{u_1,\ldots,u_{2k}\in T_n^*\\u_1\preceq v_1,\ldots,u_{2k}\preceq v_{2k}}}\mathbb{E}\left[X_{u_1}\ldots X_{u_{2k}}\right]\notag\\
&=&\sum_{\substack{u_1,\ldots,u_{2k}\in T_n^*\\u_1\preceq v_1,\ldots,u_{2k}\preceq v_{2k}}}\mathbb{E}\left[\prod_{u\in T_n^*}X_u^{\#\{i\in\llbracket1,2k\rrbracket:u_i=u\}}\right]\notag\\
&=&\sum_{\substack{u_1,\ldots,u_{2k}\in T_n^*\\u_1\preceq v_1,\ldots,u_{2k}\preceq v_{2k}}}\prod_{u\in T_n^*}\mathbb{E}\left[X_u^{\#\{i\in\llbracket1,2k\rrbracket:u_i=u\}}\right].
\end{eqnarray*}
Now, in this sum, the only non-zero terms are those for which $\#\{i\in\llbracket1,2k\rrbracket:u_i=u\}$ is even for all $u\in T_n^*$, which are non-negative.
Indeed, by the symmetry of $\nu$, we have
\[\prod_{u\in T_n^*}\mathbb{E}\left[X_u^{\#\{i\in\llbracket1,2k\rrbracket:u_i=u\}}\right]=0\]
as soon as there exists $u\in T_n^*$ for which $\#\{i\in\llbracket1,2k\rrbracket:u_i=u\}$ is odd.
This shows that $\mathbb{E}\left[S_{v_1}\ldots S_{v_{2k}}\right]\geq0$, which allows us to bound $\mathbb{E}\left[Z_n^k\right]$ as follows, using \eqref{eq:thetan}:
\[\mathbb{E}\left[Z_n^k\right]\leq\frac{1}{N_n^{2k}}\cdot\sum_{v_1,\ldots,v_{2k}\in\partial T_n}\mathbb{E}\left[S_{v_1}\ldots S_{v_{2k}}\right]\cdot\gamma_n^k=\gamma_n^k\cdot\mathbb{E}\left[\langle S\rangle_n^{2k}\right].\]
Thus, we obtain
\[\mathbb{E}\left[e^{\lambda Z_n}\right]\leq\sum_{k\geq0}\frac{\lambda^k\cdot\gamma_n^k\cdot\mathbb{E}\left[\langle S\rangle_n^{2k}\right]}{k!}=\mathbb{E}\left[\exp\left(\lambda\cdot\gamma_n\cdot\langle S\rangle_n^2\right)\right]\quad\text{for all $\lambda\in\mathbb{R}_+$.}\]
The desired result then follows from Proposition \ref{prop:concentration1}.
\end{proof}

\begin{rem}
Notice that we have used at several places the assumption that $\nu$ has compact support (in the proof of Lemma \ref{lem:In} and Lemma \ref{lem:Jn}) and is symmetric (in the proof of Lemma \ref{lem:Jn'}).
Although we do not claim that these assumptions are minimal for the result of Theorem \ref{thm:concentration2} to hold, our proof does not seem to work without them.
\end{rem}

\section{Integer-valued random height functions}\label{sec:ivrhf}

In this section, we consider our $\nu$-integer-valued random height functions on $\mathcal{T}_n^{p_n}$, and prove Theorem \ref{thm:main}.
We are back to assuming that $\nu$ is a discrete probability measure on $\mathbb{Z}$, which is symmetric and has finite support.

\begin{rem}
We are not claiming that these are the minimal assumptions under which the result of Theorem \ref{thm:main} holds.
However, we will use them at several places in the proof: in particular, they will be hidden in our use of Theorem \ref{thm:concentration2} which, as mentioned above, we do not know how to prove without these assumptions.
\end{rem}

Let $\sigma^2=\sum_{x\in\mathbb{Z}}x^2\cdot\nu(\{x\})$ be the variance of $\nu$, and set ${L=\inf\{\ell\in\mathbb{N}:\nu(\llbracket-\ell,\ell\rrbracket)=1\}}$.
We assume that $\nu(\{0\})<1$, so that $\sigma^2,L>0$ (otherwise, the statement of Theorem \ref{thm:main} is trivial, understanding $\mathcal{N}(0,0)$ as the Dirac probability measure $\delta_0$).

\subsection{Preliminaries}\label{subsec:prelivrhf}

We first collect some basic tools to estimate the height at the root of the tree for $\nu$-integer-valued random height functions.
Throughout this subsection, we fix a finite rooted tree $T_n$ with height $n\in\mathbb{N}^*$ and with no internal leaves, and we let ${H_n=(H_n(u))_{u\in T_n}}$ be a $\nu$-integer-valued random height function on $T_n$ with zero boundary conditions on $B_n=\partial T_n$.
To lighten the notation, we set ${N_n:=\#\partial T_n}$.
By definition (see \eqref{eq:ivrhf}), for any Borel function $\varphi:\mathbb{R}\rightarrow\mathbb{R}_+$, we have
\begin{equation}\label{eq:distHno}
\mathbb{E}[\varphi(H_n(\varnothing))]=\frac{\sum_{y\in\mathbb{Z}}\varphi(y)\cdot\pi_n(y)}{\sum_{y\in\mathbb{Z}}\pi_n(y)},
\end{equation}
where
\[\pi_n(y)=\sum_{\substack{h\in\mathbb{H}^0(T_n,\partial T_n)\\h_\varnothing=y}}\prod_{(u,v)\in\vec{E}(T_n)}\nu(\{h_v-h_u\})\quad\text{for all $y\in\mathbb{Z}$.}\]
To get a grasp on these $\pi_n(\cdot)$ terms, we first make the following observation (see Proposition \ref{prop:obs} below).
For each $k\in\mathbb{N}$, let $P_k:\mathbb{Z}\rightarrow[0,1]$ be the probability mass function of $X_1+\ldots+X_k$, where $X_1,\ldots,X_k$ are independent and have distribution $\nu$.
Moreover, let $(X_u)_{u\in T_n^*}$ be a collection of independent random variables with distribution $\nu$, where $T_n^*=T_n\setminus\{\varnothing\}$, and set
\[S_v=\sum_{\substack{u\in T_n^*\\u\preceq v}}X_u\quad\text{for all $v\in T_n$.}\]

\begin{prop}\label{prop:obs}
Let $\kappa_n$ be the smallest integer $k\in\llbracket0,n\rrbracket$ such that every node $u\in T_n$ with $k\leq|u|\leq n-1$ has exactly one child.\footnote{Note that this property holds vacuously for $k=n$.}
For all $y\in\mathbb{Z}$, we have
\[\pi_n(y)=\mathbb{E}\left[\prod_{v\in\partial[T_n]_{\kappa_n}}P_{n-\kappa_n}(y+S_v)\right].\]
\end{prop}
\begin{proof}
Let us denote by $\mathbb{H}_0([T_n]_{\kappa_n})$ the set of functions ${f=(f_u)_{u\in[T_n]_{\kappa_n}}\in\mathbb{Z}^{[T_n]_{\kappa_n}}}$ such that $f_\varnothing=0$.
For all $y\in\mathbb{Z}$, we have
\begin{eqnarray*}
\pi_n(y)&=&\sum_{f\in\mathbb{H}_0([T_n]_{\kappa_n})}\sum_{\substack{h\in\mathbb{H}^0(T_n,\partial T_n)\\\text{$h=y+f$ on $[T_n]_{\kappa_n}$}}}\prod_{(u,v)\in\vec{E}(T_n)}\nu(\{h_v-h_u\}).
\end{eqnarray*}
Now, we claim that for each $f\in\mathbb{H}_0([T_n]_{\kappa_n})$, we have
\begin{equation}\label{eq:cut}
\sum_{\substack{h\in\mathbb{H}^0(T_n,\partial T_n)\\\text{$h=y+f$ on $[T_n]_{\kappa_n}$}}}\prod_{(u,v)\in\vec{E}(T_n)}\nu(\{h_v-h_u\})=\left(\prod_{w\in\partial[T_n]_{\kappa_n}}P_{n-\kappa_n}(y+f_w)\right)\cdot\mathbb{P}(S=f).
\end{equation}
Indeed, for all $h\in\mathbb{H}^0(T_n,\partial T_n)$ such that $h=y+f$ on $[T_n]_{\kappa_n}$, on the one hand, we have
\[\prod_{\substack{(u,v)\in\vec{E}(T_n)\\|v|\leq\kappa_n}}\nu(\{h_v-h_u\})=\prod_{(u,v)\in\vec{E}([T_n]_{\kappa_n})}\nu(\{f_v-f_u\})=\mathbb{P}(S=f).\]
On the other hand, by the definition of $\kappa_n$, the edges ${(u,v)\in\vec{E}(T_n)}$ with $|u|\geq\kappa_n$ form a collection of $N_n$ disjoint paths with length $n-\kappa_n$, which we may index by the nodes $w\in\partial[T_n]_{\kappa_n}$.
It follows that
\begin{eqnarray*}
\lefteqn{\sum_{\substack{h\in\mathbb{H}^0(T_n,\partial T_n)\\\text{$h=y+f$ on $[T_n]_{\kappa_n}$}}}\prod_{\substack{(u,v)\in\vec{E}(T_n)\\|u|\geq\kappa_n}}\nu(\{h_v-h_u\})}\\
&=&\sum_{\substack{\left(x^w_1,\ldots,x^w_{n-\kappa_n}\right)_{w\in\partial[T_n]_{\kappa_n}}\in(\mathbb{Z}^{n-\kappa_n})^{\partial[T_n]_{\kappa_n}}\\\text{$y+f_w+x^w_1+\ldots+x^w_{n-\kappa_n}=0$ for all $w\in\partial[T_n]_{\kappa_n}$}}}\prod_{w\in\partial[T_n]_{\kappa_n}}\left(\nu(\{x^w_1\})\cdot\ldots\cdot\nu(\{x^w_{n-\kappa_n}\})\right)\\
&=&\prod_{w\in\partial[T_n]_{\kappa_n}}\left(\sum_{\substack{x_1+\ldots+x_{n-\kappa_n}\in\mathbb{Z}\\x_1+\ldots+x_{n-\kappa_n}=-(y+f_w)}}\nu(\{x_1\})\cdot\ldots\cdot\nu(\{x_{n-\kappa_n}\})\right)\\
&=&\prod_{w\in\partial[T_n]_{\kappa_n}}P_{n-\kappa_n}(-(y+f_w))=\prod_{w\in\partial[T_n]_{\kappa_n}}P_{n-\kappa_n}(y+f_w),
\end{eqnarray*}
where the last equality holds by the symmetry of $\nu$.
Altogether, we obtain \eqref{eq:cut}, and it follows that
\[\pi_n(y)=\sum_{f\in\mathbb{H}_0([T_n]_{\kappa_n})}\left(\prod_{w\in\partial[T_n]_{\kappa_n}}P_{n-\kappa_n}(y+f_w)\right)\cdot\mathbb{P}(S=f)=\mathbb{E}\left[\prod_{w\in\partial[T_n]_{\kappa_n}}P_{n-\kappa_n}(y+S_w)\right],\]
as desired.
\end{proof}

Next, we state some local central limit estimates for the random walk with step distribution $\nu$.
Let $d\in\mathbb{N}^*$ be the step of the lattice distribution $\nu$, and fix $r\in\mathbb{Z}$ such that $\nu(r+d\mathbb{Z})=1$.
For each $k\in\mathbb{N}^*$, we set
\[q_k(z)=\frac{d}{\sqrt{2\pi\sigma^2k}}\cdot\exp\left(-\frac{z^2}{2\sigma^2k}\right)\quad\text{for all $z\in rk+d\mathbb{Z}$.}\]
The estimates below, formulated in a ``multiplicative'' way for future use, quantify the approximation
\[P_k(z)\approx q_k(z)=q_k(0)\cdot\exp\left(-\frac{z^2}{2\sigma^2k}\right)\]
in the range $|z|\lesssim k$, and provide some complementary control on the ratio $P_k(z)/q_k(0)$.

\begin{lem}\label{lem:lclt}
There exists $\eta\in\mathbb{R}_+^*$, and constants $C\in\mathbb{R}_+^*$ and $c\in\left]0,\left.1\middle/\left(2\sigma^2\right)\right.\right[$ such that the following holds.
\begin{enumerate}[label=(\roman*)]
\item For all $k\in\mathbb{N}^*$, we have
\[\exp\left(-\frac{z^2}{2\sigma^2k}-C\left(\frac{1}{k}+\frac{z^4}{k^3}\right)\right)\leq\frac{P_k(z)}{q_k(0)}\leq\exp\left(-\frac{z^2}{2\sigma^2k}+C\left(\frac{1}{k}+\frac{z^4}{k^3}\right)\right)\]
for all $z\in rk+d\mathbb{Z}$ such that $|z|\leq\eta k$.

\item For all $k\in\mathbb{N}^*$, we have
\[\frac{P_k(z)}{q_k(0)}\leq\exp\left(-\frac{cz^2}{k}+\frac{C}{k}\right)\quad\text{for all $z\in rk+d\mathbb{Z}$.}\]
\end{enumerate}
\end{lem}

We give a proof of these standard estimates in Appendix \ref{sec:lclt}.
Now, as in Proposition \ref{prop:obs}, let $\kappa_n$ be the smallest integer $k\in\llbracket0,n\rrbracket$ such that every node $u\in T_n$ with $k\leq|u|\leq n-1$ has exactly one child.
The local central limit estimates of Lemma \ref{lem:lclt} will enable us to quantify the approximation
\begin{eqnarray*}
\pi_n(y)&=&q_{n-\kappa_n}(0)^{N_n}\cdot\mathbb{E}\left[\prod_{v\in\partial[T_n]_{\kappa_n}}\frac{P_{n-\kappa_n}(y+S_v)}{q_{n-\kappa_n}(0)}\right]\\
&\approx&q_{n-\kappa_n}(0)^{N_n}\cdot\mathbb{E}\left[\prod_{v\in\partial[T_n]_{\kappa_n}}\exp\left(-\frac{(y+S_v)^2}{2\sigma^2(n-\kappa_n)}\right)\right].
\end{eqnarray*}
Note that by the definition of $\kappa_n$, we have $\#\partial[T_n]_{\kappa_n}=N_n$.
For the moment, let us calculate the last expectation, in terms of the random variables
\[\langle S\rangle_{\kappa_n}:=\frac{1}{N_n}\cdot\sum_{v\in\partial[T_n]_{\kappa_n}}S_v\quad\text{and}\quad\left\langle S^2\right\rangle_{\kappa_n}:=\frac{1}{N_n}\cdot\sum_{v\in\partial[T_n]_{\kappa_n}}S_v^2.\]
Note that $\mathbb{E}[\langle S\rangle_n]=0$ and $\mathbb{E}\left[\left\langle S^2\right\rangle_n\right]=\sigma^2\cdot\kappa_n$.

\begin{lem}\label{lem:guide}
Fix $\rho\in\mathbb{R}_+^*$.
For all $y\in\mathbb{Z}$, we have
\begin{eqnarray*}
\lefteqn{\mathbb{E}\left[\prod_{v\in\partial[T_n]_{\kappa_n}}\exp\left(-\rho(y+S_v)^2\right)\right]}\\
&=&\exp\left(-\rho N_ny^2\right)\cdot\mathbb{E}\left[\exp\left(-\rho N_n\left(2y\cdot\langle S\rangle_{\kappa_n}+\left(\left\langle S^2\right\rangle_{\kappa_n}-\sigma^2\kappa_n\right)\right)\right)\right]\cdot\exp\left(-\rho N_n\sigma^2\kappa_n\right).
\end{eqnarray*}
\end{lem}
\begin{proof}
We have
\[\sum_{v\in\partial[T_n]_{\kappa_n}}(y+S_v)^2=y^2\cdot N_n+2y\cdot\sum_{v\in\partial[T_n]_{\kappa_n}}S_v+\sum_{v\in\partial[T_n]_{\kappa_n}}S_v^2=N_n\left(y^2+2y\cdot\langle S\rangle_{\kappa_n}+\left\langle S^2\right\rangle_{\kappa_n}\right),\]
hence
\begin{eqnarray*}
\lefteqn{\mathbb{E}\left[\prod_{v\in\partial[T_n]_{\kappa_n}}\exp\left(-\rho(y+S_v)^2\right)\right]}\\
&=&\mathbb{E}\left[\exp\left(-\rho N_n\left(y^2+2y\cdot\langle S\rangle_{\kappa_n}+\left\langle S^2\right\rangle_{\kappa_n}\right)\right)\right]\\
&=&\exp\left(-\rho N_ny^2\right)\cdot\mathbb{E}\left[\exp\left(-\rho N_n\left(2y\cdot\langle S\rangle_{\kappa_n}+\left(\left\langle S^2\right\rangle_{\kappa_n}-\sigma^2\kappa_n\right)\right)\right)\right]\cdot\exp\left(-\rho N_n\sigma^2\kappa_n\right),
\end{eqnarray*}
as desired.
\end{proof}

\subsection{Key technical lemmas}

Finally, we have all the ingredients at hand to prove Theorem \ref{thm:main}.
For the convenience of future use, we assemble these ingredients into the technical statements of Lemma \ref{lem:techbulk} and Lemma \ref{lem:techtails} below, which we work our way towards stating now.
In order to establish the Gaussian scaling limit for the height at the root of the tree in Theorem \ref{thm:main}, we want to quantify the approximation
\begin{equation}\label{eq:gaussianapprox}
\mathbb{P}\left(\frac{\sqrt{N_n}}{\sqrt{n}}\cdot H_n(\varnothing)\in{]a,b]}\right)\approx\frac{1}{\sqrt{2\pi\sigma^2}}\cdot\int_a^b\exp\left(-\frac{x^2}{2\sigma^2}\right)\mathrm{d}x.
\end{equation}
First, let us rewrite \eqref{eq:distHno} as follows: for any Borel function $\varphi:\mathbb{R}\rightarrow\mathbb{R}_+$, we have\footnote{Indeed, since $\nu(r+d\mathbb{Z})=1$ and since $\nu$ is symmetric, we have $\pi_n(y)>0$ only if $y\in rn+d\mathbb{Z}$, hence
\begin{eqnarray*}
\mathbb{E}\left[\varphi\left(\frac{\sqrt{N_n}}{\sqrt{n}}\cdot H_n(\varnothing)\right)\right]&=&\frac{\sum_{y\in\mathbb{Z}}\varphi\left(\frac{\sqrt{N_n}}{\sqrt{n}}\cdot y\right)\cdot\pi_n(y)}{\sum_{y\in\mathbb{Z}}\pi_n(y)}\\
&=&\frac{\sum_{z\in\mathbb{Z}}\varphi\left(\frac{\sqrt{N_n}}{\sqrt{n}}\cdot(rn+dz)\right)\cdot\pi_n(rn+dz)}{\sum_{z\in\mathbb{Z}}\pi_n(rn+dz)}\\
&=&\frac{\int_{-\infty}^\infty\varphi\left(\frac{\sqrt{N_n}}{\sqrt{n}}\cdot(rn+d\lfloor w\rfloor)\right)\cdot\pi_n(rn+d\lfloor w\rfloor)\mathrm{d}w}{\int_{-\infty}^\infty\pi_n(rn+d\lfloor w\rfloor)\mathrm{d}w}.
\end{eqnarray*}
Then, perform the change-of-variables
\[w=-\frac{rn}{d}+\frac{x\sqrt{n}}{d\sqrt{N_n}}.\]}
\begin{equation}\label{eq:igor}
\mathbb{E}\left[\varphi\left(\frac{\sqrt{N_n}}{\sqrt{n}}\cdot H_n(\varnothing)\right)\right]=\frac{\int_{-\infty}^\infty\varphi\left(\frac{\sqrt{N_n}}{\sqrt{n}}\cdot y_n(x)\right)\cdot \pi_n(y_n(x))\mathrm{d}x}{\int_{-\infty}^\infty \pi_n(y_n(x))\mathrm{d}x}.
\end{equation}
where
\[y_n(x)=rn+d\left\lfloor-\frac{rn}{d}+\frac{x\sqrt{n}}{d\sqrt{N_n}}\right\rfloor\quad\text{for all $x\in\mathbb{R}$.}\]
Note that $y_n(x)\in rn+d\mathbb{Z}$, and
\[\frac{\sqrt{N_n}}{\sqrt{n}}\cdot y_n(x)\approx x.\]
In view of \eqref{eq:gaussianapprox}, we want to identify a normalising factor $\lambda_n\in\mathbb{R}_+^*$ such that
\[\frac{1}{\lambda_n}\cdot\pi_n(y_n(x))\approx\exp\left(-\frac{x^2}{2\sigma^2}\right).\]
In hindsight, let us introduce the following appropriate renormalisation of the $\pi_n(y_n(\cdot))$ terms:
Assume that $\kappa_n\leq n-1$.
For all $x\in\mathbb{R}$, we set
\[\phi_n(x)=\exp\left(\frac{N_n\kappa_n}{2(n-\kappa_n)}\right)\cdot\frac{\pi_n(y_n(x))}{q_{n-\kappa_n}(0)^{N_n}}.\]
By \eqref{eq:igor}, for any Borel function $\varphi:\mathbb{R}\rightarrow\mathbb{R}_+$, we have
\begin{equation}\label{eq:coupdegras}
\mathbb{E}\left[\varphi\left(\frac{\sqrt{N_n}}{\sqrt{n}}\cdot H_n(\varnothing)\right)\right]=\frac{\int_{-\infty}^\infty\varphi\left(\frac{\sqrt{N_n}}{\sqrt{n}}\cdot y_n(x)\right)\cdot\phi_n(x)\mathrm{d}x}{\int_{-\infty}^\infty\phi_n(x)\mathrm{d}x}.
\end{equation}
Now, by Proposition \ref{prop:obs}, we have
\[\phi_n(x)=\exp\left(\frac{N_n\kappa_n}{2(n-\kappa_n)}\right)\cdot\mathbb{E}\left[\prod_{v\in\partial[T_n]_{\kappa_n}}\frac{P_{n-\kappa_n}(y_n(x)+S_v)}{q_{n-\kappa_n}(0)}\right].\]
Thus, in view of the local central limit estimates of Lemma \ref{lem:lclt}, we should have
\[\phi_n(x)\approx\exp\left(\frac{N_n\kappa_n}{2(n-\kappa_n)}\right)\cdot\mathbb{E}\left[\prod_{v\in\partial[T_n]_{\kappa_n}}\exp\left(-\frac{(y_n(x)+S_v)^2}{2\sigma^2(n-\kappa_n)}\right)\right].\]
Then, by the calculation of Lemma \ref{lem:guide} with $\rho=\left.1\middle/\left(2\sigma^2(n-\kappa_n)\right)\right.$, we should obtain
\[\phi_n(x)\approx\exp\left(-\frac{N_ny_n(x)^2}{2\sigma^2(n-\kappa_n)}\right)\cdot\mathbb{E}\left[\exp\left(-\frac{N_n\left(2y_n(x)\cdot\langle S\rangle_{\kappa_n}+\left(\left\langle S^2\right\rangle_{\kappa_n}-\sigma^2\kappa_n\right)\right)}{2\sigma^2(n-\kappa_n)}\right)\right].\]
Finally, counting on the concentration of $\langle S\rangle_{\kappa_n}$ and $\left(\left\langle S^2\right\rangle_{\kappa_n}-\sigma^2\kappa_n\right)$ to handle the last expectation, we should obtain
\[\phi_n(x)\approx\exp\left(-\frac{N_ny_n(x)^2}{2\sigma^2(n-\kappa_n)}\right)\approx\exp\left(-\frac{x^2}{2\sigma^2}\right),\]
which should yield the desired Gaussian approximation \eqref{eq:gaussianapprox}, by \eqref{eq:coupdegras}.

The purpose of Lemma \ref{lem:techbulk} and Lemma \ref{lem:techtails} below is to make precise the arguments that we have just unrolled.
Despite the unappealingness of their technical statements, the proof of these lemmas will be straightforward.
In view of Theorem \ref{thm:concentration2}, we set
\[\gamma_n=1\vee\left(\frac{1}{N_n}\cdot\sum_{k=1}^{\kappa_n}k^3\cdot\max_{\substack{u\in T_n\\|u|=k}}\#\{v\in\partial T_n:v\succeq u\}\right),\]
and we let
\[\delta_n=\frac{1}{32L^2\sqrt{2\gamma_n}}.\]
With the following lemma, we will control the $\phi_n(\cdot)$ terms ``in the bulk'':

\begin{lem}[Bulk estimates]\label{lem:techbulk}
Let $\eta\in\mathbb{R}_+^*$, $C\in\mathbb{R}_+^*$ and $c\in\left]0,\left.1\middle/\left(2\sigma^2\right)\right.\right[$ be as in the statement of Lemma \ref{lem:lclt}.
For all $n\in\mathbb{N}^*$ and $\Lambda\in\mathbb{R}_+$ such that the following technical condition is satisfied:
\begin{equation}\label{eq:nandA1}
\kappa_n\leq n-1,\quad\frac{\Lambda\sqrt{n}}{\sqrt{N_n}}+\kappa_n+d\leq\eta(n-\kappa_n)\quad\text{and}\quad\frac{N_n}{2\sigma^2(n-\kappa_n)}\leq\frac{\delta_n}{2},
\end{equation}
the following holds.
For all $x\in[-\Lambda,\Lambda]$, we have
\[\phi_n(x)\geq\exp\left(-\alpha_1(T_n)\cdot x^2-\beta_1(T_n)\cdot|x|-\theta_1(T_n,\Lambda)\right),\]
where $\alpha_1(\cdot)$, $\beta_1(\cdot)$ and $\theta_1(\cdot,\cdot)$ are given by
\begin{equation}\label{eq:alpha1}
\alpha_1(T_n)=\frac{n}{2\sigma^2(n-\kappa_n)},
\end{equation}
\begin{equation}\label{eq:beta1}
\beta_1(T_n)=\frac{dn}{\sigma^2(n-\kappa_n)}\cdot\sqrt{\frac{N_n}{n}},
\end{equation}
and
\begin{equation}\label{eq:theta1}
\theta_1(T_n,\Lambda)=C\left(\frac{N_n}{n-\kappa_n}+\frac{N_n}{(n-\kappa_n)^3}\cdot\left(\frac{\Lambda\sqrt{n}}{\sqrt{N_n}}+\kappa_n+d\right)^4\right)+\frac{d^2N_n}{2\sigma^2(n-\kappa_n)}.
\end{equation}
Moreover, for all $x\in[-\Lambda,\Lambda]$, we have
\[\phi_n(x)\leq\exp\left(-\alpha_2(T_n)\cdot x^2+\beta_2(T_n)\cdot|x|+\theta_2(T_n,\Lambda)\right),\]
where $\alpha_2(\cdot)$, $\beta_2(\cdot)$ and $\theta_2(\cdot,\cdot)$ are given by
\begin{equation}\label{eq:alpha2}
\alpha_2(T_n)=\alpha_1(T_n)-L^2\gamma_n\cdot\left(\frac{n}{\sigma^2(n-\kappa_n)}\right)^2\cdot\frac{N_n}{n},
\end{equation}
\begin{equation}\label{eq:beta2}
\beta_2(T_n)=\beta_1(T_n)+2dL^2\gamma_n\cdot\left(\frac{n}{\sigma^2(n-\kappa_n)}\right)^2\cdot\left(\frac{N_n}{n}\right)^{3/2},
\end{equation}
and
\begin{equation}\label{eq:theta2}
\theta_2(T_n,\Lambda)=\theta_1(T_n,\Lambda)+d^2L^2\gamma_n\cdot\left(\frac{n}{\sigma^2(n-\kappa_n)}\right)^2\cdot\left(\frac{N_n}{n}\right)^2+\frac{2}{\delta_n^2}\cdot\left(\frac{N_n}{2\sigma^2(n-\kappa_n)}\right)^2.
\end{equation}
\end{lem}

Even though the $\alpha_i$, $\beta_i$ and $\theta_i$ terms have rather cumbersome expressions, they all simplify in the limit $n\to\infty$: if $N_n/n\rightarrow0$ and $\kappa_n/n\rightarrow0$, and $\gamma_n\lesssim1$ and $\delta_n\gtrsim1$, then $\alpha_i\rightarrow\left.1/\left(2\sigma^2\right)\right.$ and $\beta_i,\theta_i\rightarrow0$.

\begin{proof}[Proof of Lemma \ref{lem:techbulk}]
Fix $n\in\mathbb{N}^*$ and $\Lambda\in\mathbb{R}_+$, and assume that the technical condition \eqref{eq:nandA1} is satisfied.
First, by Proposition \ref{prop:obs}, we have
\[\phi_n(x)=\exp\left(\frac{N_n\kappa_n}{2(n-\kappa_n)}\right)\cdot\mathbb{E}\left[\prod_{v\in\partial[T_n]_{\kappa_n}}\frac{P_{n-\kappa_n}(y_n(x)+S_v)}{q_{n-\kappa_n}(0)}\right]\quad\text{for all $x\in\mathbb{R}$.}\]
For future reference, note that since $|\lfloor w\rfloor-w|\leq1$ for all $w\in\mathbb{R}$, we have
\begin{equation}\label{eq:ynx}
\left|y_n(x)-\frac{x\sqrt{n}}{\sqrt{N_n}}\right|\leq d\quad\text{for all $x\in\mathbb{R}$.}
\end{equation}
Thus, by \eqref{eq:nandA1}, the following holds: for all $x\in[-\Lambda,\Lambda]$, and for all $v\in\partial[T_n]_{\kappa_n}$, we have
\[|y_n(x)+S_v|\leq\left(\left|\frac{x\sqrt{n}}{\sqrt{N_n}}\right|+d\right)+\kappa_n\leq\frac{\Lambda\sqrt{n}}{\sqrt{N_n}}+\kappa_n+d\leq\eta(n-\kappa_n).\]
Therefore, by Lemma \ref{lem:lclt}, we get that
\begin{eqnarray*}
\lefteqn{\frac{P_{n-\kappa_n}(y_n(x)+S_v)}{q_{n-\kappa_n}(0)}}\\
&\leq&\exp\left(-\frac{(y_n(x)+S_v)^2}{2\sigma^2(n-\kappa_n)}+C\left(\frac{1}{n-\kappa_n}+\frac{(y_n(x)+S_v)^4}{(n-\kappa_n)^3}\right)\right)\\
&\leq&\exp\left(-\frac{(y_n(x)+S_v)^2}{2\sigma^2(n-\kappa_n)}+C\left(\frac{1}{n-\kappa_n}+\frac{1}{(n-\kappa_n)^3}\cdot\left(\frac{\Lambda\sqrt{n}}{\sqrt{N_n}}+\kappa_n+d\right)^4\right)\right).
\end{eqnarray*}
Similarly, we have
\begin{eqnarray*}
\lefteqn{\frac{P_{n-\kappa_n}(y_n(x)+S_v)}{q_{n-\kappa_n}(0)}}\\
&\geq&\exp\left(-\frac{(y_n(x)+S_v)^2}{2\sigma^2(n-\kappa_n)}-C\left(\frac{1}{n-\kappa_n}+\frac{1}{(n-\kappa_n)^3}\cdot\left(\frac{\Lambda\sqrt{n}}{\sqrt{N_n}}+\kappa_n+d\right)^4\right)\right).
\end{eqnarray*}
Taking the product over $v\in\partial[T_n]_{\kappa_n}$, taking expectations and multiplying through by the ad-hoc exponential, we obtain that
\[\phi_n(x)\leq\psi_n(x)\cdot\exp\left(C\left(\frac{N_n}{n-\kappa_n}+\frac{N_n}{(n-\kappa_n)^3}\cdot\left(\frac{\Lambda\sqrt{n}}{\sqrt{N_n}}+\kappa_n+d\right)^4\right)\right),\]
and
\[\phi_n(x)\geq\psi_n(x)\cdot\exp\left(-C\left(\frac{N_n}{n-\kappa_n}+\frac{N_n}{(n-\kappa_n)^3}\cdot\left(\frac{\Lambda\sqrt{n}}{\sqrt{N_n}}+\kappa_n+d\right)^4\right)\right),\]
where
\[\psi_n(x)=\exp\left(\frac{N_n\kappa_n}{2(n-\kappa_n)}\right)\cdot\mathbb{E}\left[\prod_{v\in\partial[T_n]_{\kappa_n}}\exp\left(-\frac{(y_n(x)+S_v)^2}{2\sigma^2(n-\kappa_n)}\right)\right].\]
Now, by the calculation of Lemma \ref{lem:guide} with $\rho=\left.1\middle/\left(2\sigma^2(n-\kappa_n)\right)\right.$, we get that
\[\psi_n(x)=\exp\left(-\frac{N_ny_n(x)^2}{2\sigma^2(n-\kappa_n)}\right)\cdot\mathbb{E}\left[\exp\left(-\frac{N_n\left(2y_n(x)\cdot\langle S\rangle_{\kappa_n}+\left(\left\langle S^2\right\rangle_{\kappa_n}-\sigma^2\kappa_n\right)\right)}{2\sigma^2(n-\kappa_n)}\right)\right].\]
Next, by \eqref{eq:nandA1} and by Corollary \ref{cor:zeta}, we have
\begin{eqnarray*}
\mathbb{E}\left[\exp\left(-\frac{N_n\left(2y_n(x)\cdot\langle S\rangle_{\kappa_n}+\left(\left\langle S^2\right\rangle_{\kappa_n}-\sigma^2\kappa_n\right)\right)}{2\sigma^2(n-\kappa_n)}\right)\right]\geq1,
\end{eqnarray*}
and
\begin{eqnarray*}
\lefteqn{\mathbb{E}\left[\exp\left(-\frac{N_n\left(2y_n(x)\cdot\langle S\rangle_{\kappa_n}+\left(\left\langle S^2\right\rangle_{\kappa_n}-\sigma^2\kappa_n\right)\right)}{2\sigma^2(n-\kappa_n)}\right)\right]}\\
&\leq&\exp\left(L^2\gamma_n\cdot\left(\frac{N_ny_n(x)}{\sigma^2(n-\kappa_n)}\right)^2+\frac{2}{\delta_n^2}\cdot\left(\frac{N_n}{2\sigma^2(n-\kappa_n)}\right)^2\right).
\end{eqnarray*}
Thus, we obtain
\[\phi_n(x)\geq\exp\left(-\frac{N_ny_n(x)^2}{2\sigma^2(n-\kappa_n)}-C\left(\frac{N_n}{n-\kappa_n}+\frac{N_n}{(n-\kappa_n)^3}\cdot\left(\frac{\Lambda\sqrt{n}}{\sqrt{N_n}}+\kappa_n+d\right)^4\right)\right),\]
and
\begin{eqnarray*}
\phi_n(x)&\leq&\exp\left(-\frac{N_ny_n(x)^2}{2\sigma^2(n-\kappa_n)}+C\left(\frac{N_n}{n-\kappa_n}+\frac{N_n}{(n-\kappa_n)^3}\cdot\left(\frac{\Lambda\sqrt{n}}{\sqrt{N_n}}+\kappa_n+d\right)^4\right)\right.\\
&&\left.+L^2\gamma_n\cdot\left(\frac{N_ny_n(x)}{\sigma^2(n-\kappa_n)}\right)^2+\frac{2}{\delta_n^2}\cdot\left(\frac{N_n}{2\sigma^2(n-\kappa_n)}\right)^2\right).
\end{eqnarray*}
Finally, by \eqref{eq:ynx}, we have
\begin{equation}\label{eq:ynx2}
\left|y_n(x)^2-\frac{x^2n}{N_n}\right|\leq d\left(\frac{2|x|\sqrt{n}}{\sqrt{N_n}}+d\right).
\end{equation}
Plugging this into the above bounds, and grouping together the terms in $x^2$, the terms in $|x|$ and the constant terms, we obtain the desired result (let us spare the reader these tedious details).
\end{proof}

\begin{rem}\label{rem:lemtech1}
For future reference, observe that if furthermore the following technical condition is satisfied:
\begin{equation}\label{eq:nandA3}
L^2\gamma_n\cdot\frac{N_n}{\sigma^2(n-\kappa_n)}\leq\frac{1}{8},
\end{equation}
then we have
\begin{eqnarray*}
\alpha_2(T_n)&=&\frac{n}{2\sigma^2(n-\kappa_n)}-L^2\gamma_n\cdot\left(\frac{n}{\sigma^2(n-\kappa_n)}\right)^2\cdot\frac{N_n}{n}\\
&=&\frac{n}{2\sigma^2(n-\kappa_n)}-\left(L^2\gamma_n\cdot\frac{N_n}{\sigma^2(n-\kappa_n)}\right)\cdot\frac{n}{\sigma^2(n-\kappa_n)}\\
&\geq&\frac{n}{2\sigma^2(n-\kappa_n)}-\frac{1}{8}\cdot\frac{n}{\sigma^2(n-\kappa_n)}\\
&=&\frac{3n}{8\sigma^2(n-\kappa_n)}.
\end{eqnarray*}
\end{rem}

With the next lemma, we will control the ``tails'' of the $\phi_n(\cdot)$ terms:

\begin{lem}[Tail estimates]\label{lem:techtails}
For all $n\in\mathbb{N}^*$ and $\Lambda\in\mathbb{R}_+$ such that the following technical condition is satisfied:
\begin{equation}\label{eq:nandA2}
\kappa_n\leq n-1,\quad\frac{\Lambda\sqrt{n}}{\sqrt{N_n}}-d\geq\sqrt{\frac{\left(1-2c\sigma^2\right)\kappa_n}{c}}\quad\text{and}\quad\frac{N_n}{2\sigma^2(n-\kappa_n)}\leq\frac{\delta_n}{2},
\end{equation}
the following holds.
For all $x\in\mathbb{R}$ such that $|x|\geq\Lambda$, we have
\[\phi_n(x)\leq\exp\left(-\alpha_3(T_n)\cdot x^2+\beta_3(T_n)\cdot|x|+\theta_3(T_n)\right),\]
where $\alpha_3(\cdot)$, $\beta_3(\cdot)$ and $\theta_3(\cdot)$ are given by
\begin{equation}\label{eq:alpha3}
\alpha_3(T_n)=\frac{cn}{2(n-\kappa_n)}-L^2\gamma_n\cdot\left(\frac{2cn}{n-\kappa_n}\right)^2\cdot\frac{N_n}{n},
\end{equation}
\begin{equation}\label{eq:beta3}
\beta_3(T_n)=\frac{cdn}{n-\kappa_n}\cdot\sqrt{\frac{N_n}{n}}+2dL^2\gamma_n\cdot\left(\frac{2cn}{n-\kappa_n}\right)^2\cdot\left(\frac{N_n}{n}\right)^{3/2},
\end{equation}
and
\begin{equation}\label{eq:theta3}
\theta_3(T_n)=\frac{CN_n}{n-\kappa_n}+\frac{cd^2N_n}{2(n-\kappa_n)}+d^2L^2\gamma_n\cdot\left(\frac{2cn}{n-\kappa_n}\right)^2\cdot\left(\frac{N_n}{n}\right)^2\\
+\frac{2}{\delta_n^2}\cdot\left(\frac{cN_n}{n-\kappa_n}\right)^2.
\end{equation}
\end{lem}

Again, even though $\alpha_3$, $\beta_3$ and $\theta_3$ have rather cumbersome expressions, they all simplify in the limit $n\to\infty$: if $N_n/n\rightarrow0$ and $\kappa_n/n\rightarrow0$, and $\gamma_n\lesssim1$ and $\delta_n\gtrsim1$, then $\alpha_3\rightarrow c/2$ and $\beta_3,\theta_3\rightarrow0$.

\begin{proof}[Proof of Lemma \ref{lem:techtails}]
Fix $n\in\mathbb{N}^*$ and $\Lambda\in\mathbb{R}_+$, and assume that the technical condition \eqref{eq:nandA2} is satisfied.
First, by Proposition \ref{prop:obs}, we have
\[\phi_n(x)=\exp\left(\frac{N_n\kappa_n}{2(n-\kappa_n)}\right)\cdot\mathbb{E}\left[\prod_{v\in\partial[T_n]_{\kappa_n}}\frac{P_{n-\kappa_n}(y_n(x)+S_v)}{q_{n-\kappa_n}(0)}\right]\quad\text{for all $x\in\mathbb{R}$.}\]
Now, by Lemma \ref{lem:lclt}, the following holds: for all $x\in\mathbb{R}$, and for all $v\in\partial[T_n]_{\kappa_n}$, we have
\[\frac{P_{n-\kappa_n}(y_n(x)+S_v)}{q_{n-\kappa_n}(0)}\leq\exp\left(-\frac{c(y_n(x)+S_v)^2}{n-\kappa_n}+\frac{C}{n-\kappa_n}\right).\]
Taking the product over $v\in\partial[T_n]_{\kappa_n}$, taking expectations and multiplying through by the ad-hoc exponential, we get
\[\phi_n(x)\leq\psi_n(x)\cdot\exp\left(\frac{CN_n}{n-\kappa_n}\right),\]
where
\[\psi_n(x)=\exp\left(\frac{N_n\kappa_n}{2(n-\kappa_n)}\right)\cdot\mathbb{E}\left[\prod_{v\in\partial[T_n]_{\kappa_n}}\exp\left(-\frac{c(y_n(x)+S_v)^2}{n-\kappa_n}\right)\right].\]
Then, by the calculation of Lemma \ref{lem:guide} with $\rho=c/(n-\kappa_n)$, we get
\begin{eqnarray*}
\psi_n(x)&=&\exp\left(\frac{N_n\kappa_n}{2(n-\kappa_n)}\right)\cdot\exp\left(-\frac{cN_ny_n(x)^2}{n-\kappa_n}\right)\\
&&\cdot\mathbb{E}\left[\exp\left(-\frac{cN_n\left(2y_n(x)\cdot\langle S\rangle_{\kappa_n}+\left(\left\langle S^2\right\rangle_{\kappa_n}-\sigma^2\kappa_n\right)\right)}{n-\kappa_n}\right)\right]\cdot\exp\left(-\frac{cN_n\sigma^2\kappa_n}{n-\kappa_n}\right)\\
&=&\exp\left(\frac{\left(1-2c\sigma^2\right)N_n\kappa_n}{2(n-\kappa_n)}\right)\cdot\exp\left(-\frac{cN_ny_n(x)^2}{n-\kappa_n}\right)\\
&&\cdot\mathbb{E}\left[\exp\left(-\frac{cN_n\left(2y_n(x)\cdot\langle S\rangle_{\kappa_n}+\left(\left\langle S^2\right\rangle_{\kappa_n}-\sigma^2\kappa_n\right)\right)}{n-\kappa_n}\right)\right].
\end{eqnarray*}
Next, by \eqref{eq:ynx} and by \eqref{eq:nandA2}, the following holds: for all $x\in\mathbb{R}$ such that $|x|\geq\Lambda$, we have
\[|y_n(x)|\geq\left|\frac{x\sqrt{n}}{\sqrt{N_n}}\right|-d\geq\frac{\Lambda\sqrt{n}}{\sqrt{N_n}}-d\geq\sqrt{\frac{\left(1-2c\sigma^2\right)\kappa_n}{c}}.\]
Thus, we get
\begin{eqnarray*}
\lefteqn{\exp\left(\frac{\left(1-2c\sigma^2\right)N_n\kappa_n}{2(n-\kappa_n)}\right)\cdot\exp\left(-\frac{cN_ny_n(x)^2}{n-\kappa_n}\right)}\\
&\leq&\exp\left(\frac{cN_ny_n(x)^2}{2(n-\kappa_n)}\right)\cdot\exp\left(-\frac{cN_ny_n(x)^2}{n-\kappa_n}\right)=\exp\left(-\frac{cN_ny_n(x)^2}{2(n-\kappa_n)}\right).
\end{eqnarray*}
Moreover, since $c<\left.1\middle/\left(2\sigma^2\right)\right.$, by \eqref{eq:nandA2} and by Corollary \ref{cor:zeta}, we have
\begin{eqnarray*}
\lefteqn{\mathbb{E}\left[\exp\left(-\frac{cN_n\left(2y_n(x)\cdot\langle S\rangle_{\kappa_n}+\left(\left\langle S^2\right\rangle_{\kappa_n}-\sigma^2\kappa_n\right)\right)}{n-\kappa_n}\right)\right]}\\
&\leq&\exp\left(L^2\gamma_n\cdot\left(\frac{2cN_ny_n(x)}{n-\kappa_n}\right)^2+\frac{2}{\delta_n^2}\cdot\left(\frac{cN_n}{n-\kappa_n}\right)^2\right).
\end{eqnarray*}
Thus, we obtain
\[\phi_n(x)\leq\exp\left(-\frac{cN_ny_n(x)^2}{2(n-\kappa_n)}+\frac{CN_n}{n-\kappa_n}+L^2\gamma_n\cdot\left(\frac{2cN_ny_n(x)}{n-\kappa_n}\right)^2+\frac{2}{\delta_n^2}\cdot\left(\frac{cN_n}{n-\kappa_n}\right)^2\right).\]
Plugging \eqref{eq:ynx2} into the above bound, and grouping together the terms in $x^2$, the terms in $|x|$ and the constant terms, we obtain the desired result (let us spare the reader these tedious details).
\end{proof}

\begin{rem}\label{rem:lemtech2}
For future reference, observe that if furthermore the technical condition \eqref{eq:nandA3} is satisfied, then we have, since $c<\left.1\middle/\left(2\sigma^2\right)\right.$:
\begin{eqnarray*}
\alpha_3(T_n)&=&\frac{cn}{2(n-\kappa_n)}-L^2\gamma_n\cdot\left(\frac{2cn}{n-\kappa_n}\right)^2\cdot\frac{N_n}{n}\\
&=&\frac{cn}{2(n-\kappa_n)}-\left(L^2\gamma_n\cdot\frac{2cN_n}{n-\kappa_n}\right)\cdot\frac{2cn}{n-\kappa_n}\\
&\geq&\frac{cn}{2(n-\kappa_n)}-\frac{1}{8}\cdot\frac{2cn}{n-\kappa_n}\\
&=&\frac{cn}{4(n-\kappa_n)}.
\end{eqnarray*}
\end{rem}

\subsection{Gaussian scaling limit}

Now armed with the previous lemmas, let us dive into the proof of Theorem \ref{thm:main}.
For each $n\in\mathbb{N}^*$, fix $p_n\in{]0,1]}$.
Then, conditionally on the random tree $\mathcal{T}_n^{p_n}$, let $H_n=(H_n(u))_{u\in\mathcal{T}_n^{p_n}}$ be a $\nu$-integer-valued random height function on $\mathcal{T}_n^{p_n}$ with zero boundary conditions on $\partial\mathcal{T}_n^{p_n}$.
In this subsection, we establish the Gaussian scaling limit in Theorem \ref{thm:main}: we show that if
\begin{equation}\label{eq:assdelocivrhf}
\frac{m^n\cdot p_n}{n}\longrightarrow0\quad\text{as $n\to\infty$,}
\end{equation}
then we have
\begin{equation}\label{eq:resdelocivrhf}
\mathbb{P}\left(\frac{\sqrt{\#\partial\mathcal{T}_n^{p_n}}}{\sqrt{n}}\cdot H_n(\varnothing)\in\cdot\,\middle|\,\mathcal{T}_n^{p_n}\right)\longrightarrow\mathcal{N}\left(0,\sigma^2\right)\quad\text{in probability as $n\to\infty$,}
\end{equation}
in the Polish space of Borel probability measures on $\mathbb{R}$ endowed with the topology of weak convergence.
Recall from Proposition \ref{prop:leavesTnp'} that \eqref{eq:assdelocivrhf} entails that $\left.\#\partial\mathcal{T}_n^{p_n}\middle/n\right.\rightarrow0$ in probability as $n\to\infty$.

\begin{proof}[Proof of the Gaussian scaling limit in Theorem \ref{thm:main}]
To lighten the notation, we set $N_n=\#\partial\mathcal{T}_n^{p_n}$.
To prove \eqref{eq:resdelocivrhf}, it suffices to show that for each $a<b\in\mathbb{R}$, we have
\[\mathbb{P}\left(\frac{\sqrt{N_n}}{\sqrt{n}}\cdot H_n(\varnothing)\in{]a,b]}\,\middle|\,\mathcal{T}_n^{p_n}\right)\longrightarrow\frac{1}{\sqrt{2\pi\sigma^2}}\cdot\int_a^b\exp\left(-\frac{x^2}{2\sigma^2}\right)\mathrm{d}x\]
in probability as $n\to\infty$ (see Remark \ref{rem:cvM1R} in the devoted Appendix \ref{sec:cvM1R}).
For all $y\in\mathbb{Z}$, we set
\[\Pi_n(y)=\sum_{\substack{h\in\mathbb{H}^0\left(\mathcal{T}_n^{p_n},\partial\mathcal{T}_n^{p_n}\right)\\h_\varnothing=y}}\prod_{(u,v)\in\vec{E}\left(\mathcal{T}_n^{p_n}\right)}\nu(\{h_v-h_u\}).\]
Let $K_n$ be the smallest integer $k\in\llbracket0,n\rrbracket$ such that every node $u\in\mathcal{T}_n^{p_n}$ has exactly one child.
On the good event $(K_n\leq n-1)$, we set
\[\Phi_n(x)=\exp\left(\frac{N_nK_n}{2(n-K_n)}\right)\cdot\frac{\Pi_n(Y_n(x))}{q_{n-K_n}(0)^{N_n}}\quad\text{for all $x\in\mathbb{R}$,}\]
where
\[Y_n(x)=rn+d\left\lfloor-\frac{rn}{d}+\frac{x\sqrt{n}}{d\sqrt{N_n}}\right\rfloor.\]
By \eqref{eq:coupdegras}, we have
\[\mathbb{P}\left(\frac{\sqrt{N_n}}{\sqrt{n}}\cdot H_n(\varnothing)\in{]a,b]}\,\middle|\,\mathcal{T}_n^{p_n}\right)=\frac{\int_{A_n}^{B_n}\Phi_n(x)\mathrm{d}x}{\int_{-\infty}^\infty\Phi_n(x)\mathrm{d}x},\]
where
\[A_n=\inf\left\{x\in\mathbb{R}:\frac{\sqrt{N_n}}{\sqrt{n}}\cdot Y_n(x)>a\right\}=\frac{\sqrt{N_n}}{\sqrt{n}}\cdot\left(rn+d\left(\left\lfloor-\frac{rn}{d}+\frac{a\sqrt{n}}{d\sqrt{N_n}}\right\rfloor+1\right)\right),\]
and
\[B_n=\sup\left\{x\in\mathbb{R}:\frac{\sqrt{N_n}}{\sqrt{n}}\cdot Y_n(x)\leq b\right\}=\frac{\sqrt{N_n}}{\sqrt{n}}\cdot\left(rn+d\left(\left\lfloor-\frac{rn}{d}+\frac{b\sqrt{n}}{d\sqrt{N_n}}\right\rfloor+1\right)\right).\]
From the explicit expressions of $A_n$ and $B_n$, let us only keep the following in mind: since $|(\lfloor y\rfloor+1)-y|\leq1$ for all $y\in\mathbb{R}$, we have
\begin{equation}\label{eq:AnBn}
|A_n-a|\leq\frac{d\sqrt{N_n}}{\sqrt{n}}\quad\text{and}\quad|B_n-b|\leq\frac{d\sqrt{N_n}}{\sqrt{n}}.
\end{equation}
Now, to complete the proof, it suffices to show that
\begin{equation}\label{eq:goalthm2}
\int_{A_n}^{B_n}\Phi_n(x)\mathrm{d}x\longrightarrow\int_a^b\exp\left(-\frac{x^2}{2\sigma^2}\right)\mathrm{d}x\quad\text{and}\quad\int_{-\infty}^\infty\Phi_n(x)\mathrm{d}x\longrightarrow\sqrt{2\pi\sigma^2}
\end{equation}
in probability as $n\to\infty$.
To this end, we will use Lemma \ref{lem:techbulk} and Lemma \ref{lem:techtails} to estimate $\Phi_n(\cdot)$.
Set
\[\Gamma_n=1\vee\left(\frac{1}{N_n}\cdot\sum_{k=1}^{K_n}k^3\cdot\max_{\substack{u\in\mathcal{T}_n^{p_n}\\|u|=k}}\#\left\{v\in\partial\mathcal{T}_n^{p_n}:v\succeq u\right\}\right),\]
and let
\[\Delta_n=\frac{1}{32L^2\sqrt{2\Gamma_n}}.\]
Then, let $\eta\in\mathbb{R}_+^*$, $C\in\mathbb{R}_+^*$ and $c\in\left]0,\left.1\middle/\left(2\sigma^2\right)\right.\right[$ be as in the statement of Lemma \ref{lem:lclt}.
Moreover, set $\Lambda_n=\ln(n)$, and consider the good event (in view of the technical conditions \eqref{eq:nandA1} and \eqref{eq:nandA2})
\begin{eqnarray*}
E_n&=&\left(K_n\leq n-1\,;\,\frac{\Lambda_n\sqrt{n}}{\sqrt{N_n}}+K_n+d\leq\eta(n-K_n)\,;\,\frac{N_n}{2\sigma^2(n-K_n)}\leq\frac{\Delta_n}{2}\,;\right.\\
&&\left.\,\frac{\Lambda_n\sqrt{n}}{\sqrt{N_n}}-d\geq\sqrt{\frac{\left(1-2c\sigma^2\right)K_n}{c}}\right).
\end{eqnarray*}
First, we claim that $\mathbb{P}(E_n)\rightarrow1$ as $n\to\infty$.
Indeed, by Lemma \ref{lem:cut-height}, there exists a constant $C'\in\mathbb{R}_+^*$ such that
\[\mathbb{P}(K_n>2\log_m(n))\leq C'\cdot\frac{(m^n\cdot p_n)\vee(m^n\cdot p_n)^2}{m^{2\log_m(n)}}=C'\cdot\frac{(m^n\cdot p_n)\vee(m^n\cdot p_n)^2}{n^2}.\]
Thus, by \eqref{eq:assdelocivrhf}, we get that
\begin{equation}\label{eq:fact1}
\mathbb{P}(K_n>2\log_m(n))\longrightarrow0\quad\text{as $n\to\infty$.}
\end{equation}
Moreover, by \eqref{eq:assdelocivrhf} and by Proposition \ref{prop:leavesTnp'}, we have
\begin{equation}\label{eq:fact2}
\frac{N_n}{n}\longrightarrow0\quad\text{in probability as $n\to\infty$.}
\end{equation}
Finally, by Lemma \ref{lem:cut-height}, the sequence of random variables $\left(\Gamma_n\right)_{n\in\mathbb{N}^*}$ is tight: we have
\begin{equation}\label{eq:fact3}
\text{$\varlimsup_{n\to\infty}\mathbb{P}(\Gamma_n\geq\gamma)\longrightarrow0$ as $\gamma\to\infty$, \quad hence \quad $\varlimsup_{n\to\infty}\mathbb{P}(\Delta_n\leq\delta)\longrightarrow0$ as $\delta\to0^+$.}
\end{equation}
Combining \eqref{eq:fact1}, \eqref{eq:fact2} and \eqref{eq:fact3}, we get (recall that $\Lambda_n=\ln(n)$):
\[K_n\ll n,\quad\frac{\Lambda_n\sqrt{n}}{\sqrt{N_n}}\ll n,\quad\frac{N_n}{n}\ll\Delta_n,\quad\text{and}\quad\sqrt{K_n}\ll\frac{\Lambda_n\sqrt{n}}{\sqrt{N_n}}\quad\text{in probability as $n\to\infty$.}\]
This entails that $E_n$ is realised with high probability as $n\to\infty$, as desired.
Now, by construction, on the good event $E_n$, the technical conditions \eqref{eq:nandA1} and \eqref{eq:nandA2} are satisfied.
Thus, by Lemma \ref{lem:techbulk}, the following holds: for all $x\in[-\Lambda_n,\Lambda_n]$, we have
\[\Phi_n(x)\geq\exp\left(-\alpha_1\left(\mathcal{T}_n^{p_n}\right)\cdot x^2-\beta_1\left(\mathcal{T}_n^{p_n}\right)\cdot|x|-\theta_1\left(\mathcal{T}_n^{p_n},\Lambda_n\right)\right),\]
where $\alpha_1(\cdot)$, $\beta_1(\cdot)$ and $\theta_1(\cdot,\cdot)$ are given by \eqref{eq:alpha1}, \eqref{eq:beta1} and \eqref{eq:theta1}, and
\[\Phi_n(x)\leq\exp\left(-\alpha_2\left(\mathcal{T}_n^{p_n}\right)\cdot x^2+\beta_2\left(\mathcal{T}_n^{p_n}\right)\cdot|x|+\theta_2\left(\mathcal{T}_n^{p_n},\Lambda_n\right)\right),\]
where $\alpha_2(\cdot)$, $\beta_2(\cdot)$ and $\theta_2(\cdot,\cdot)$ are given by \eqref{eq:alpha2}, \eqref{eq:beta2} and \eqref{eq:theta2}.
Moreover, by Lemma \ref{lem:techtails}, the following holds: for all $x\in\mathbb{R}$ such that $|x|\geq\Lambda_n$, we have
\[\Phi_n(x)\leq\exp\left(-\alpha_3\left(\mathcal{T}_n^{p_n}\right)\cdot x^2+\beta_3\left(\mathcal{T}_n^{p_n}\right)\cdot|x|+\theta_3\left(\mathcal{T}_n^{p_n}\right)\right),\]
where $\alpha_3(\cdot)$, $\beta_3(\cdot)$ and $\theta_3(\cdot)$ are given by \eqref{eq:alpha3}, \eqref{eq:beta3} and \eqref{eq:theta3}.
Now, combining \eqref{eq:fact1}, \eqref{eq:fact2} and \eqref{eq:fact3}, we get (recall that $\Lambda_n=\ln(n)$):
\[K_n\ll n,\quad\frac{N_n}{n}\ll1,\quad\frac{N_n}{n^3}\cdot\left(\frac{\Lambda_n\sqrt{n}}{\sqrt{N_n}}+K_n+d\right)^4\ll1,\quad\text{and}\quad\Gamma_n,\frac{1}{\Delta_n^2}\lesssim1\]
in probability as $n\to\infty$.
This entails that for each $i\in\{1,2\}$, we have
\[\alpha_i\left(\mathcal{T}_n^{p_n}\right)\longrightarrow\frac{1}{2\sigma^2}\quad\text{and}\quad\beta_i\left(\mathcal{T}_n^{p_n}\right),\theta_i\left(\mathcal{T}_n^{p_n},\Lambda_n\right)\longrightarrow0\quad\text{in probability as $n\to\infty$,}\]
and
\[\alpha_3\left(\mathcal{T}_n^{p_n}\right)\longrightarrow\frac{c}{2}\quad\text{and}\quad\beta_3\left(\mathcal{T}_n^{p_n}\right),\theta_3\left(\mathcal{T}_n^{p_n}\right)\longrightarrow0\quad\text{in probability as $n\to\infty$.}\]
Finally, by \eqref{eq:AnBn} and by \eqref{eq:fact2}, we have $A_n\rightarrow a$ and $B_n\rightarrow b$ in probability as $n\to\infty$.
Lastly, recall that $\Lambda_n=\ln(n)\rightarrow\infty$ as $n\to\infty$.
Altogether, by the dominated convergence theorem, we obtain that for each $\varepsilon\in\mathbb{R}_+^*$, we have
\[\mathbb{P}\left(\left|\int_{A_n}^{B_n}\Phi_n(x)\mathrm{d}x-\int_a^b\exp\left(-\frac{x^2}{2\sigma^2}\right)\mathrm{d}x\right|\geq\varepsilon\,;\,E_n\right)\longrightarrow0\quad\text{as $n\to\infty$,}\]
and
\[\mathbb{P}\left(\left|\int_{-\infty}^\infty\Phi_n(x)\mathrm{d}x-\int_{-\infty}^\infty\exp\left(-\frac{x^2}{2\sigma^2}\right)\mathrm{d}x\right|\geq\varepsilon\,;\,E_n\right)\longrightarrow0\quad\text{as $n\to\infty$.}\]
Since $\mathbb{P}(E_n)\rightarrow1$ as $n\to\infty$, this yields \eqref{eq:goalthm2}, which concludes the proof.
\end{proof}

\subsection{Localisation}\label{subsec:locivrhf}

In this subsection, we prove the localisation in Theorem \ref{thm:main}.
From now on, we assume that $\nu$ is the uniform distribution on $\{-1,1\}$.
In that case, the absolute value of $\nu$-integer-valued random height functions has the following monotonicity property:

\begin{prop}\label{prop:monotonicity}
Let $T_n$ be a finite rooted tree with height $n\in\mathbb{N}^*$ and with no internal leaves, and let $B_n'\subset B_n\subset\partial T_n$ be nested non-empty subsets of leaves.
Then, let $H_n=(H_n(u))_{u\in T_n}$ be a $\nu$-integer-valued random height function on $T_n$ with zero boundary conditions on $B_n$, and let $H_n'=(H_n'(u))_{u\in T_n}$ be a $\nu$-integer-valued random height function on $T_n$ with zero boundary conditions on $B_n'$.
For all non-decreasing functions ${\Phi:\mathbb{N}^{T_n}\rightarrow\mathbb{R}_+}$ (i.e, such that $\Phi(f)\leq\Phi(f')$ whenever ${f_u\leq f'_u}$ for all $u\in T_n$), we have $\mathbb{E}[\Phi(|H_n|)]\leq\mathbb{E}[\Phi(|H_n'|)]$.
In particular, we get that $\mathbb{V}(H_n(\varnothing))\leq\mathbb{V}(H_n'(\varnothing))$.
\end{prop}

This is not a new observation: see \cite[Proposition 2.2]{logvariance} for a variant with more general boundary conditions, and see \cite[Section 2]{randomgraphhom} for other variants.
For the convenience of the reader, we provide a proof of Proposition \ref{prop:monotonicity} in Appendix \ref{sec:monotonicity}.

\begin{rem}
\begin{itemize}
\item This monotonicity property fails for general distributions $\nu$: fix $p\in{]1/3,1/2[}$, and let $\nu$ be the discrete probability measure on $\mathbb{Z}$ given by ${\nu(\{1\}),\nu(\{-1\})=p}$ and $\nu(\{0\})=1-2p$.
Then, let $T_1$ be the finite rooted tree with height $1$ in which the root $\varnothing$ has two children $1$ and $2$, and set $B_1=\{1,2\}$ and $B_1'=\{1\}$.
Next, let $H_1=(H_1(u))_{u\in T_1}$ be a $\nu$-integer-valued random height function on $T_1$ with zero boundary conditions on $B_1$, and let $H_1'=(H_1'(u))_{u\in T_1}$ be a $\nu$-integer-valued random height function on $T_1$ with zero boundary conditions on $B_1'$.
We have
\[\mathbb{V}(H_1(\varnothing))=\frac{2p^2}{2p^2+(1-2p)^2}>2p=\mathbb{V}(H_1'(\varnothing)).\]

\item Under different assumptions on the distribution $\nu$, the FKG lattice condition proved by Lammers and Ott \cite[Theorem 2.8]{pietseb} entails that the absolute value of $\nu$-integer-valued random height functions has the same monotonicity property as in Proposition \ref{prop:monotonicity}.
Namely, their result applies to distributions $\nu$ such that $\nu(\{x\})=e^{-V(x)}$ for all $x\in\mathbb{Z}$, where $V:\mathbb{Z}\rightarrow\mathbb{R}$ is convex, symmetric, and ``super-Gaussian'', meaning that
\[x\in\mathbb{N}\mapsto V(x+1)-2V(x)+V(x-1)\]
is non-increasing.
\end{itemize}
\end{rem}

Back to our uniform distribution $\nu$ on $\{-1,1\}$ and equipped with this monotonicity property, we are now ready to prove the localisation in Theorem \ref{thm:main}.
For each $n\in\mathbb{N}^*$, fix $p_n\in{]0,1]}$, and conditionally on the random tree $\mathcal{T}_n^{p_n}$, let $H_n=(H_n(u))_{u\in\mathcal{T}_n^{p_n}}$ be a $\nu$-integer-valued random height function on $\mathcal{T}_n^{p_n}$ with zero boundary conditions on $\partial\mathcal{T}_n^{p_n}$.
We show below that if
\begin{equation}\label{eq:asslocivrhf}
\inf_{n\geq1}\frac{m^n\cdot p_n}{n}>0,
\end{equation}
then the sequence of random variables $(\mathbb{V}(H_n(\varnothing)\,|\,\mathcal{T}_n^{p_n}))_{n\in\mathbb{N}^*}$ is tight.

\begin{proof}[Proof of the localisation in Theorem \ref{thm:main}]
Let
\[\varepsilon_0=\inf_{n\geq1}\frac{m^n\cdot p_n}{n}>0.\]
For each $\varepsilon\in{]0,\varepsilon_0]}$, we set
\[p_n'(\varepsilon)=\frac{\varepsilon\cdot n}{m^n}\leq p_n\quad\text{for all $n\in\mathbb{N}^*$,}\]
and we denote by $V_n^\varepsilon$ the variance of the height at the root of the tree for a $\nu$-integer-valued random height function on $\mathcal{T}_n^{p_n'(\varepsilon)}$ with zero boundary conditions on $\partial\mathcal{T}_n^{p_n'(\varepsilon)}$.
We claim that for all $A\in\mathbb{R}_+^*$, we have
\begin{equation}\label{eq:comparison}
\varliminf_{n\to\infty}\mathbb{P}\left(\mathbb{V}\left(H_n(\varnothing)\,\middle|\,\mathcal{T}_n^{p_n}\right)\leq A\right)\geq\varliminf_{n\to\infty}\mathbb{P}\left(V_n^\varepsilon\leq A\right).
\end{equation}
In fact, let us show that\footnote{At first, we thought that $\mathbb{V}\left(H_n(\varnothing)\,\middle|\,\mathcal{T}_n^{p_n}\right)$ would be stochastically dominated by $V_n^\varepsilon$.
However, because of the annealed conditioning in the definition of the pruned tree $\mathcal{T}_n^p$, it is not true in general that if $p,p'\in{]0,1]}$ are such that $p'\leq p$, then there is a coupling of $\mathcal{T}_n^p$ and $\mathcal{T}_n^{p'}$ such that $\mathcal{T}_n^{p'}$ is a subtree of $\mathcal{T}_n^p$.
Indeed, this reduces to asking whether $\mu_p$ is stochastically increasing in $p\in{]0,1]}$, where $\mu_p$ denotes the distribution of a random variable with distribution $\mathrm{Binomial}(\xi,p)$ conditionally on $\xi$, conditioned on being strictly positive:
\[\mu_p(\{k\})=\frac{\mathbb{E}\left[\binom{\xi}{k}p^k(1-p)^{\xi-k}\right]}{1-G(1-p)}\quad\text{for all $k\in\mathbb{N}^*$.}\]
Although $\mu_p$ is indeed stochastically increasing in $p$ if $\xi=m$ almost surely for some $m\in\mathbb{N}^*$ (then, for all ${p\leq q\in{]0,1]}}$, the Radon--Nikodym derivative
\[k\in\llbracket1,m\rrbracket\longmapsto\frac{\mu_q(\{k\})}{\mu_p(\{k\})}=\left(\frac{1-q}{1-p}\right)^m\cdot\left(\frac{q/(1-q)}{p/(1-p)}\right)^k\]
of $\mu_q$ with respect to $\mu_p$ is non-decreasing, hence $\mu_p$ is stochastically dominated by $\mu_q$, by virtue of the Harris--FKG inequality), or if $\xi\in\{1,2\}$ almost surely (then, it suffices to check that
\[p\in{]0,1]}\longmapsto\mu_p(\{2\})=\frac{p^2\cdot\mu(\{2\})}{1-G(1-p)}\]
is non-decreasing, and by the convexity of $G$, the function ${p\in{]0,1]}\mapsto(1-G(1-p))/p}$ is non-increasing), this is false in general.
For instance, if $\mu$ is given by ${\mu(\{1\}),\mu(\{3\})=1/2}$, then the function
\[p\in{]0,1]}\longmapsto\mu_p(\{1\})=\frac{p\cdot G'(1-p)}{1-G(1-p)}=\frac{p\left(1+3(1-p)^2\right)}{1+p-(1-p)^3}\]
is not non-increasing.}
\begin{equation}\label{eq:comparison'}
\mathbb{P}\left(\mathbb{V}\left(H_n(\varnothing)\,\middle|\,\mathcal{T}_n^{p_n}\right)\leq A\right)\geq\frac{1-G^n(1-p_n'(\varepsilon))}{1-G^n(1-p_n)}\cdot\mathbb{P}\left(V_n^\varepsilon\leq A\right).
\end{equation}
To this end, let $(U_n(v))_{v\in\partial[\mathbb{U}]_n}$ be a collection of independent random variables with uniform distribution on $[0,1]$, independent of $\mathcal{T}$.
Then, set
\[\mathcal{B}_n=\{v\in\partial[\mathcal{T}]_n:U_n(v)\leq p_n\}\quad\text{and}\quad\mathcal{B}_n'=\{v\in\partial[\mathcal{T}]_n:U_n(v)\leq p_n'(\varepsilon)\},\]
and let us introduce the following temporary random variables:
\begin{itemize}
\item On the event $(\mathcal{B}_n\neq\emptyset)$, let $V_n$ be the variance of the height at the root of the tree for a $\nu$-integer-valued random height function on $[\mathcal{T}]_n$ with zero boundary conditions on $\mathcal{B}_n$,
\item On the event $(\mathcal{B}_n'\neq\emptyset)$, let $V_n'$ be the variance of the height at the root of the tree for a $\nu$-integer-valued random height function on $[\mathcal{T}]_n$ with zero boundary conditions on $\mathcal{B}_n'$.
\end{itemize}
By the definition of the pruned trees $\mathcal{T}_n^{p_n}$ and $\mathcal{T}_n^{p_n'(\varepsilon)}$, and by Proposition \ref{prop:prune}, we have
\[\mathbb{P}\left(\mathbb{V}\left(H_n(\varnothing)\,\middle|\,\mathcal{T}_n^{p_n}\right)\leq A\right)=\frac{\mathbb{P}(V_n\leq A\,;\,\mathcal{B}_n\neq\emptyset)}{1-G^n(1-p_n)},\]
and
\[\mathbb{P}\left(V_n^\varepsilon\leq A\right)=\frac{\mathbb{P}(V_n'\leq A\,;\,\mathcal{B}_n'\neq\emptyset)}{1-G^n(1-p_n'(\varepsilon))}.\]
Then, since $\mathcal{B}_n'\subset\mathcal{B}_n$, on the event $(\mathcal{B}_n'\neq\emptyset)$, by Proposition \ref{prop:monotonicity}, we have $V_n\leq V_n'$.
In particular, we get that
\[\mathbb{P}\left(V_n'\leq A\,;\,\mathcal{B}_n'\neq\emptyset\right)\leq\mathbb{P}(V_n\leq A\,;\,\mathcal{B}_n\neq\emptyset),\]
and it follows that
\begin{eqnarray*}
\mathbb{P}\left(\mathbb{V}\left(H_n(\varnothing)\,\middle|\,\mathcal{T}_n^{p_n}\right)\leq A\right)&=&\frac{\mathbb{P}(V_n\leq A\,;\,\mathcal{B}_n\neq\emptyset)}{1-G^n(1-p_n)}\\
&\geq&\frac{\mathbb{P}\left(V_n'\leq A\,;\,\mathcal{B}_n'\neq\emptyset\right)}{1-G^n(1-p_n)}\\
&=&\frac{1-G^n(1-p_n'(\varepsilon))}{1-G^n(1-p_n)}\cdot\mathbb{P}\left(V_n^\varepsilon\leq A\right),
\end{eqnarray*}
as claimed in \eqref{eq:comparison'}.
Then, by Proposition \ref{prop:leavesTnp'} and by \eqref{eq:asslocivrhf}, we have $1-G^n(1-p_n'(\varepsilon))\rightarrow1$ and $1-G^n(1-p_n)\rightarrow1$ as $n\to\infty$.
This yields the desired inequality \eqref{eq:comparison}.
Now, to complete the proof, it suffices to show that for each $\zeta\in{]0,1[}$, there exists $\varepsilon\in{]0,\varepsilon_0]}$ and $A\in\mathbb{R}_+^*$ such that
\begin{equation}\label{eq:locivrhfgoal}
\varliminf_{n\to\infty}\mathbb{P}\left(V_n^\varepsilon\leq A\right)\geq1-\zeta.
\end{equation}
To lighten the notation, we set $N_n^\varepsilon:=\#\partial\mathcal{T}_n^{p_n'(\varepsilon)}$.
For all $y\in\mathbb{Z}$, let
\[\Pi_n^\varepsilon(y)=\sum_{\substack{h\in\mathbb{H}^0\left(\mathcal{T}_n^{p_n'(\varepsilon)},\partial\mathcal{T}_n^{p_n'(\varepsilon)}\right)\\h_\varnothing=y}}\prod_{(u,v)\in\vec{E}\left(\mathcal{T}_n^{p_n'(\varepsilon)}\right)}\nu(\{h_v-h_u\}).\]
Then, let $K_n^\varepsilon$ be the smallest integer $k\in\llbracket0,n\rrbracket$ such that every node $u\in\mathcal{T}_n^{p_n'(\varepsilon)}$ has exactly one child.
On the good event $\left(K_n^\varepsilon\leq n-1\right)$, we set
\[\Phi_n^\varepsilon(x)=\exp\left(\frac{N_n^\varepsilon K_n^\varepsilon}{2\left(n-K_n^\varepsilon\right)}\right)\cdot\frac{\Pi_n^\varepsilon\left(Y_n^\varepsilon(x)\right)}{q_{n-K_n^\varepsilon}(0)^{N_n^\varepsilon}}\quad\text{for all $x\in\mathbb{R}$,}\]
where
\[Y_n^\varepsilon(x)=rn+d\left\lfloor-\frac{rn}{d}+\frac{x\sqrt{n}}{d\sqrt{N_n^\varepsilon}}\right\rfloor.\]
By \eqref{eq:coupdegras}, we have
\[V_n^\varepsilon=\frac{\int_{-\infty}^{\infty}Y_n^\varepsilon(x)^2\cdot\Phi_n^\varepsilon(x)\mathrm{d}x}{\int_{-\infty}^\infty\Phi_n^\varepsilon(x)\mathrm{d}x}.\]
By \eqref{eq:ynx}, using the inequality $(a+b)^2\leq2\left(a^2+b^2\right)$ for all $a,b\in\mathbb{R}$, we get
\[V_n^\varepsilon\leq\frac{\int_{-\infty}^\infty2\left(\frac{x^2n}{N_n^\varepsilon}+d^2\right)\cdot\Phi_n^\varepsilon(x)\mathrm{d}x}{\int_{-\infty}^\infty\Phi_n^\varepsilon(x)\mathrm{d}x}=2\left(\frac{n}{N_n^\varepsilon}\cdot\frac{\int_{-\infty}^\infty x^2\cdot\Phi_n^\varepsilon(x)\mathrm{d}x}{\int_{-\infty}^\infty\Phi_n^\varepsilon(x)\mathrm{d}x}+d^2\right).\]
Now, let us estimate the right-hand side.
First, by Proposition \ref{prop:leavesTnp'}, we have
\begin{equation}\label{eq:fact2'}
\frac{N_n^\varepsilon}{n}\longrightarrow\varepsilon\cdot W\quad\text{in distribution as $n\to\infty$,}
\end{equation}
where $W>0$ almost surely.
Next, we will use Lemma \ref{lem:techbulk} and Lemma \ref{lem:techtails} to estimate $\Phi_n^\varepsilon(\cdot)$.
Set
\[\Gamma_n^\varepsilon=1\vee\left(\frac{1}{N_n^\varepsilon}\cdot\sum_{k=1}^{K_n^\varepsilon}k^3\cdot\max_{\substack{u\in\mathcal{T}_n^{p_n'(\varepsilon)}\\|u|=k}}\#\left\{v\in\partial\mathcal{T}_n^{p_n'(\varepsilon)}:v\succeq u\right\}\right),\]
and let
\[\Delta_n^\varepsilon=\frac{1}{32L^2\sqrt{2\Gamma_n^\varepsilon}}.\]
Then, let $\eta\in\mathbb{R}_+^*$, $C\in\mathbb{R}_+^*$ and $c\in\left]0,\left.1\middle/\left(2\sigma^2\right)\right.\right[$ be as in the statement of Lemma \ref{lem:lclt}.
Moreover, set $\Lambda_n=\ln(n)$, and consider the good events (in view of the technical conditions \eqref{eq:nandA1}, \eqref{eq:nandA2} and \eqref{eq:nandA3}):
\begin{eqnarray*}
E_n^\varepsilon&=&\left(K_n^\varepsilon\leq n-1\,;\,\frac{\Lambda_n\sqrt{n}}{\sqrt{N_n^\varepsilon}}+K_n^\varepsilon+d\leq\eta\left(n-K_n^\varepsilon\right)\,;\,\frac{\Lambda_n\sqrt{n}}{\sqrt{N_n^\varepsilon}}-d\geq\sqrt{\frac{\left(1-2c\sigma^2\right)K_n^\varepsilon}{c}}\right),
\end{eqnarray*}
and
\begin{eqnarray*}
F_n^\varepsilon&=&E_n^\varepsilon\cap\left(\frac{N_n^\varepsilon}{2\sigma^2\left(n-K_n^\varepsilon\right)}\leq\frac{\Delta_n^\varepsilon}{2}\,;\,L^2\Gamma_n^\varepsilon\cdot\frac{N_n^\varepsilon}{\sigma^2\left(n-K_n^\varepsilon\right)}\leq\frac{1}{8}\right).
\end{eqnarray*}
We claim that
\begin{equation}\label{eq:Fnepsilon}
\varliminf_{n\to\infty}\mathbb{P}\left(F_n^\varepsilon\right)\longrightarrow1\quad\text{as $\varepsilon\to0$.}
\end{equation}
First, let us show that for each $\varepsilon\in{]0,\varepsilon_0]}$, we have $\mathbb{P}\left(E_n^\varepsilon\right)\rightarrow1$ as $n\to\infty$.
By Lemma \ref{lem:cut-height}, there exists a constant $C'\in\mathbb{R}_+^*$ such that
\begin{equation}\label{eq:fact1'}
\mathbb{P}\left(K_n^\varepsilon>3\log_m(n)\right)\leq C'\cdot\frac{(m^n\cdot p_n'(\varepsilon))\vee(m^n\cdot p_n'(\varepsilon))^2}{m^{3\log_m(n)}}=C'\cdot\frac{(\varepsilon\cdot n)\vee(\varepsilon\cdot n)^2}{n^3}\underset{n\to\infty}{\longrightarrow}0.
\end{equation}
Combining \eqref{eq:fact1'} and \eqref{eq:fact2'}, we get (recall that $\Lambda_n=\ln(n)$):
\[K_n^\varepsilon\ll n,\quad\frac{\Lambda_n\sqrt{n}}{\sqrt{N_n^\varepsilon}}\ll n,\quad\text{and}\quad\sqrt{K_n^\varepsilon}\ll\frac{\Lambda_n\sqrt{n}}{\sqrt{N_n^\varepsilon}}\quad\text{in probability as $n\to\infty$.}\]
This entails that $\mathbb{P}(E_n^\varepsilon)\rightarrow1$ as $n\to\infty$, as desired.
Then, to prove \eqref{eq:Fnepsilon}, it suffices to show that $\varlimsup_{n\to\infty}\mathbb{P}(E_n^\varepsilon\setminus F_n^\varepsilon)\rightarrow0$ as $\varepsilon\to0$.
For this purpose, fix $\zeta\in{]0,1[}$.
For all $\delta\in\mathbb{R}_+^*$ and $\gamma\in\mathbb{R}_+^*$, we have
\begin{eqnarray*}
\mathbb{P}\left(E_n^\varepsilon\setminus F_n^\varepsilon\right)&\leq&\mathbb{P}\left(\frac{N_n^\varepsilon}{2\sigma^2\left(n-K_n^\varepsilon\right)}\geq\frac{\Delta_n^\varepsilon}{2}\right)+\mathbb{P}\left(L^2\Gamma_n^\varepsilon\cdot\frac{N_n^\varepsilon}{\sigma^2\left(n-K_n^\varepsilon\right)}\geq\frac{1}{8}\right)\\
&\leq&\mathbb{P}(\Delta_n^\varepsilon\leq\delta)+\mathbb{P}\left(\frac{N_n^\varepsilon}{2\sigma^2\left(n-K_n^\varepsilon\right)}\geq\frac{\delta}{2}\right)\\
&&+\mathbb{P}\left(\Gamma_n^\varepsilon\geq\gamma\right)+\mathbb{P}\left(L^2\gamma\cdot\frac{N_n^\varepsilon}{\sigma^2\left(n-K_n^\varepsilon\right)}\geq\frac{1}{8}\right)
\end{eqnarray*}
Now, by Lemma \ref{lem:cut-height}, the collection of random variables $\left(\Gamma_n^\varepsilon\right)_{\substack{n\in\mathbb{N}^*\\\varepsilon\in{]0,\varepsilon_0]}}}$ is tight: we have
\begin{equation}\label{eq:fact3'}
\text{$\sup_{\substack{n\in\mathbb{N}^*\\\varepsilon\in{]0,\varepsilon_0]}}}\mathbb{P}\left(\Gamma_n^\varepsilon\geq\gamma\right)\longrightarrow0$ as $\gamma\to\infty$, \quad hence \quad $\sup_{\substack{n\in\mathbb{N}^*\\\varepsilon\in{]0,\varepsilon_0]}}}\mathbb{P}\left(\Delta_n^\varepsilon\leq\delta\right)\longrightarrow0$ as $\delta\to0^+$.}
\end{equation}
In particular, we can fix $\delta\in\mathbb{R}_+^*$ and $\gamma\in\mathbb{R}_+^*$ such that
\[\sup_{\substack{n\in\mathbb{N}^*\\\varepsilon\in{]0,\varepsilon_0]}}}\mathbb{P}\left(\Delta_n^\varepsilon\leq\delta\right)\leq\zeta\quad\text{and}\quad\sup_{\substack{n\in\mathbb{N}^*\\\varepsilon\in{]0,\varepsilon_0]}}}\mathbb{P}\left(\Gamma_n^\varepsilon\geq\gamma\right)\leq\zeta.\]
Next, combining \eqref{eq:fact1'} and \eqref{eq:fact2'}, by Slutsky's lemma, we get that
\[\frac{N_n^\varepsilon}{n-K_n^\varepsilon}\longrightarrow\varepsilon\cdot W\quad\text{in distribution as $n\to\infty$.}\]
By the Portmanteau theorem, it follows that
\[\varlimsup_{n\to\infty}\mathbb{P}\left(\frac{N_n^\varepsilon}{2\sigma^2\left(n-K_n^\varepsilon\right)}\geq\frac{\delta}{2}\right)\leq\mathbb{P}\left(\frac{\varepsilon\cdot W}{2\sigma^2}\geq\frac{\delta}{2}\right),\]
and
\[\varlimsup_{n\to\infty}\mathbb{P}\left(L^2\gamma\cdot\frac{N_n^\varepsilon}{\sigma^2\left(n-K_n^\varepsilon\right)}\geq\frac{1}{8}\right)\leq\mathbb{P}\left(L^2\gamma\cdot\frac{\varepsilon\cdot W}{\sigma^2}\geq\frac{1}{8}\right).\]
Thus, we obtain
\[\varlimsup_{n\to\infty}\mathbb{P}(E_n^\varepsilon\setminus F_n^\varepsilon)\leq\zeta+\mathbb{P}\left(\frac{\varepsilon\cdot W}{2\sigma^2}\geq\frac{\delta}{2}\right)+\zeta+\mathbb{P}\left(L^2\gamma\cdot\frac{\varepsilon\cdot W}{\sigma^2}\geq\frac{1}{8}\right).\]
Finally, since $W>0$ almost surely, we have
\[\mathbb{P}\left(\frac{\varepsilon\cdot W}{2\sigma^2}\geq\frac{\delta}{2}\right)\longrightarrow0\quad\text{and}\quad\mathbb{P}\left(L^2\gamma\cdot\frac{\varepsilon\cdot W}{\sigma^2}\geq\frac{1}{8}\right)\longrightarrow0\quad\text{as $\varepsilon\to0$.}\]
Thus, we get
\[\varlimsup_{\varepsilon\to0}\left(\varlimsup_{n\to\infty}\mathbb{P}(E_n^\varepsilon\setminus F_n^\varepsilon)\right)\leq\zeta+0+\zeta+0=2\zeta.\]
Since $\zeta$ is arbitrary, this shows that $\varlimsup_{n\to\infty}\mathbb{P}(E_n^\varepsilon\setminus F_n^\varepsilon)\rightarrow0$ as $\varepsilon\to0$, which completes the proof of \eqref{eq:Fnepsilon}.

Now, to prove \eqref{eq:locivrhfgoal}, fix $\zeta\in{]0,1[}$, and fix $\varepsilon\in{]0,\varepsilon_0]}$ such that $\varliminf_{n\to\infty}\mathbb{P}(F_n^\varepsilon)\geq1-\zeta$.
On the good event $F_n^\varepsilon$, by construction, the technical conditions \eqref{eq:nandA1}, \eqref{eq:nandA2} and \eqref{eq:nandA3} are satisfied.
Thus, by Lemma \ref{lem:techbulk}, the following holds: for all $x\in[-\Lambda_n,\Lambda_n]$, we have
\[\Phi_n^\varepsilon(x)\geq\exp\left(-\alpha_1\left(\mathcal{T}_n^{p_n'(\varepsilon)}\right)\cdot x^2-\beta_1\left(\mathcal{T}_n^{p_n'(\varepsilon)}\right)\cdot|x|-\theta_1\left(\mathcal{T}_n^{p_n'(\varepsilon)},\Lambda_n\right)\right),\]
where $\alpha_1(\cdot)$, $\beta_1(\cdot)$ and $\theta_1(\cdot,\cdot)$ are given by \eqref{eq:alpha1}, \eqref{eq:beta1} and \eqref{eq:theta1}, and
\[\Phi_n^\varepsilon(x)\leq\exp\left(-\alpha_2\left(\mathcal{T}_n^{p_n'(\varepsilon)}\right)\cdot x^2+\beta_2\left(\mathcal{T}_n^{p_n'(\varepsilon)}\right)\cdot|x|+\theta_2\left(\mathcal{T}_n^{p_n'(\varepsilon)},\Lambda_n\right)\right),\]
where $\alpha_2(\cdot)$, $\beta_2(\cdot)$ and $\theta_2(\cdot,\cdot)$ are given by \eqref{eq:alpha2}, \eqref{eq:beta2} and \eqref{eq:theta2}.
Furthermore, by Remark \ref{rem:lemtech1}, since the technical condition \eqref{eq:nandA3} is satisfied, we have
\[\alpha_2\left(\mathcal{T}_n^{p_n'(\varepsilon)}\right)\geq\frac{3n}{8\sigma^2(n-K_n^\varepsilon)}.\]
Moreover, by Lemma \ref{lem:techtails}, the following holds: for all $x\in\mathbb{R}$ such that $|x|\geq\Lambda_n$, we have
\[\Phi_n^\varepsilon(x)\leq\exp\left(-\alpha_3\left(\mathcal{T}_n^{p_n'(\varepsilon)}\right)\cdot x^2+\beta_3\left(\mathcal{T}_n^{p_n'(\varepsilon)}\right)\cdot|x|+\theta_3\left(\mathcal{T}_n^{p_n'(\varepsilon)}\right)\right),\]
where $\alpha_3(\cdot)$, $\beta_3(\cdot)$ and $\theta_3(\cdot)$ are given by \eqref{eq:alpha3}, \eqref{eq:beta3} and \eqref{eq:theta3}.
Furthermore, by Remark \ref{rem:lemtech2}, since the technical condition \eqref{eq:nandA3} is satisfied, we have
\[\alpha_3\left(\mathcal{T}_n^{p_n'(\varepsilon)}\right)\geq\frac{cn}{4(n-K_n^\varepsilon)}.\]
It follows that
\begin{eqnarray*}
\int_{-\infty}^\infty\Phi_n^\varepsilon(x)\mathrm{d}x&\geq&\int_{-\Lambda_n}^{\Lambda_n}\exp\left(-\alpha_1\left(\mathcal{T}_n^{p_n'(\varepsilon)}\right)\cdot x^2-\beta_1\left(\mathcal{T}_n^{p_n'(\varepsilon)}\right)\cdot|x|-\theta_1\left(\mathcal{T}_n^{p_n'(\varepsilon)},\Lambda_n\right)\right)\mathrm{d}x\\
&=:&I_n^1,
\end{eqnarray*}
and
\begin{eqnarray*}
\lefteqn{\int_{-\infty}^\infty x^2\cdot\Phi_n^\varepsilon(x)\mathrm{d}x}\\
&\leq&\int_{-\Lambda_n}^{\Lambda_n}x^2\cdot\exp\left(-\frac{3n}{8\sigma^2(n-K_n^\varepsilon)}\cdot x^2+\beta_2\left(\mathcal{T}_n^{p_n'(\varepsilon)}\right)\cdot|x|+\theta_2\left(\mathcal{T}_n^{p_n'(\varepsilon)},\Lambda_n\right)\right)\mathrm{d}x\\
&&+\int_{|x|\geq\Lambda_n}x^2\cdot\exp\left(-\frac{cn}{4(n-K_n^\varepsilon)}\cdot x^2+\beta_3\left(\mathcal{T}_n^{p_n'(\varepsilon)}\right)\cdot|x|+\theta_3\left(\mathcal{T}_n^{p_n'(\varepsilon)},\Lambda_n\right)\right)\mathrm{d}x\\
&=:&I_n^2+I_n^3.
\end{eqnarray*}
Altogether, we get that on the good event $F_n^\varepsilon$, we have
\[V_n^\varepsilon\leq2\left(\frac{n}{N_n^\varepsilon}\cdot\frac{I_n^2+I_n^3}{I_n^1}+d^2\right)=:\Upsilon_n.\]
Finally, we claim that the sequence $(\Upsilon_n)_{n\in\llbracket2,\infty\llbracket}$ is tight.
Indeed:
\begin{itemize}
\item By \eqref{eq:fact1'}, we have
\[\alpha_1\left(\mathcal{T}_n^{p_n'(\varepsilon)}\right)\longrightarrow\frac{1}{2\sigma^2}\quad\text{in probability as $n\to\infty$,}\]
and by \eqref{eq:fact1'} and \eqref{eq:fact2'}, the sequences
\[\left(\beta_1\left(\mathcal{T}_n^{p_n'(\varepsilon)}\right)\right)_{n\in\mathbb{N}^*}\quad\text{and}\quad\left(\theta_1\left(\mathcal{T}_n^{p_n'(\varepsilon)}\right)\right)_{n\in\mathbb{N}^*}\]
are tight (in fact, converge in distribution).
This entails that $\varlimsup_{n\to\infty}\mathbb{P}\left(I_n^1\leq\iota\right)\rightarrow0$ as $\iota\to0^+$.

\item By \eqref{eq:fact1'}, we have
\[\frac{3n}{8\sigma^2(n-K_n^\varepsilon)}\longrightarrow\frac{3}{8\sigma^2}\quad\text{in probability as $n\to\infty$,}\]
and by \eqref{eq:fact1'}, \eqref{eq:fact2'} and \eqref{eq:fact3'}, the sequences
\[\left(\beta_2\left(\mathcal{T}_n^{p_n'(\varepsilon)}\right)\right)_{n\in\mathbb{N}^*}\quad\text{and}\quad\left(\theta_2\left(\mathcal{T}_n^{p_n'(\varepsilon)}\right)\right)_{n\in\mathbb{N}^*}\]
are tight.
This entails that ${\varlimsup_{n\to\infty}\mathbb{P}\left(I_n^2\geq A\right)\rightarrow0}$ as $A\to\infty$.

\item By \eqref{eq:fact1'}, we have
\[\frac{cn}{4(n-K_n^\varepsilon)}\longrightarrow\frac{c}{4}\quad\text{in probability as $n\to\infty$,}\]
and by \eqref{eq:fact1'}, \eqref{eq:fact2'} and \eqref{eq:fact3'}, the sequences
\[\left(\beta_3\left(\mathcal{T}_n^{p_n'(\varepsilon)}\right)\right)_{n\in\mathbb{N}^*}\quad\text{and}\quad\left(\theta_3\left(\mathcal{T}_n^{p_n'(\varepsilon)}\right)\right)_{n\in\mathbb{N}^*}\]
are tight.
Since $\Lambda_n=\ln(n)\rightarrow\infty$ as $n\to\infty$, this entails that $I_n^3\rightarrow0$ in probability as $n\to\infty$.
\end{itemize}
The previous points together with \eqref{eq:fact2'} entail that $\varlimsup_{n\to\infty}\mathbb{P}(\Upsilon_n\geq A)\rightarrow0$ as $A\to\infty$.
In particular, there exists $A\in\mathbb{R}_+^*$ such that $\varliminf_{n\to\infty}\mathbb{P}(\Upsilon_n\leq A)\geq1-\zeta$.
Since ${\varliminf_{n\to\infty}\mathbb{P}\left(F_n^\varepsilon\right)\geq1-\zeta}$, it follows that
\[\varliminf_{n\to\infty}\mathbb{P}\left(V_n^\varepsilon\leq A\right)\geq\varliminf_{n\to\infty}\mathbb{P}\left(\Upsilon_n\leq A\,;\,F_n^\varepsilon\right)\geq1-2\zeta.\]
This proves \eqref{eq:locivrhfgoal}, which completes the proof of (the localisation in) Theorem \ref{thm:main}.
\end{proof}

\appendix

\section{Appendix}

\subsection{The discrete Gaussian free field and electric networks}\label{sec:thomson}

In this appendix, we provide a self-contained proof of Proposition \ref{prop:thomson}, which we recall is a direct consequence of Thomson's principle for electric networks, given that $\mathbb{V}(G_n(\varnothing))$ corresponds to the effective resistance between $\varnothing$ and $\partial T_n$.

\begin{proof}[Proof of Proposition \ref{prop:thomson}]
We keep the notation introduced in the proof of Proposition \ref{prop:green}: recall that $C_n(u)=\mathbb{E}[G_n(\varnothing)G_n(u)]$ for all $u\in T_n$.
Equation \eqref{eq:green} entails that
\[I_n:(u,v)\in\vec{E}(T_n)\longmapsto C_n(u)-C_n(v)\]
is a unit flow from $\varnothing$ to $\partial T_n$, i.e, belongs to $\Phi(T_n)$.
Now, we claim that
\begin{equation}\label{eq:thomson1}
C_n(\varnothing)=\sum_{(u,v)\in\vec{E}(T_n)}I_n(u,v)^2,
\end{equation}
and that for any $\varphi\in\Phi(T_n)$, we have
\begin{equation}\label{eq:thomson2}
\sum_{(u,v)\in\vec{E}(T_n)}I_n(u,v)\cdot(\varphi(u,v)-I_n(u,v))=0.
\end{equation}
These two equations entail that for any $\varphi\in\Phi(T_n)$, we have
\[\sum_{(u,v)\in\vec{E}(T_n)}\varphi(u,v)^2=\sum_{(u,v)\in\vec{E}(T_n)}I_n(u,v)^2+\sum_{(u,v)\in\vec{E}(T_n)}(\varphi(u,v)-I_n(u,v))^2\geq C_n(\varnothing),\]
with equality if and only if $\varphi=I_n$, which yields the result of the proposition.
Now, let us prove \eqref{eq:thomson1} and \eqref{eq:thomson2}.
In fact, both equations follow from the fact that for any flow $\psi$ on $T_n$ from $\varnothing$ to $\partial T_n$, i.e, for any function $\psi:\vec{E}(T_n)\rightarrow\mathbb{R}$ such that \eqref{eq:flow} holds (with $\psi$ in place of $\varphi$), we have
\begin{equation}\label{eq:thomson3}
\sum_{(u,v)\in\vec{E}(T_n)}(C_n(u)-C_n(v))\cdot\psi(u,v)=C_n(\varnothing)\cdot\sum_{\substack{v\in T_n\\\text{$v$ child of $\varnothing$}}}\psi(\varnothing,v).
\end{equation}
Indeed, taking $\psi=I_n$ in \eqref{eq:thomson3} yields \eqref{eq:thomson1}, and \eqref{eq:thomson2} is obtained from \eqref{eq:thomson3} by taking $\psi=\varphi-I_n$.
To complete the proof of the proposition, let us prove \eqref{eq:thomson3}.
For any flow $\psi:\vec{E}(T_n)\rightarrow\mathbb{R}$ from $\varnothing$ to $\partial T_n$, we have
\begin{eqnarray*}
\lefteqn{\sum_{(u,v)\in\vec{E}(T_n)}(C_n(u)-C_n(v))\cdot\psi(u,v)}\\
&=&\sum_{u\in T_n\setminus\partial T_n}\sum_{\substack{v\in T_n\\\text{$v$ child of $u$}}}(C_n(u)-C_n(v))\cdot\psi(u,v)\\
&=&\sum_{u\in T_n\setminus\partial T_n}C_n(u)\cdot\sum_{\substack{v\in T_n\\\text{$v$ child of $u$}}}\psi(u,v)-\sum_{u\in T_n\setminus\partial T_n}\sum_{\substack{v\in T_n\\\text{$v$ child of $u$}}}C_n(v)\cdot\psi(u,v)\\
&=&C_n(\varnothing)\cdot\sum_{\substack{v\in T_n\\\text{$v$ child of $\varnothing$}}}\psi(\varnothing,v)+\sum_{u\in T_n^*\setminus\partial T_n}C_n(u)\cdot\psi(\mathrm{parent}(u),u)\\
&&-\sum_{v\in T_n^*}C_n(v)\cdot\psi(\mathrm{parent}(v),v)\\
&=&C_n(\varnothing)\cdot\sum_{\substack{v\in T_n\\\text{$v$ child of $\varnothing$}}}\psi(\varnothing,v)-\sum_{v\in\partial T_n}C_n(v)\cdot\psi(\mathrm{parent}(v),v)\\
&=&C_n(\varnothing)\cdot\sum_{\substack{v\in T_n\\\text{$v$ child of $\varnothing$}}}\psi(\varnothing,v),
\end{eqnarray*}
where $T_n^*=T_n\setminus\{\varnothing\}$ (the third equality holds since $\psi$ is a flow, and the last equality holds since $C_n(v)=0$ for all $v\in\partial T_n$).
\end{proof}

\subsection{Local central limit estimates}\label{sec:lclt}

In this appendix, we give a proof of the local central limit estimates of Lemma \ref{lem:lclt}.

\begin{proof}[Proof of Lemma \ref{lem:lclt}]
(i) is \cite[Theorem 2.3.11]{lawlerlimic}.\footnote{To be perfectly accurate, \cite[Theorem 2.3.11]{lawlerlimic} is stated for an aperiodic step distribution.
As usual, taking into account the non-aperodicity of our step distribution $\nu$ does not affect the proof, except the estimate becomes valid for $z\in rk+d\mathbb{Z}$.}
Now, let us prove (ii).
Set
\[\delta=\left(\frac{1}{2\sigma\sqrt{C}}\right)\wedge\eta,\]
so that for all $k\in\mathbb{N}^*$, we have
\[C\cdot\frac{z^4}{k^3}\leq\frac{z^2}{4\sigma^2k}\quad\text{for all $z\in\mathbb{R}$ such that $|z|\leq\delta k$.}\]
By the previous point, it follows that for all $z\in rk+d\mathbb{Z}$ such that $|z|\leq\delta k$, we have
\begin{equation}\label{eq:regime1}
\frac{P_k(z)}{q_k(0)}\leq\exp\left(-\frac{z^2}{2\sigma^2k}+\frac{C}{k}+\frac{z^2}{4\sigma^2k}\right)=\exp\left(-\frac{z^2}{4\sigma^2k}+\frac{C}{k}\right).
\end{equation}
On the other hand, by Hoeffding's inequality, we have
\[P_k(z)\leq\exp\left(-\frac{z^2}{2L^2k}\right)\quad\text{for all $z\in rk+d\mathbb{Z}$.}\]
Now, fix $k_0\in\mathbb{N}^*$ large enough so that
\[\frac{1}{q_k(0)}=\frac{\sqrt{2\pi\sigma^2k}}{d}\leq\exp\left(\frac{\delta^2k}{4L^2}\right)\quad\text{for all $k\in\llbracket k_0,\infty\llbracket$.}\]
By the previous two inequalities, for all $k\in\llbracket k_0,\infty\llbracket$, the following holds: for all $z\in rk+d\mathbb{Z}$ such that $|z|\geq\delta k$, we have
\begin{equation}\label{eq:regime2}
\frac{P_k(z)}{q_k(0)}\leq\exp\left(-\frac{z^2}{2L^2k}+\frac{\delta^2k}{4L^2}\right)\leq\exp\left(-\frac{z^2}{2L^2k}+\frac{z^2}{4L^2k}\right)=\exp\left(-\frac{z^2}{4L^2k}\right).
\end{equation}
Combining \eqref{eq:regime1} and \eqref{eq:regime2}, it follows that for all $k\in\llbracket k_0,\infty\llbracket$, we have
\begin{equation}\label{eq:lcltii}
\frac{P_k(z)}{q_k(0)}\leq\exp\left(-\frac{cz^2}{k}+\frac{C}{k}\right)\quad\text{for all $z\in rk+d\mathbb{Z}$,}
\end{equation}
where
\[c=\left(\frac{1}{4\sigma^2}\right)\wedge\left(\frac{1}{4L^2}\right).\]
By adjusting the value of the constant $C$ if necessary, we can ensure that \eqref{eq:lcltii} holds for all $k\in\mathbb{N}^*$ (since $\nu$ has finite support, so does $P_k$ for each $k\in\llbracket1,k_0-1\rrbracket$).
This completes the proof.
\end{proof}

\subsection{The space of Borel probability measures on $\mathbb{R}$}\label{sec:cvM1R}

In this appendix, we discuss convergence of random variables with values in the Polish space $\mathcal{M}_1(\mathbb{R})$ of Borel probability measures on $\mathbb{R}$ endowed with the topology of weak convergence.
By definition, the topology of weak convergence on $\mathcal{M}_1(\mathbb{R})$ is generated by the following basic open subsets:
\[O_{\mu_0,f,r}:=\{\mu\in\mathcal{M}_1(\mathbb{R}):|\mu(f)-\mu_0(f)|<r\},\]
where $\mu_0\in\mathcal{M}_1(\mathbb{R})$, where $f:\mathbb{R}\rightarrow\mathbb{R}$ is continuous and bounded, where $r\in\mathbb{R}_+^*$, and where ${\mu(f)=\int_\mathbb{R}f(x)\mu(\mathrm{d}x)}$.
Since $\mathbb{R}$ is separable and complete, the space $\mathcal{M}_1(\mathbb{R})$ endowed with the topology of weak convergence is Polish: there exists a metric on $\mathcal{M}_1(\mathbb{R})$ (e.g, the Prohorov metric) which induces the topology of weak convergence, and for which $\mathcal{M}_1(\mathbb{R})$ is separable and complete \cite[Theorem 6.8]{billingsley}.
Thus, it makes sense to speak of convergence in probability in that space (and of course, any other metric that induces the topology of weak convergence has the same convergent sequences).
In practice, one uses the following characterisation:

\begin{prop}\label{prop:cvM1R}
Let $(M_n)_{n\in\mathbb{N}^*}$ be a sequence of random variables with values in $\mathcal{M}_1(\mathbb{R})$, and let $\mu_\infty\in\mathcal{M}_1(\mathbb{R})$.
The following assertions are equivalent:
\begin{enumerate}[label=(\roman*)]
\item $M_n\rightarrow\mu_\infty$ in probability as $n\to\infty$,
\item For every bounded continuous function $f:\mathbb{R}\rightarrow\mathbb{R}$, we have
\[M_n(f)\longrightarrow\mu_\infty(f)\quad\text{in probability as $n\to\infty$,}\]
\item For every $a<b\in\mathbb{R}\setminus\{x\in\mathbb{R}:\mu_\infty(\{x\})>0\}$, we have
\[M_n(]a,b])\longrightarrow\mu_\infty(]a,b])\quad\text{in probability as $n\to\infty$.}\]
\end{enumerate}
\end{prop}

\begin{rem}\label{rem:cvM1R}
\begin{itemize}
\item In particular, the following holds: if for each $a<b\in\mathbb{R}$, we have
\[M_n(]a,b])\longrightarrow\frac{1}{\sqrt{2\pi\sigma^2}}\cdot\int_a^b\exp\left(-\frac{x^2}{2\sigma^2}\right)\mathrm{d}x\quad\text{in probability as $n\to\infty$,}\]
then $M_n\rightarrow\mathcal{N}\left(0,\sigma^2\right)$ in probability as $n\to\infty$.

\item Moreover, if $V_n$ is a sequence of positive random variables such that $V_n\rightarrow1$ in probability as $n\to\infty$, then $M_n=\mathcal{N}(0,V_n)\rightarrow\mathcal{N}(0,1)$ in probability as $n\to\infty$.
Indeed, for each $a<b\in\mathbb{R}$, we have
\begin{eqnarray*}
M_n(]a,b])&=&\frac{1}{\sqrt{2\pi V_n}}\cdot\int_a^b\exp\left(-\frac{x^2}{2V_n}\right)\mathrm{d}x\\
&=&\frac{1}{\sqrt{2\pi}}\cdot\int_{a\sqrt{V_n}}^{b\sqrt{V_n}}\exp\left(-\frac{y^2}{2}\right)\mathrm{d}y\\
&\longrightarrow&\frac{1}{\sqrt{2\pi}}\cdot\int_a^b\exp\left(-\frac{y^2}{2}\right)\mathrm{d}y\quad\text{in probability as $n\to\infty$.}
\end{eqnarray*}
\end{itemize}
\end{rem}

\begin{proof}[Proof of Proposition \ref{prop:cvM1R}]
Assume that (i) holds, and let us prove (ii).
Fix a bounded continuous function $f:\mathbb{R}\rightarrow\mathbb{R}$.
Since the mapping $\mu\in\mathcal{M}_1(\mathbb{R})\mapsto\mu(f)$ is continuous by the definition of the topology of weak convergence, the continuous mapping theorem entails that $M_n(f)\rightarrow\mu_\infty(f)$ in probability as $n\to\infty$.

\quad

Now, assume that (ii) holds, and let us prove (iii).
Fix $a<b\in\mathbb{R}\setminus\{x\in\mathbb{R}:\mu_\infty(\{x\})>0\}$.
For each $\delta\in\mathbb{R}_+^*$, consider the bounded continuous functions $f_\delta,g_\delta:\mathbb{R}\rightarrow\mathbb{R}$ given by
\[f_\delta(x)=\begin{cases}
0&\text{if $x\leq a$,}\\
\frac{x-a}{\delta}&\text{if $x\in[a,a+\delta]$,}\\
1&\text{if $x\in[a+\delta,b-\delta]$,}\\
\frac{b-x}{\delta}&\text{if $x\in[b-\delta,b]$,}\\
0&\text{if $x\geq b$,}
\end{cases}\quad\text{for all $x\in\mathbb{R}$,}\]
and
\[g_\delta(x)=\begin{cases}
0&\text{if $x\leq a-\delta$,}\\
\frac{x-a+\delta}{\delta}&\text{if $x\in[a-\delta,a]$,}\\
1&\text{if $x\in[a,b]$,}\\
\frac{b+\delta-x}{\delta}&\text{if $x\in[b,b+\delta]$,}\\
0&\text{if $x\geq b+\delta$,}
\end{cases}\quad\text{for all $x\in\mathbb{R}$.}\]
We have $0\leq f_\delta(x)\leq\mathbf{1}_{]a,b]}(x)\leq g_\delta(x)\leq1$ for all $x\in\mathbb{R}$, and ${f_\delta(x)\rightarrow\mathbf{1}_{]a,b[}(x)}$ and ${g_\delta(x)\rightarrow\mathbf{1}_{[a,b]}(x)}$ as $\delta\to0$.
By the bounded convergence theorem, we get that
\[\mu_\infty(f_\delta)\longrightarrow\mu_\infty(]a,b[)=\mu_\infty(]a,b])\quad\text{and}\quad\mu_\infty(g_\delta)\longrightarrow\mu_\infty([a,b])=\mu_\infty(]a,b])\]
as $\delta\to0$.
In particular, for each $\varepsilon\in\mathbb{R}_+^*$, there exists $\delta\in\mathbb{R}_+^*$ such that
\[\mu_\infty(f_\delta)\geq\mu_\infty(]a,b])-\varepsilon\quad\text{and}\quad\mu_\infty(g_\delta)\leq\mu_\infty(]a,b])+\varepsilon.\]
Then, since
\[M_n(f_\delta)\longrightarrow\mu_\infty(f_\delta)\quad\text{and}\quad M_n(g_\delta)\longrightarrow\mu_\infty(g_\delta)\quad\text{in probability as $n\to\infty$,}\]
we get that
\begin{eqnarray*}
\lefteqn{\mathbb{P}(|M_n(]a,b])-\mu_\infty(]a,b])|\geq2\varepsilon)}\\
&=&\mathbb{P}(M_n(]a,b])\leq\mu_\infty(]a,b])-2\varepsilon)+\mathbb{P}(M_n(]a,b])\geq\mu_\infty(]a,b])+2\varepsilon)\\
&\leq&\mathbb{P}(M_n(f_\delta)\leq\mu_\infty(f_\delta)-\varepsilon)+\mathbb{P}(M_n(g_\delta)\geq\mu_\infty(g_\delta)+\varepsilon)\\
&\longrightarrow&0\quad\text{as $n\to\infty$,}
\end{eqnarray*}
which is the desired result.

\quad

Next, assume that (iii) holds, and let us prove (ii).
For any function $f:\mathbb{R}\rightarrow\mathbb{R}$ of the form $f=\sum_{i=1}^ky_i\cdot\mathbf{1}_{]a_i,b_i]}$, where $y_i\in\mathbb{R}$ and $a_i<b_i\in\mathbb{R}\setminus\{x\in\mathbb{R}:\mu_\infty(\{x\})>0\}$, we have $M_n(f)\rightarrow\mu_\infty(f)$ in probability as $n\to\infty$.
Then, since any continuous function with compact support can be approximated uniformly by such functions $f$ (recall that $\{x\in\mathbb{R}:\mu_\infty(\{x\})>0\}$ is countable), it easily follows that for any continuous function $g:\mathbb{R}\rightarrow\mathbb{R}$ with compact support, we have $M_n(g)\rightarrow\mu_\infty(g)$ in probability as $n\to\infty$.
Finally, let $h:\mathbb{R}\rightarrow\mathbb{R}$ be a bounded continuous function.
For each $\varepsilon\in{]0,1[}$, there exists $a<b\in\mathbb{R}\setminus\{x\in\mathbb{R}:\mu_\infty(\{x\})>0\}$ such that $\mu_\infty(]a,b])\geq1-\varepsilon$ (recall that $\{x\in\mathbb{R}:\mu_\infty(\{x\})>0\}$ is countable).
Then, consider the function $\varphi:\mathbb{R}\rightarrow\mathbb{R}$ given by
\[\varphi(x)=\begin{cases}
0&\text{if $x\leq a-1$,}\\
x-a+1&\text{if $x\in[a-1,a]$,}\\
1&\text{if $x\in[a,b]$,}\\
b+1-x&\text{if $x\in[b,b+1]$,}\\
0&\text{if $x\geq b+1$,}
\end{cases}\quad\text{for all $x\in\mathbb{R}$.}\]
Notice that $g=h\cdot\varphi$ is continuous and has compact support, and that
\[|h-g|=|h\cdot(1-\varphi)|\leq\|h\|_\infty\cdot\left(1-\mathbf{1}_{]a,b]}\right).\]
By the triangle inequality, we have
\begin{eqnarray*}
|M_n(h)-\mu_\infty(h)|&\leq&|M_n(h)-M_n(g)|+|M_n(g)-\mu_\infty(g)|+|\mu_\infty(g)-\mu_\infty(h)|\\
&\leq&\|h\|_\infty\cdot(1-M_n(]a,b]))+|M_n(g)-\mu_\infty(g)|+\|h\|_\infty\cdot(1-\mu_\infty(]a,b])).
\end{eqnarray*}
Since as $n\to\infty$ the right-hand side converges in probability towards
\[\|h\|_\infty\cdot(1-\mu_\infty(]a,b]))+0+\|h\|_\infty\cdot(1-\mu_\infty(]a,b]))\leq2\|h\|_\infty\cdot\varepsilon,\]
we conclude that
\[\mathbb{P}(|M_n(h)-\mu_\infty(h)|\geq(2\|h\|_\infty+1)\cdot\varepsilon)\longrightarrow0\quad\text{as $n\to\infty$,}\]
which yields the desired result.

\quad

Finally, assume that (ii) holds, and let us prove (i).
By the Portmanteau theorem, it suffices to show that for any open subset $O\subset\mathcal{M}_1(\mathbb{R})$ containing $\mu_\infty$, we have $\mathbb{P}(M_n\in O)\rightarrow1$ as $n\to\infty$.
By the definition of the topology of weak convergence on $\mathcal{M}_1(\mathbb{R})$, there exists a bounded continuous function $f:\mathbb{R}\rightarrow\mathbb{R}$ and $r\in\mathbb{R}_+^*$ such that $\mu_\infty\in O_{\mu_\infty,f,r}\subset O$.
It follows that
\[\mathbb{P}(M_n\notin O)\leq\mathbb{P}(|M_n(f)-\mu_\infty(f)|\geq r)\longrightarrow0\quad\text{as $n\to\infty$,}\]
as desired.
The proof of the proposition is now complete.
\end{proof}

\subsection{Monotonicity for the absolute value of integer-valued random height functions}\label{sec:monotonicity}

In this appendix, we give a proof of the monotonicity property of Proposition \ref{prop:monotonicity} for the absolute value of $\nu$-integer-valued random height functions, where $\nu$ is the uniform distribution on $\{-1,1\}$.

\begin{proof}[Proof of Proposition \ref{prop:monotonicity}]
By the Holley argument, it suffices to check that for any functions ${f=(f_u)_{u\in T_n}\in\mathbb{N}^{T_n}}$ and ${f'=(f'_u)_{u\in T_n}\in\mathbb{N}^{T_n}}$ such that
\begin{equation}\label{eq:admissible}
\mathbb{P}(|H_n|=f)>0\quad\text{and}\quad\mathbb{P}(|H_n'|=f')>0
\end{equation}
and
\begin{equation}\label{eq:order}
f_u\leq f_u'\quad\text{for all $u\in T_n$,}
\end{equation}
and for any node $u\in T_n$, the probability measure
\[\pi:=\mathbb{P}(|H_n(u)|\in\cdot\,|\,\text{$|H_n(v)|=f_v$ for all $v\in T_n\setminus\{u\}$})\]
is stochastically dominated by
\[\pi':=\mathbb{P}(|H_n'(u)|\in\cdot\,|\,\text{$|H_n'(v)|=f'_v$ for all $v\in T_n\setminus\{u\}$}).\]
If $u\in B_n$, then for any coupling $(\xi,\xi')$ of $\pi$ and $\pi'$, we have $\xi=0\leq\xi'$ almost surely, hence $\pi$ is stochastically dominated by $\pi'$.
Next, let us assume that $u\notin B_n$, hence $u\notin B_n'$.
\begin{itemize}
\item If there exists two neighbours $v_1,v_2\in T_n$ of $u$ such that $f_{v_1}<f_{v_2}$, then by \eqref{eq:admissible}, we must have $f_{v_2}=f_{v_1}+2$, and by \eqref{eq:order}, we get that $f'_{v_2}\geq f_{v_2}=f_{v_1}+2$.
It follows that for any coupling $(\xi,\xi')$ of $\pi$ and $\pi'$, we have $\xi=f_{v_1}+1\leq\xi'$ almost surely, hence $\pi$ is stochastically dominated by $\pi'$.

Similarly, if there exists two neighbours $v_1,v_2\in T_n$ of $u$ such that $f'_{v_1}<f'_{v_2}$, then by \eqref{eq:admissible}, we must have $f'_{v_2}=f'_{v_1}+2$, and by \eqref{eq:order}, we get that $f_{v_1}\leq f'_{v_1}$.
It follows that for any coupling $(\xi,\xi')$ of $\pi$ and $\pi'$, we have $\xi\leq f'_{v_1}+1=\xi'$ almost surely, hence $\pi'$ is stochastically dominated by $\pi$.

\item In the remaining cases, by \eqref{eq:admissible} and by \eqref{eq:order}, there must exist $y\leq y'\in\mathbb{N}$ such that $f_v=y$ and $f'_v=y'$ for all neighbours $v\in T_n$ of $u$.
\begin{itemize}
\item If $y<y'$, then since $y$ and $y'$ must have the same parity, we get that $y'-y\geq2$.
Then, for any coupling $(\xi,\xi')$ of $\pi$ and $\pi'$, we get that $\xi\leq y+1\leq y'-1\leq\xi'$ almost surely, hence $\pi$ is stochastically dominated by $\pi'$.

\item In the remaining cases, we must have $y=y'$.
\begin{itemize}
\item If $y=0$, then for any coupling $(\xi,\xi')$ of $\pi$ and $\pi'$, we have $\xi=1=\xi'$ almost surely, hence $\pi$ is stochastically dominated by $\pi'$.

\item For the last two cases, let us observe that
\[\mathbb{P}(|H_n|=f)=\frac{\#\left\{h\in\mathbb{H}^0(T_n,B_n):|h|=f\right\}}{\#\mathbb{H}^0(T_n,B_n)}=\frac{2^{\kappa(f)}}{\#\mathbb{H}^0(T_n,B_n)},\]
where $\kappa(f)$ denotes the number of connected components of $\{v\in T_n:f_v>0\}$.
Now, let $f^\pm=\left(f^\pm_v\right)_{v\in T_n}\in\mathbb{N}^{T_n}$ be the function given by
\[f^\pm_v=\begin{cases}
y\pm1&\text{if $v=u$,}\\
f_v&\text{otherwise,}
\end{cases}\quad\text{for all $v\in T_n$.}\]
The distribution of $\xi$ is given by
\begin{equation}\label{eq:distxi}
\pi(\{y\pm1\})=\frac{\mathbb{P}\left(|H_n|=f^\pm\right)}{\mathbb{P}\left(|H_n|=f^+\right)+\mathbb{P}\left(|H_n|=f^-\right)}=\frac{2^{\kappa\left(f^\pm\right)}}{2^{\kappa\left(f^+\right)}+2^{\kappa\left(f^-\right)}}.
\end{equation}
We are now ready to treat the last two cases.

\item If $y\geq2$, then we have $\kappa\left(f^+\right)=\kappa\left(f^-\right)$.
By \eqref{eq:distxi}, it follows that $\pi$ is the uniform distribution on $\{y-1,y+1\}$.
The same holds for $\pi'$, hence $\pi$ is stochastically dominated by $\pi'$.

\item Finally, if $y=1$, then we have $\kappa\left(f^-\right)-\kappa\left(f^+\right)=\deg(u)-1$.
By \eqref{eq:distxi}, it follows that
\[\pi(\{2\})=\frac{1}{1+2^{\deg(u)-1}}\quad\text{and}\quad\pi(\{0\})=\frac{2^{\deg(u)-1}}{1+2^{\deg(u)-1}}.\]
The same holds for $\pi'$, hence $\pi$ is stochastically dominated by $\pi'$.
\end{itemize}
\end{itemize}
\end{itemize}

The Holley argument then entails that for any non-decreasing function $\Phi:\mathbb{N}^{T_n}\rightarrow\mathbb{R}_+$, we have $\mathbb{E}[\Phi(|H_n|)]\leq\mathbb{E}[\Phi(|H_n'|)]$.
In particular, we get that
\[\mathbb{V}(H_n(\varnothing))=\mathbb{E}\left[H_n(\varnothing)^2\right]\leq\mathbb{E}\left[H_n'(\varnothing)^2\right]=\mathbb{V}(H_n'(\varnothing)).\]
\end{proof}

\bibliographystyle{siam}
\bibliography{biblio.bib}

@book{concentration,
author = {Boucheron, Stéphane and Lugosi, G\'{a}bor and Massart, Pascal},
title = {{Concentration Inequalities: A Nonasymptotic Theory of Independence}},
publisher = {Oxford University Press},
year = {2013}}

@article{irenequentinising,
author = {Ayuso Ventura, Irene and Berger, Quentin},
title = {{Ising model on a Galton--Watson tree with a sparse random external field}},
journal = {Electronic Journal of Probability},
volume = {29},
year = {2024},
pages = {1--43}}

@article{irenequentinnonlinear,
author = {Ayuso Ventura, Irene and Berger, Quentin},
title = {{Non-linear conductances of Galton--Watson trees and applications to the (near) critical random cluster model}},
journal = {The Annals of Applied Probability},
year = {2025+}}

@article{lammerstoninelli,
author = {Lammers, Piet and Toninelli, Fabio},
title = {{Height function localisation on trees}},
journal = {Combinatorics, Probability and Computing},
volume = {33},
year = {2024},
pages = {50--64}}

@article{peledsamotijyehudayoff,
author = {Peled, Ron and Samotij, Wojciech and Yehudayoff, Amir},
title = {{Grounded Lipschitz functions on trees are typically flat}},
journal = {Electronic Communications in Probability},
volume = {18},
year = {2013},
pages = {1--9}}

@incollection{wernerpowellgff,
author = {Werner, Wendelin and Powell, Ellen},
title = {{Lecture Notes on the Gaussian Free Field}},
booktitle = {Cours spécialisés},
volume = {28},
year = {2021}}

@book{lawlerlimic,
author = {Lawler, Gregory F. and Limic, Vlada},
title = {{Random Walk: A Modern Introduction}},
series = {Cambridge Studies in Advanced Mathematics},
volume = {128},
publisher = {Cambridge University Press},  
year = {2010}}

@article{localconv,
author = {Butler, Nathaniel and Krishnan, Kesav and Ray, Gourab and Spinka, Yinon},
title = {{On the local convergence of integer-valued Lipschitz functions on regular trees}},
journal = {arXiv preprint},
volume = {2410.05542},
year = {2024}}

@article{chenhulin,
author = {Chen, Dayue and Hu, Yueyun and Lin, Shen},
title = {{Resistance growth of branching random networks}},
journal = {Electronic Journal of Probability},
volume = {23},
year = {2018},
pages = {1--17}}

@article{pietseb,
author = {Lammers, Piet and Ott, Sébastien},
title = {{Delocalisation and absolute-value-FKG in the solid-on-solid model}},
journal = {Probability Theory and Related Fields},
volume = {188},
year = {2024},
pages = {63--87}}

@book{lyonsperes,
author = {Lyons, Russell and Peres, Yuval},
title = {{Probability on Trees and Networks}},
publisher = {Cambridge University Press},
year = {2017},
series = {Cambridge Series in Statistical and Probabilistic Mathematics}}

@article{logvariance,
author = {Duminil-Copin, Hugo and Harel, Matan and Laslier, Benoit and Raoufi, Aran and Ray, Gourab},
title = {{Logarithmic Variance for the Height Function of Square-Ice}},
journal = {Communications in Mathematical Physics},
volume = {396},
year = {2022},
pages = {867--902}}

@book{legall,
author = {Le Gall, Jean-François},
title = {{Measure Theory, Probability, and Stochastic Processes}},
publisher = {Springer},
year = {2022},
series = {Graduate Texts in Mathematics}}

@article{hellerspinka,
author = {Heller, Alon and Spinka, Yinon},
title = {{Random homomorphisms and Lipschitz functions on trees}},
journal = {arXiv preprint},
volume = {2606.29426},
year = {2026}}

@book{billingsley,
author = {Billinglsey, Patrick},
title = {{Convergence of Probability Measures}},
publisher = {Wiley},
year = {1999},
series = {Wiley Series in Probability and Statistics}}

@article{wu,
author = {Wu, Wei},
title = {{A Central Limit Theorem for Square Ice}},
journal = {Acta Mathematica Sinica, English Series},
volume = {42},
year = {2026},
pages = {1609--1620}}

@article{coquilledarioleny,
author = {Coquille, Loren and Dario, Paul and Le Ny, Arnaud},
title = {{Quantitative delocalisation for the Gaussian and $q$-SOS long-range chain}},
journal = {Annales Henri Lebesgue},
volume = {9},
year = {2026},
pages = {1379--1432}}

@article{chalhoubdariofaipeurleny,
author = {Chalhoub, Christopher and Dario, Paul and Faipeur, Corentin and Le Ny, Arnaud},
title = {{Delocalisation and scaling limit for the disordered long-range Discrete Gaussian Chain}},
journal = {arXiv preprint},
volume = {2609.18786},
year = {2026}}

@article{randomgraphhom,
author = {Benjamini, Itai and Häggström, Olle and Mossel, Elchanan},
title = {{On Random Graph Homomorphisms into $\mathbb{Z}$}},
journal = {Journal of Combinatorial Theory, Series B},
volume = {78},
year = {2000},
pages = {86--114}}

\end{document}